\documentclass[12 pt]{amsart}
\usepackage{amssymb,amscd,amsmath,latexsym}
\usepackage[mathscr]{euscript}
\usepackage[all]{xy}
\usepackage{layout}
\usepackage{pb-diagram}
\usepackage{pstricks}

\def\A{\mathcal{A}}

\def\C{\mathcal{C}}
\def\D{\mathcal{D}}
\def\Do{\mathbb{D}}

\def\E{\mathcal{E}}
\def\M{\mathcal{M}}
\def\Mi{\textbf{M}}
\def\Fi{\textbf{F}}
\def\F{\mathcal{F}}

\def\O{\mathcal{O}}
\def\Os{\mathbb{O}}

\def\Zr{\mathcal{Z}}

\def\T{\mathbb{T}}

\def\P{\mathbb{P}}
\def\Per{\mathbb{P}}

\def\FS{\mathcal{FS}}
\def\CRS{\Sigma\mathcal{CRS}}
\def\Sym{\Sigma \mathcal{S}}
\def\W{\mathcal{W}}
\def\Po{\mathbb{P}^{{\E}}}
\def\Mo{\mathcal{M}^{{\E}}}
\def\V{\bigvee}
\def\x{\times}
\def\ot{\otimes}
\def\oti{\displaystyle\mathop{\otimes}}
\def\cat{\textbf{Cat}}
\def\op{\oplus}

\def\u{\underline{u}}

\def\m{\underline{m}}
\def\n{\underline{n}}
\def\ux{\underline{x}}
\def\y{\underline{y}}
\def\z{\underline{z}}
\def\c{\underline{c}}
\def\e{\underline{e}}
\def\d{\underline{d}}
\def\v{\underline{v}}
\def\f{\underline{f}}
\def\g{\underline{g}}

\def\nin{\newline\indent}
\def\Proof{\medskip\noindent{\bf Proof: }}
\newtheorem{theorem}{Theorem}
\newtheorem{proposition}[theorem]{Proposition}
\newtheorem{lemma}[theorem]{Lemma}

\newtheorem{corollary}[theorem]{Corollary}
\newtheorem{definition}[theorem]{Definition}
\newtheorem{remark}[theorem]{Remark}

\begin{document}
\title[From fibered symmetric bimonoidal categories to symmetric spectra]{From fibered symmetric bimonoidal categories to symmetric spectra}
\author[Jos\'e Manuel G\'omez]{Jos\'e Manuel G\'omez}
\address{Department of Mathematics, University of Michigan, Ann Arbor, MI 48109, USA}
\email{josmago@math.ubc.ca,josmago@umich.edu}
\thanks{The author was supported in part by NSF RTG Grant \# 0602191}

\begin{abstract}
In here we define the concept of fibered symmetric bimonoidal categories. These are roughly speaking fibered categories $\Lambda:\D\to \C$ whose fibers are symmetric monoidal categories parametrized by $\C$ and such that both $\D$ and $\C$ have a further structure of a symmetric monoidal category that satisfy certain coherences that we describe. Our goal is to show that we can correspond to a fibered symmetric bimonoidal category an $E_{\infty}$-ring spectrum in a functorial way. 
\end{abstract}

\maketitle

\section{Introduction}

The goal of this paper is to study in general what we denote by fibered symmetric bimonoidal categories. These are fibered categories $\Lambda:\D\to \C$, such that $\D$ and $\C$ are topological categories and such that each fiber $\D_{c}$ has the structure of a symmetric monoidal category $(\D_{c},\op_{c},0_{c})$. In addition, both $\D$ and $\C$ are symmetric monoidal categories, $(\D,\ot,1)$ and $(\C,\ot 1)$, and the functor $\Lambda$ is a continuous symmetric monoidal functor. The operations $\ot$ and $\op_{c}$, are compatible in the sense that they satisfy some coherences. These coherences are similar to those satisfied by a symmetric bimonoidal category plus some new coherences that we require. (For the precise definition see Definition \ref{fiberedsymmetricbimonoidal} below). We show that given a fibered symmetric bimonoidal category, we can correspond an $E_{\infty}$-ring spectrum  in a functorial way. 

Our motivation to study fibered symmetric bimonoidal categories comes from an attempt to produce a geometric based model for elliptic cohomology. In \cite{HuKriz}, Hu and Kriz proposed a construction of an elliptic cohomology type spectrum based on a concept that they referred to as elliptic bundles. These are defined to be stringy bundles on a fixed elliptic curve that are invariant under translations. Roughly speaking, a stringy bundle is a holomorphic analogue of Segal's elliptic objects (see \cite{HuKriz} for definitions). The category of elliptic bundles and isomorphisms between them is a symmetric monoidal category under a product $\ot$. One difficulty in the construction of \cite{HuKriz} is, that as defined, there is no obvious way to add elliptic bundles.  Hu and Kriz get around this inconvenience by using some machinery from homotopy theory. In \cite{GomezHuKriz}, G\'omez, Hu and Kriz propose a geometric model of elliptic cohomology using a modified version of stringy bundles over $\mathbb{C}$ with compact supports. This produces an example of a fibered symmetric bimonoidal category which is our motivating example. More precisely, our main motivation for studying fibered symmetric bimonoidal categories is the following construction. Let $\C$ be the SPCMC of worldsheets over the stack of finite dimensional complex manifolds. Given a Riemann surface $X$ we can associate a SPCMC $\C_{X}$ (see \cite{GomezHuKriz} for the definition). We have the topological category $\textbf{A}$ of conformal field theories (CFT's) on $\C$ which is a fibered symmetric bimonoidal category over the discrete category $\mathcal{MF}$ of modular functors over $\C$. In addition, we have the topological category $\textbf{B}$ of CFT's $\C_{\mathbb{C}}$ over modular functors pullback from $\C$. This is also fibered symmetric bimonoial category. Using the machinery that we describe here we can correspond to this data, two $E_{\infty}$-ring spectra $\E_{\textbf{A}}$ and  $\E_{\textbf{B}}$, together with a map of spectra $\phi:\E_{\textbf{A}}\to \E_{\textbf{B}}$. The fiber of this map is our proposed elliptic cohomology type spectrum. (See \cite{GomezHuKriz} for more details).

We begin our study of fibered symmetric bimonoidal categories by first studying the discrete case; that is, we first study the case where all the structure insight is discrete. By applying a streefication process we show that in general a discrete fibered symmetric bimonoidal category can be replaced by an equivalent fibered category $\Lambda:\D^{s}\to \C^{s}$ such that $(\D^{s},\ot,1)$ and $(\C^{s},\ot,1)$ are permutative categories and $\Lambda^{s}:\D^{s}\to \C^{s}$ is a strict map. In addition for every object $c$ of $\C^{s}$ the fibers $\D^{s}_{c}$ have the structure of a permutative category $(\D^{s}_{c},\op_{c},0_{c})$. We also have distributivity maps satisfying coherences similar to those satisfied by a bipermutative category. Thus $\Lambda^{s}:\D^{s}\to \C^{s}$ is what we call a discrete fibered bipermutative category. (See Definition \ref{bipermutativeonfibers} below).

In general, having a fibered category $\Lambda:\D\to \C$ is equivalent to having a contravariant lax 2-functor $\C^{op}\to \cat$ and the latter can be replaced by an equivalent functor $\C^{op}\to \cat$. This is shown for example in \cite[Sections 3.1.2 and 3.1.3]{Vistolo}. We use this idea as our starting point. Thus we can correspond to a fibered bipermutative category $\Lambda:\D\to \C$, a functor $\C^{op}\to \P$, where $\P$ is the category of small permutative categories. However, in order to capture the multiplicative structure present we need to enlarge the category $\C^{op}$. Thus given a fibered bipermutative category $\Lambda:\D\to \C$ we consider the permutative category $\A$ that is defined as the wreath product category $\textbf{Inj}\int (\C^{op})^{*}$, for a functor 
\[
(\C^{op})^{*}:\textbf{Inj}\to \P.
\] 
We show that given a fibered bipermutative category $\Lambda:\D\to \C$, we can correspond a functor $\Psi:\A\to \P$. This functor preserves the multiplicative structure in the sense that for objects $u$ and $v$ of $\A$, we have a functor 
\[
\ot_{\u,\v}:\Psi(\u)\x\Psi(\v)\to \Psi(\u\odot\v)
\]
that satisfies certain conditions. (These are conditions $(\textbf{c.1})-(\textbf{c.14})$ of Theorem \ref{multiEsigma}).

The functor $\Psi:\A\to \P$ induced by a fibered bipermutative category motivates the study in general of the category of functors $\E\to \P$, where $(\E,\ot,1)$ is a general permutative category. We define a multicategory structure on the this category. This is a special case of a more general construction; that is, if $\M$ is a multicategory, then we show that we can find a multicategory $\M^{\E}$ whose objects are the functors $\E\to \M$. In addition, we show that for the functor $\Psi:\A\to \P$ induced by a fibered bipermutative category the additional multiplicative structure given by the permutative structure on $\D$ and the coherences gives rise to an enriched multifunctor over $\cat$ 
\[
T_{1}:E\Sigma_{*}\to \P^{\A}.
\] 
Here $E\Sigma_{*}$ is the category valued operad, whose value at $k\ge 0$, is the translation category $E\Sigma_{k}$ that has as object set $\Sigma_{k}$ and there is only one morphism between any two objects of $E\Sigma_{k}$.

We also show that given an enriched multifunctor $\F:\M_{1}\to \M_{2}$ between two multicategories $\M_{1}$ and $\M_{2}$ enriched over $\cat$, then composition with $\F$ gives rise to a multifunctor $\F_{*}:\M_{1}^{\E}\to\M_{2}^{\E}$. In particular, if we take the multifunctor $K:\P\to \Sym$ defined by Elemendorf and Mandell \cite[Theorem 1.1]{Elmendorf}, we get a multifunctor $K_{*}:\Po\to \Sym^{\E}$. By composing the multifunctors $T_{1}$ and $K_{*}$ in the case of $\E=\A$, we obtain a multifunctor $T:E\Sigma_{*}\to \Sym^{\A}$. 
In addition, we show that in general, given a multicategory $\M$ enriched over simplicial sets, then the positive model structure on $\Sym$ can be lifted as to get a closed model category structure on the category of multifunctors $\M\to \Sym^{\E}$. We use this model structure to show that a multifunctor $T:E\Sigma_{*}\to \Sym^{\E}$ can be rectified as to obtain a multifunctor $T':*\to \Sym^{\E}$. It turns out that having such a multifunctor is equivalent as having a lax map $\E\to \Sym$. Thus in the case of $\A$, we obtain a lax map $\vartheta:\A\to \Sym$. 

We have a canonical map $\W:\A\to \C^{op}$ that is a strict map. Using the model structure mentioned above, we show that this map induces, via a left adjoint in a Quillen adjunction, a lax map $\phi:\C^{op}\to \Sym$. In a similar way, we show that such a functor gives rise to a lax map $\phi':(\C^{op})^{-1}\C^{op}\to \Sym$. In particular, the image of the unit give us a strictly commutative symmetric ring spectrum. 
The outlined construction is functorial, thus we obtain a functor 
\[
\Zr:\FS\to \CRS, 
\]
where $\FS$ is the category of discrete fibered categories and $\CRS$ is the category of strictly commutative symmetric ring spectra. 

Finally, given a topological fibered category, by applying the singular functor, we can see it as a simplicial fibered category; that is, given a topological fibered category we correspond a functor
\[
\Delta^{op}\to \FS
\]
and by composing this functor with the functor $\Zr$ we obtain a simplicial commutative ring spectrum whose realization is the desired spectrum.

A big part of this work was inspired from the ideas of  Elmendorf and Mandell in their beautiful work \cite{Elmendorf} and \cite{Elmendorf2}. Some propositions and theorems here are direct adaptations of their ideas to our settings. Also I would like to thank Igor Kriz and Po Hu  for all their help and useful comments throughout this work.

\section{From fibered symmetric bimonoidal categories to fibered bipermutative categories}

In this section we define fibered symmetric bimonoidal categories and fibered bipermutative categories. Then we show, using a standard procedure in category theory, that every fibered symmetric bimonoidal category can be replaced by an equivalent fibered bipermutative category. The latter will be used as input for the machine that we construct and that produces an $E_{\infty}$-ring spectrum. 

\begin{definition}\label{fiberedsymmetricbimonoidal} 
A discrete fibered symmetric bimonoidal category is a fibered category $\Lambda:\D\to \C$, (see Definition \ref{fiberedcategory}), where $(\D,\ot,1,\gamma^{\ot})$ and $(\C,\ot,1,\gamma^{\ot})$ are two small symmetric monoidal categories and $\Lambda$ is a symmetric monoidal functor. In addition, for each object $c$ of $\C$, the fiber $\D_{c}$ has the structure of a symmetric monoidal category $(\D_{c},\op_{c},0_{c},\gamma^{\op}_{c})$. These structures are compatible in the sense that we have natural distributivity maps 
\begin{align*}
d^{l}:(x\ot y)\op(x'\ot y)\to(x\op x')\ot y,\\
d^{r}:(x\ot y)\op(x\ot y')\to x\ot(y\op y),
\end{align*}
defined whenever $x$ and $x'$, $y$ and $y'$ are in fibers $\D_{c}$, $\D_{d}$ respectively. The morphisms $d^{l}$ and $d^{r}$ and are morphisms in $\D_{\Lambda(x\ot y)}=\D_{\Lambda((x\op x')\ot y)}$. We also have natural isomorphisms
\begin{align*}
\lambda'_{c}:&0_{c}\op x \to x, \ \ \ \ \ \ \rho'_{c}:x\op 0_{c}\to x_{c}'\\
\lambda^{*}_{c}:&0_{d}\ot x \to 0_{d\ot c},\ \ \rho_{d}^{*}:x\ot 0_{d}\to 0_{c\ot d},
\end{align*}
for an object $x$ in $\D_{c}$. We require similar coherences than those described by Laplaza for symmetric bimonoidal categories as described in \cite{Laplaza} to hold whenever they make sense with the zero and null morphisms replaced by the above morphisms. In addition, we can add morphisms in $\D$ over a same morphism in $\C$; that is, given morphisms $\alpha:x\to y$ and $\beta:x'\to y'$ with $\Lambda(\alpha)=\Lambda(\beta)=f$, for $f:c\to d$, then we can find a morphism 
\[
\alpha\op \beta:x\op x'\to y\op y',
\]
with $\Lambda(\alpha\op \beta)=f$. In the case that $f=\text{id}_{c}$, then $\alpha\op \beta=\alpha\op_{c} \beta$. We require the following diagrams to be commutative:
\[
\xymatrix{
x\op x'\ar[r]^-{\alpha\op \beta}\ar[d]_-{\gamma_{c}^{\op}} &y\op y'\ar[d]^-{\gamma^{\op}_{d}}\\
x'\op x\ar[r]^-{\beta\op \alpha}&y'\op y
}
\]
\[
\xymatrix{
(x\op x')\op x''\ \ \ar[r]^-{(\alpha\op \beta)\op \delta}\ar[d]_-{\cong}&\ \ (y\op y')\op y''\ar[d]^-{\cong}\\
x\op(x'\op x'')\ \ \ar[r]^-{\alpha\op (\beta\op \delta)}&\ \ y\op (y'\op y'')
}
\]
Also, given any morphism $f:c\to d$ in $\C$ we can find a morphism $0_{f}:0_{c}\to 0_{d}$ over $f$ such that
\[
\xymatrix{
x\op 0_{c}\ar[d]_-{\cong}\ar[r]^-{\alpha\op 0_{f}}&y\op 0_{d}\ar[d]^-{\cong}\\ 
x\ar[r]_-{\alpha}&y.
}
\]
The morphisms $0_{f}$'s are defined in a functorial way; that is, given $f:c\to d$ and $g:d\to e$ morphisms in $\C$, then 
$0_{g\circ f}=0_{g}\circ 0_{f}$ and if $f=\text{id}:c\to c$, then $0_{id_{c}}=\text{id}_{0_{c}}$.
\end{definition}

\begin{remark} The naturality statement for all the morphisms in the above definition means naturality whenever this makes sense, for example for the distributivity morphisms $d^{l}$ and $d^{r}$, the naturality means that given $g:x\to x_{1}$ and $g':x'\to x_{1}'$ morphisms over $f:c\to c'$ and a morphism $h:y\to y'$ over $f':d\to d'$, then the following diagram is commutative
\[
\xymatrix{
(x\ot y)\op (x'\ot y)\ar[r]^-{d^{l}}\ar[d]_-{(g\ot h)\op (g'\ot h)}&(x\op x')\ot y\ar[d]^-{(g\op g')\ot h}\\
(x_{1}\ot y_{1})\op (x_{1}'\ot y_{1})\ar[r]_-{d^{l}}&(x_{1}\op x_{1}')\ot y_{1}
} 
\]
and similarly for $d^{r}$.
\end{remark}

\begin{remark} If $\C$ is the trivial category with only one object and one morphism, then a discrete fibered symmetric bimonoidal category over $\C$ is just a symmetric bimonoidal category.
\end{remark}

In the case of symmetric bimonoidal categories, it is convenient to work in a more rigid scenario by requiring the operations in sight to be strictly associative and to have strict units. Similarly, we have a strict version of a fibered symmetric bimonoidal category. We call these fibered bipermutative categories. Before giving the precise definition, we recall the definition of a permutative category.

\begin{definition} A permutative category consists of a category $\C$, a functor $\ot:\C\x\C\to \C$ that is strictly associative, a strict unit $1$; that is, $x\ot 1=1=1\ot x$ for all objects $x$ and a natural isomorphism $\gamma:x\ot y\to y\ot x$ such that the following diagrams are commutative
\[
\xymatrix{
x\ot 1\ar[rr]^-{\gamma}_\cong\ar[dr]_{=}&&1\ot x\ar[dl]^-{=} & &x\ot y\ar[rr]^-{=}\ar[dr]^\cong_\gamma&&x\ot y\\
&x & & & &y\ot x,\ar[ur]^\cong_\gamma\\
}
\]
\[
\xymatrix{
x\ot y\ot z\ar[rr]^\gamma\ar[dr]_{\textsl{id}\ot\gamma}&&z\ot
x\ot y\\ &x\ot z\ot y\ar[ur]_{\gamma\ot \textsl{id}}.\\
}
\]
\end{definition}

We have three different notions of morphisms between permutative categories that will be useful for us. Suppose that $(\D,\ot,1)$ and $(\E,\ot,1)$ are permutative categories, then: 

\begin{itemize}
\item A functor $f:\C\to \D$ is said to be a strict map if $f(x\ot y)=f(x)\ot f(y)$, $f(1)=1$ and 
\[
\xymatrix{
f(x\ot y)\ar[r]^-{=}\ar[d]_-{f(\gamma)}   &f(x)\ot f(y)\ar[d]^-{\gamma}\\
f(y\ot x)\ar[r]_-{=}        &f(y)\ot f(x)
}
\]
is a commutative diagram. 
\item A functor $f:\C\to \D$ is said to be a lax$_{*}$ map if $f(1)=1$  and there exists a natural transformation 
\[
\lambda=\lambda_{f}:f(x)\ot f(y)\to f(x\ot y) 
\] 
such that $\lambda=\text{id}$ whenever $x=1$ or $y=1$ and the following diagrams are commutative	
\[
\xymatrix{
f(x)\ot f(y)\ot f(z)\ar[r]^-{\text{id}\ot \lambda}\ar[d]_-{\lambda\ot \text{id}}   &f(x)\ot f(y\ot z)\ar[d]^-{\lambda}\\
f(x\ot y)\ot f(z)\ar[r]_-{\lambda}        &f(x\ot y\ot z),
}
\]
\[
\xymatrix{
f(x)\ot f(y)\ar[r]^-{\lambda}\ar[d]_-{\gamma}   &f(x\ot y)\ar[d]^-{f(\gamma)}\\
 f(y)\ot f(x)\ar[r]_-{\lambda}        &f(y\ot x).
}
\]

\item A functor $f:\C\to \D$ is said to be a lax map, if we can find a map $\eta:1\to f(1)$ and a natural transformation 
\[
\lambda=\lambda_{f}:f(x)\ot f(y)\to f(x\ot y) 
\] 
such that the following similar coherences are satisfied together with those coherences involving the unit.
\end{itemize}

From now on we will denote by $\P$ the category whose objects are the small permutative categories and whose morphisms are the lax$_{*}$ maps. This category $\P$ will play a crucial roll throughout this work. 

\begin{definition}\label{bipermutative} A bipermutative category is category that has two permutative category structures $(\C,\op,0,\gamma^{\op})$ and $(\C,\ot,1,\gamma^{\ot})$. These satisfy that $x\ot 0=0=0\ot x$ and we have natural distributivity isomorphisms
\begin{align*}
d^{l}:(x\ot y)\op(x'\ot y)\to(x\op x')\ot y,\\
d^{r}:(x\ot y)\op(x\ot y')\to x\ot(y\op y),
\end{align*}
that satisfy the following coherences:
\[
\tag{\textbf{a.1}}
\xymatrix{
(x\ot y)\oplus(x'\ot y)\oplus(x''\ot
y)\ar[r]^-{d^l\op\textsl{id}}\ar[d]_-{\textsl{id}\op d^{l}}
&((x\op x')\ot y)\op(x''\ot y)\ar[d]^-{d^l}\\
(x\ot y)\op((x'\op x'')\ot y)\ar[r]_-{d^l}
&(x\op x'\op x'')\ot y,
}
\]
\[
\tag{\textbf{a.2}}
\xymatrix{
(x\ot y)\op(x'\ot
y)\ar[r]^-{d^l}\ar[d]_-{\gamma^{\op}}
&(x\op x')\ot y\ar[d]^-{\gamma^{\op}\ot\textsl{id}}\\
(x'\ot y)\op(x\ot y)\ar[r]_-{d^l}&(x'\op x)\ot
y, }
\]
\[
\tag{\textbf{a.3}}
\xymatrix{
(x\ot y\ot z)\op(x'\ot y\ot z)\ar[r]^-{d^l}
\ar[d]_-{d^l}&(x\op x')\ot y\ot z\\
((x\ot y)\op(x'\ot y))\ot z\ar[ur]_-{d^l\ot\textsl{id}}
&
}
\]
are commutative diagrams. Also we have similar commutative diagrams with $d^{l}$ replaced by $d^{r}$.
Also the following diagrams are commutative:
\[
\tag{\textbf{a.4}}
\xymatrix@C-100pt @R-6pt{
&(x\ot(y\oplus y'))\op(x'\ot(y\op y'))\ar[ddr]^-{d^l}\\
(x\ot y)\op(x\ot y')\op(x'\ot y)\op(x'\ot
y')\ar[ur]^-{d^r\op d^r}
\ar[dd]_-{\textsl{id}\op\gamma^{\op}\op\textsl{id}}\\
&&(x\op x')\ot(y\op y'),\\
(x\ot y)\op(x'\ot y)\op(x\ot y')\op(x'\ot
y')\ar[dr]_-{d^l\op
d^l}\\
&((x\op x')\ot y)\op((x\op x')\ot y')\ar[uur]_-{d^r}}
\]
\[
\tag{\textbf{a.5}}
\xymatrix{
(x\ot y)\op(x'\ot y)\ar[r]^-{d^l}\ar[d]_-{\gamma^\ot\op\gamma^\ot}
&(x\op x')\ot
y\ar[d]^-{\gamma^\ot}\\
(y\ot x)\op(y\ot x')\ar[r]_-{d^r}&y\ot(x\op x').\\
}
\]
\end{definition}

\begin{definition}\label{bipermutativeonfibers}
A fibered category $\Lambda:\D\to \C$ is said to be a fibered bipermutative category if it satisfies the following properties: $\Lambda:\D\to \C$ is a fibered category such that $(\D,\ot,1)$ and $(\C,\ot,1)$ are small permutative categories and $\Lambda:\D\to \C$ is a strict map. In addition, we require that the fibers of $\Lambda$ have an additional structure of a permutative category $(\D_{c},\op_{c},0_{c}\gamma^{\op}_{c})$. We assume that $x\ot 0_{d}=0_{c\ot d}$ and $0_{c}\ot y=0_{c\ot d}$ for all objects $x$ and $y$ of $\D_{c}$ and $\D_{d}$ respectively. We also require the operations $\ot$ and $\op_{c}$ to be compatible in the sense that we can find natural distributivity maps which are isomorphisms
\begin{align*}
d^{l}:(x\ot y)\op(x'\ot y)\to(x\op x')\ot y,\\
d^{r}:(x\ot y)\op(x\ot y')\to x\ot(y\op y),
\end{align*}
that are defined for objects $x$ and $x'$, $y$ and $y'$ of $\D_{c}$ and $\D_{c'}$ respectively. The distributive maps are required to be maps in the category $\D_{c\ot c'}$ and also to satisfy the coherences $(\textbf{a.1})-(\textbf{a.5})$ where they make sense. We name these conditions $(\textbf{b.1})-(\textbf{b.5})$.
In addition, given $g:x\to y$ and $g':x'\to y'$ morphisms over $f:c\to c'$, then we require the existence of a morphism $g\op g':x\op x'\to y\op y'$ defined in a natural way; that is, if we have morphisms $h:y\to z$ and $h':y'\to z'$ over $f':c'\to c''$, then the following diagram is commutative
\[
\tag{\textbf{b.6}}
\xymatrix{
x\op x'\ar[d]_-{g\op g'}\ar[rd]^-{(h\circ h)\op (h'\circ g')}&\\
y\op y'\ar[r]_{h\op h'}&z\op z'.
}
\]
In the case that $f=\text{id}:c\to c$, then $g\op g'=g\op_{c} g'$. Also, we require this addition to be associative; that is,
\[
\tag{\textbf{b.7}}
(g\op g')\op g''=g\op(g'\op g'')
\]
for morphism $g,g'$ and $g''$ over $f$.
Moreover, we need the following diagram to commute 
\[
\tag{\textbf{b.8}}
\xymatrix{
x\op x'\ar[r]^-{g\op g'}\ar[d]_-{\gamma^{\op}_{c}}&y\op y'\ar[d]^-{\gamma^{\op}_{c'}}\\
x'\op x\ar[r]^-{g'\op g}&y'\op y.
}
\]
Finally, for each morphism $f:c\to c'$ in $\C$ we can find a morphism $0_{f}:0_{c}\to 0_{c'}$ over $f$, defined in a natural way in the sense that if $f':c'\to c''$ is another morphism in $\C$ then $0_{f'}\circ 0_{f}=0_{f'\circ f}$ and the $0_{f}$ satisfy that 
\[
\tag{\textbf{b.9}}
g\op 0_{f}=g=0_{f}\op g. 
\]
If $f=\text{id}:c\to c$, then $0_{\text{id}}=\text{id}_{0_{c}}$. Note that the naturality of the distributivity morphisms $d^{l}$ and $d^{r}$ means that given $g:x\to x_{1}$ and $g':x'\to x_{1}'$ morphisms over $f:c\to c'$ and a morphism $h:y\to y'$ over $f':d\to d'$, then the following diagram is commutative
\[
\tag{\textbf{b.10}}
\xymatrix{
x\ot y\op x'\ot y\ar[r]^-{d^{l}}\ar[d]_-{g\ot h\op g'\ot h}&(x\op x')\ot y\ar[d]^-{(g\op g')\ot h}\\
x_{1}\ot y_{1}\op x_{1}'\ot y_{1}\ar[r]_-{d^{l}}&(x_{1}\op x_{1}')\ot y_{1}
} 
\]
and similarly for $d^{r}$.
\end{definition}

\begin{remark} A fibered bipermutative category over the trivial category is just a bipermutative category as defined above.
\end{remark}

For symmetric bimonoidal categories, it is well known that every symmetric bimonoidal category is equivalent to a bipermutative category. In our case this is also true; that is, every fibered symmetric bimonoidal category is equivalent to a fibered bipermutative category. We show this in the following theorem.

\begin{theorem}
Given a fibered symmetric bimonoidal category $\Lambda:\D\to \C$, we can find equivalent categories $\D^{s}$ and $\C^{s}$ together with a functor $\Lambda^{s}:\D^{s}\to \C^{s}$ that is a fibered bipermutative category. Moreover, there are equivalences $\Phi':\C\to \C^{s}$ and $\Theta':\D\to \D^{s}$ that are symmetric monoidal functors, and such that the following diagram is commutative
\[
\xymatrix{
\D\ar[r]^-{\Theta'}\ar[d]_-{\Lambda}&\D^{s}\ar[d]^{\Lambda^{s}}\\
\C\ar[r]_-{\Phi'}&\C^{s}.
}
\]
\end{theorem}
\Proof
To begin, let us replace the symmetric monoidal category $\C$ by an equivalent permutative category $\C^{s}$. This is an standard construction in category theory. The objects of $\C^{s}$ are formal products of the form 
\[
\c=c_{1}\boxtimes\cdots \boxtimes c_{n},
\]
where $n\ge 0$. When $n=0$, $\c=()$ is the empty product. Given such a sequence we define 
\[
\Phi(c_{1}\boxtimes\cdots \boxtimes c_{n})=c_{1}\ot(\cdots \ot(c_{n-1}\ot c_{n})\cdots )
\]
for $n>0$, and for $n=0$ we define
\[
\Phi()=1.
\]
Here $1$ is the unit of the symmetric monoidal category $(\C,\ot,1,\gamma^{\ot})$. The morphisms of $\C^{s}$ are defined in the following way. Suppose 
\begin{align*}
\c&=c_{1}\boxtimes\cdots \boxtimes c_{n},\\
\d&=d_{1}\boxtimes\cdots \boxtimes d_{m}
\end{align*}
are two objects of $\C^{s}$. Then a morphism $f:\c\to \d$ in $\C^{s}$ is a morphism in $\C$ $f:\Phi(\c)\to \Phi(\d)$. This way we obtain a category and $\Phi:\C^{s}\to \C$ is an equivalence of categories. The inverse $\Phi':\C\to \C^{s}$ is defined on objects by $\Phi'(c)=c$, where on the right $c$ is the string of objects of $\C$ of length 1. Similarly $\Phi'$ is defined on morphisms. 
\nin
Note that $\C^{s}$ is a permutative category. Indeed, the product $\boxtimes$ is given by juxtaposition; that is, given $\c$ and $\d$ objects of $\C^{s}$ as before, then 
\[
\c\boxtimes\d=c_{1}\boxtimes\cdots \boxtimes c_{n}\boxtimes d_{1}\boxtimes\cdots \boxtimes d_{m}.
\]
Similarly, if $f:\c\to \c'$ and $g:\d\to \d'$ are two morphisms in $\C^{s}$, then we define
\[
f\boxtimes g:\c\boxtimes \d\to \c'\boxtimes\d'
\]
to be the following composite in $\C$
\[
\Phi(\c\boxtimes \d)\stackrel{\cong}{\rightarrow}\Phi(\c)\ot \Phi(\d)\stackrel{f\ot g}{\rightarrow}\Phi(\c')\ot \Phi(\d')\stackrel{\cong}{\rightarrow}\Phi(\c'\boxtimes \d').
\]
Here the outer isomorphisms, are the coherent isomorphisms in $\C$ arising from a rearrangement of parenthesis in a given product. 

Before continuing we introduce some notation. Given a formal sequence 
\[
\ux=x_{1}\boxtimes\cdots \boxtimes x_{k}
\] 
of objects of $\D$, then we denote 
\[
\Delta(\ux)=\Delta(x_{1}\boxtimes\cdots \boxtimes x_{k})=x_{1}\otimes(\cdots\ot(x_{k-1}\ot x_{k})\cdots ).
\]
Thus $\Delta(x_{1}\boxtimes\cdots \boxtimes x_{k})$ is the object of $\D$ obtained by multiplying the elements $x_{1},...,x_{k}$ in a consistent way.
\nin 
On the other hand, since $\Lambda:\D\to \C$ is a fibered category, as explained in Theorem \ref{2lax}, for every morphism $f:c\to c'$ we can correspond a functor $f^{*}:\D_{c'}\to \D_{c}$  in such a way that the correspondence $c\mapsto \D_{c}$, $f\mapsto f^{*}$ is a contravariant lax $2$-functor $\C\to \cat$. We will fix from now on such an assignment.

Now we want to replace the category $\D$ by an equivalent category $\D^{s}$. The objects of $\D^{s}$ will be the set of formal sequences of the form
\[
X=(c_{1}\boxtimes\cdots \boxtimes c_{n},(\ux_{1},f_{1})\boxplus\cdots \boxplus(\ux_{m},f_{m})),
\]
where $n\ge 0$, $m\ge 0$, each $\ux_{i}$ is a formal product
\[
\ux_{i}=x_{i1}\boxtimes\cdots \boxtimes x_{ik_{i}}
\]
and
\[
f_{i}:\Phi(c_{1}\boxtimes\cdots\boxtimes c_{n})\to \Lambda (\Delta(\ux_{i}))
\]
is an isomorphism in $\C$.
\nin
For an object 
\[
X=(c_{1}\boxtimes\cdots \boxtimes c_{n},(\ux_{1},f_{1})\boxplus\cdots \boxplus(\ux_{m},f_{m})),
\]
of $\D^{s}$ we define $\Theta(X)$ in the following way. If $m>0$,
\[
\Theta(X)=f_{1}^{*}(\Delta(\ux_{1}))\op( \cdots \op (f_{m-1}^{*}(\Delta(\ux_{m-1})) \op f_{m}^{*}(\Delta(\ux_{m})))).
\]
When $m=0$ and $n>0$, then $X=(c_{1}\boxtimes\cdots \boxtimes c_{n},)$ and define
\[
\Theta(X)=0_{\Phi(\c)}.
\]
Finally, in $\D^{s}$ we have an object of the form $X=((),1)$  which will be the multiplicative unit and we define
\[
\Theta((),1)=1,
\]
the unit of $\D$.
\nin 
Note that $\Theta(X)$ is a well defined object in $\D$, this is because for every $1\le i\le m$, each $f_{i}^{*}(\Delta(\ux_{i}))$ is an object in  $\D_{\Phi(\c)}$ and addition is well defined on fibers. 
\nin
To define the morphisms of $\D^{s}$ we use $\Theta$ as follows. Take two objects in  $\D^{s}$
\begin{align*}
X&=(c_{1}\boxtimes\cdots \boxtimes c_{n},(\ux_{1},f_{1})\boxplus\cdots \boxplus(\ux_{m},f_{m})),\\
Y&=(d_{1}\boxtimes\cdots \boxtimes d_{r},(\y_{1},g_{1})\boxplus\cdots \boxplus(\y_{s},g_{s})).
\end{align*}
Then a morphism $f:X\to Y$ in $\D^{s}$, is precisely a morphism $f:\Theta(X)\to \Theta(Y)$ in $\D$. In this way we obtain a category $\D^{s}$ and we trivially see that $\Theta:\D^{s}\to \D$ is a functor. It's easy to see that the $\Theta:\D^{s}\to \D$ defines an equivalence of categories with inverse $\Theta':\D\to \D^{s}$ the functor that for an object $x$ of $\D_{c}$ corresponds $\Theta'(x)=(c,(x,\text{id}_{c}))$.

We continue now with the definition of the functor $\Lambda^{s}:\D^{s}\to \C^{s}$. Given an object of $\D^{s}$
\[
X=(c_{1}\boxtimes\cdots \boxtimes c_{n},(\ux_{1},f_{1})\boxplus\cdots \boxplus(\ux_{m},f_{m})),
\]
define
\[
\Lambda^{s}(X)=c_{1}\boxtimes\cdots \boxtimes c_{n}.
\]
If $f:X\to Y$ is a morphism in $\D^{s}$, then $f:\Theta(X)\to \Theta(Y)$ is a morphism in $\D$ and thus $\Lambda(f):\Lambda(\Theta(X))=\Phi(\Lambda^{s}(X))\to \Lambda(\Theta(Y))=\Phi(\Lambda^{s}(Y))$ is a morphism in $\C$. We define then 
\[
\Lambda^{s}(f)=\Lambda(f).
\]

We want to see that $\Lambda^{s}:\D^{s}\to \C^{s}$ satisfies the required properties. To begin, note that by the definition we can easily see that 
\[
\xymatrix{
\D\ar[r]^-{\Theta'}\ar[d]_-{\Lambda}&\D^{s}\ar[d]^{\Lambda^{s}}\\
\C\ar[r]_-{\Phi'}&\C^{s}
}
\]
is a commutative diagram. 

Our next step is to show that each fiber $\D^{s}_{\c}$ has the structure of a permutative category. Suppose then that 
\begin{align*}
X&=(c_{1}\boxtimes\cdots \boxtimes c_{n},(\ux_{1},f_{1})\boxplus\cdots \boxplus(\ux_{m},f_{m})),\\
Y&=(c_{1}\boxtimes\cdots \boxtimes c_{n},(\y_{1},g_{1})\boxplus\cdots \boxplus(\y_{s},g_{s}))
\end{align*}
are two objects in the fiber $\D^{s}_{\c}$.
Define 
\[
X\boxplus_{\c} Y=(c_{1}\boxtimes\cdots \boxtimes c_{n},(\ux_{1},f_{1})\boxplus\cdots \boxplus(\ux_{m},f_{m})\boxplus (\y_{1},g_{1})\boxplus\cdots \boxplus(\y_{s},g_{s})).
\] 
We want to see that this defines a permutative category on $\D^{s}_{\c}$. First of all, it is clear that $\boxplus_{\c}$ is strictly associative. In addition, if we denote by
\[
0_{\c}=(c_{1}\boxtimes\cdots \boxtimes c_{n},)
\]
then it is clear that $0_{\c}$ is strict unit for $\boxplus_{\c}$. We need to find a symmetry isomorphism $\gamma^{\op}_{\c}$. Given objects $X$ and $Y$ in the fiber $\D^{s}_{\c}$ as before, then we define a symmetry isomorphism
\[
\gamma^{\boxplus}_{\c}:X\boxplus_{\c}Y\to Y\boxplus_{\c}X
\]
to be the composite
\[
\Theta(X\boxplus_{\c}Y)\stackrel{\cong}{\rightarrow}\Theta(X)\oplus\Theta(Y)\stackrel{\gamma^{\op}_{c_{1}\ot(\cdots(c_{n-1}\ot c_{n}))}}{\rightarrow}\Theta(Y)\oplus \Theta(X)\stackrel{\cong}{\rightarrow}\Theta(Y\boxplus_{\c}X)
\]
where the outer isomorphism are the coherent isomorphisms in $\D_{\Phi(\c)}$ coming from the associativity of $\op_{\Phi(\c)}$. Then using some coherent theory it follows that each $(\D^{s}_{\c},\boxplus_{\c},0_{\c},\gamma^{\boxplus}_{\c})$ is a permutative category. 

Let us show now that the category $\D^{s}$ has the structure of a permutative category under a product $\boxtimes$. We begin by defining the functor $\boxtimes$. Suppose then that 
\begin{align*}
X&=(c_{1}\boxtimes\cdots \boxtimes c_{n},(\ux_{1},f_{1})\boxplus\cdots \boxplus(\ux_{m},f_{m})),\\
Y&=(d_{1}\boxtimes\cdots \boxtimes d_{r},(\y_{1},g_{1})\boxplus\cdots \boxplus(\y_{s},g_{s}))
\end{align*}
are two objects of $\D^{s}$ with $n, m, r, s>0$. We define  
\begin{align*}
X\boxtimes Y=&\\
(c_{1}\boxtimes\cdots \boxtimes c_{n}\boxtimes d_{1}\boxtimes\cdots \boxtimes d_{r}&, (\ux_{1}\boxtimes \y_{1},f_{1}\boxtimes g_{1})\boxplus\cdots \boxplus (\ux_{1}\boxtimes\y_{s},f_{1}\boxtimes g_{s})\boxplus \\
\cdots\boxplus &(\ux_{m}\boxtimes\y_{1},f_{n}\boxtimes g_{1})\boxplus\cdots \boxplus (\ux_{m}\boxtimes\y_{s},f_{m}\boxtimes g_{s})).
\end{align*}
Here, if
\begin{align*}
\ux_{i}&=x_{i1}\boxtimes\cdots \boxtimes x_{ik_{i}},\\
\y_{j}&=y_{j1}\boxtimes\cdots \boxtimes y_{jt_{j}}
\end{align*}
then
\[
\ux_{i}\boxtimes \y_{j}=x_{i1}\boxtimes\cdots \boxtimes x_{ik_{i}}\boxtimes y_{j1}\boxtimes\cdots \boxtimes y_{jt_{j}},
\]
and $f_{i}\boxtimes g_{j}$ is defined to be the composite in $\D$
\begin{align*}
\Phi(\c\boxtimes \d)\stackrel{\cong}{\rightarrow}&\Phi(\c)\ot \Phi(\d)\stackrel{f_{i}\ot g_{j}}{\rightarrow}\Lambda(\Delta(\ux_{i}))\ot\Lambda (\Delta(\y_{j}))\stackrel{\lambda_{\Lambda}}{\rightarrow}\\
&\Lambda(\Delta(\ux_{i})\ot \Delta(\y_{j}))\stackrel{\cong}{\rightarrow}\Lambda(\Delta(\ux_{i}\boxtimes \y_{j})).
\end{align*}
The outer maps are the ones obtained by a rearrangement of parenthesis. 
Also we we define 
\[
((),1)\boxtimes X=X=X\boxtimes ((),1)
\]
for all objects $X$ of $\D^{s}$. Finally, if $X$ is any object in $\D^{s}_{\c}$, then 
\begin{align*}
0_{\d}\boxtimes X= 0_{\d\boxtimes \c}\\
X\boxtimes 0_{\d}=0_{\c\boxtimes \d}.
\end{align*}
In a similar way we define $\boxtimes$ on morphisms of $\D^{s}$.
\nin
Let us show that this defines a permutative category on $\D^{s}$. We begin by checking the associativity property. Suppose that 
\begin{align*}
X&=(c_{1}\boxtimes\cdots \boxtimes c_{n},(\ux_{1},f_{1})\boxplus\cdots \boxplus(\ux_{m},f_{m})),\\
Y&=(d_{1}\boxtimes\cdots \boxtimes d_{r},(\y_{1},g_{1})\boxplus\cdots \boxplus(\y_{s},g_{s})),\\
Z&=(e_{1}\boxtimes\cdots\boxtimes e_{l},(\z_{1},h_{1})\boxplus\cdots \boxplus(\z_{a},h_{a}))
\end{align*}
are three objects of $\D^{s}$. Then by definition we have that
\begin{align*}
(X\boxtimes Y)\boxtimes Z=(\c\boxtimes \d\boxtimes \e, &(\ux_{1}\boxtimes \y_{1}\boxtimes \z_{1},(f_{1}\boxtimes g_{1})\boxtimes h_{1})\boxplus\\
 &\cdots\boxplus (\ux_{m}\boxtimes \y_{s}\boxtimes \z_{a},(f_{m}\boxtimes g_{s})\boxtimes h_{a})),
\end{align*} 
\begin{align*}
X\boxtimes (Y\boxtimes Z)=(\c\boxtimes \d\boxtimes \e,& (\ux_{1}\boxtimes \y_{1}\boxtimes \z_{1},f_{1}\boxtimes (g_{1}\boxtimes h_{1}))\boxplus\\ &\cdots\boxplus (\ux_{m}\boxtimes \y_{s}\boxtimes \z_{a},f_{m}\boxtimes (g_{s}\boxtimes h_{a})))
\end{align*}
so the only question to be answered is whether or not $(f_{u}\boxtimes g_{v})\boxtimes h_{w}$ and $f_{u}\boxtimes (g_{v}\boxtimes h_{w})$ agree for all $u,v$ and $w$. But these morphisms agree by coherence. Thus $\boxtimes$ is strictly associative. By definition , $((),1)$ is a strict unit for $\boxtimes$. 
\nin
We construct now a symmetry isomorphism $\gamma^{\boxtimes}:X\boxtimes Y\to Y\boxtimes X$. To do so, take objects of $\D^{s}$ of the form
\begin{align*}
X&=(c_{1}\boxtimes\cdots \boxtimes c_{n},(\ux,f)),\\
Y&=(d_{1}\boxtimes\cdots\boxtimes d_{r},(\y,g))
\end{align*}
where
\begin{align*}
\ux&=x_{1}\boxtimes\cdots \boxtimes x_{m},\\
\y&=y_{1}\boxtimes\cdots\boxtimes y_{s}.
\end{align*}
We will define first the symmetry isomorphism in this case
\[
\gamma^{\boxtimes}:X\boxtimes Y\to Y\boxtimes X.
\]
By definition, we have that 
\begin{align*}
\Theta(X\boxtimes Y)&=(f\boxtimes g)^{*}(\Delta(\ux\boxtimes \y)),\\
\Theta(Y\boxtimes X)&=(g\boxtimes f)^{*}(\Delta(\y\boxtimes \ux)).
\end{align*}
We want to define a natural isomorphism
\[
\gamma^{\boxtimes}_{f,g}:(f\boxtimes g)^{*}(\Delta(\ux\boxtimes \y))\to (g\boxtimes f)^{*}(\Delta(\y\boxtimes \ux)).
\]
By definition of a fibered category, there is a unique morphism $\gamma^{\boxtimes}_{f,g}$ making the following diagram commutative
\[
\xymatrix{
(f\boxtimes g)^{*}(\Delta(\ux\boxtimes\y))\ar[rr]^-{} \ar@{-->}[rd]^-{\gamma^{\boxtimes}_{f,g}}\ar@{|->}[dd]&&\Delta(\ux\boxtimes \y)\ar@{|->}[dd]\ar[rd]^-{\Delta(\gamma^{\boxtimes})}& \\
&(g\boxtimes f)^{*}(\Delta(\y\boxtimes\ux))\ar[rr]^-{}\ar@{|->}[dd]&&\Delta(\y\boxtimes \ux)\ar@{|->}[dd]\\
\Phi(\c\boxtimes \d)\ar[rr]^-{f\boxtimes g}\ar[dr]^{\gamma^{\boxtimes}}&&\Lambda(\Delta(\ux\boxtimes \y))\ar[dr]^{\Lambda(\Delta(\gamma^{\boxtimes}))}&\\
&\Phi(\d\boxtimes \c)\ar[rr]^{g\boxtimes f}&&\Lambda(\Delta(\y\boxtimes \ux))
}
\]
Because of the uniqueness condition, the morphism $\gamma^{\boxtimes}_{f,g}$ is natural and satisfies the required coherences. We extend the definition of the $\gamma^{\boxtimes}$ to all the objects of $\D^{s}$ in the following way. Suppose that 
\begin{align*}
X=&(c_{1}\boxtimes\cdots \boxtimes c_{n},(\ux_{1},f_{1})\boxplus\cdots \boxplus(\ux_{m},f_{m})),\\
Y=&(d_{1}\boxtimes\cdots \boxtimes d_{r},(\y_{1},g_{1})\boxplus\cdots \boxplus(\y_{s},g_{s}))
\end{align*}
are two general objects of $\D^{s}$. Then  we define $\gamma^{\boxtimes}$ to be the composite
\begin{align*}
\Theta(X\boxtimes Y)&=\sum_{i=1}^{m}\sum_{j=1}^{r}(f_{i}\boxtimes g_{j})^{*}(\Delta(\ux_{i}\boxtimes \y_{j}))\stackrel{\sum_{i=1}^{m}\sum_{j=1}^{r}\gamma^{\boxtimes}_{f_{i},g_{j}}}{\rightarrow}\\
&\sum_{i=1}^{m}\sum_{j=1}^{r}(g_{j}\boxtimes f_{i})^{*}(\Delta(\y_{j}\boxtimes \ux_{i}))
\stackrel{\cong}{\rightarrow} \sum_{j=1}^{r}\sum_{i=1}^{m}(g_{j}\boxtimes f_{i})^{*}(\Delta(\y_{j}\boxtimes \ux_{i}))\\
&=\Theta(Y\boxtimes X)
\end{align*}
Here by $\sum_{i=1}^{m}\sum_{j=1}^{r}x_{ij}$ we mean the sum with a consistent way of inserting parenthesis as before and the unlabeled isomorphism is the isomorphism obtained by rearranging the terms of summation using the isomorphism $\gamma^{\op}$. This way defined we see by using some coherence theory that $(\D^{s},\boxtimes,((),1),\gamma^{\boxtimes})$ has the structure of a permutative category.

Our next step is to construct distributivity maps
\begin{align*}
d^{l}:(X\boxtimes Y)\boxplus(X'\boxtimes Y)\to(X\boxplus X')\boxtimes Y,\\
d^{r}:(X\boxtimes Y)\boxplus(X\boxtimes Y')\to X\boxtimes(Y\boxplus Y),
\end{align*}
wherever they make sense; that is, for objects $X$ and $X'$ of $\D^{s}_{\c}$ and $Y$ and $Y'$ of $\D^{s}_{\d}$. By a straight forward computation we can see that 
\[
d^{l}=\text{id}:(X\boxtimes Y)\boxplus(X'\boxtimes Y)=(X\boxplus X')\boxtimes Y.
\]
On the other hand, if we write 
\begin{align*}
X&=(c_{1}\boxtimes\cdots \boxtimes c_{n},(\ux_{1},f_{1})\boxplus\cdots \boxplus(\ux_{m},f_{m})),\\
Y&=(d_{1}\boxtimes\cdots \boxtimes d_{r},(\y_{1},g_{1})\boxplus\cdots \boxplus(\y_{s},g_{s})),\\
Y'&=(d_{1}\boxtimes\cdots\boxtimes d_{l},(\y'_{1},g'_{1})\boxplus\cdots \boxplus(\y'_{s'},g'_{s'}))
\end{align*}
then we have that 
\begin{align*}
X\boxtimes Y\boxplus X\boxtimes Y'=&\\
(\c\boxtimes \d,&(\ux_{1}\boxtimes \y_{1},f_{1}\boxtimes g_{1})\boxplus\cdots\boxplus(\ux_{m}\boxtimes \y_{s},f_{m}\boxtimes g_{s})\\
&\boxplus(\ux_{1}\boxtimes \y'_{1},f_{1}\boxtimes g'_{1})\boxplus\cdots\boxplus(\ux_{m}\boxtimes \y'_{s'},f_{m}\boxtimes g'_{s'}))
\end{align*}
and 
\begin{align*}
X\boxtimes (Y\boxplus Y')=&\\
(\c\boxtimes \d,&(\ux_{1}\boxtimes \y_{1},f_{1}\boxtimes g_{1})\boxplus\cdots\boxplus(\ux_{1}\boxtimes \y'_{s'},f_{1}\boxtimes g'_{s'})\\
&\boxplus (\ux_{m}\boxtimes \y'_{1},f_{1}\boxtimes g'_{1})\boxplus\cdots\boxplus(\ux_{m}\boxtimes \y'_{s'},f_{m}\boxtimes g'_{s'})).
\end{align*}
Thus we define 
\[
d^{r}:(X\boxtimes Y)\boxplus(X\boxtimes Y')\to X\boxtimes(Y\boxplus Y)
\]
as an iteration of $\gamma_{\c\boxtimes \d}^{\boxplus}$. By a trivial but long computation one can see that these distributivity maps are natural and satisfy the coherences of a fibered bipermutative category as in Definition \ref{bipermutativeonfibers}.

To finish we need to show that we have addition on morphisms in $\D^{s}$ over the same morphism in $\C^{s}$. Thus suppose that $\alpha:\c\to \d$ is a morphism in $\C^{s}$ and $\beta:X\to Y$, $\delta:X'\to Y'$ are two morphism in $\D^{s}$ over $f$, where $X$, $X'$ are objects in $\D^{s}_{\c}$ and $Y$, $Y'$ are objects in $\D^{s}_{\d}$. As noted before, we have a natural coherent isomorphism
\[
\mu_{X,X'}:\Theta(X\boxplus_{\c}X')\stackrel{\cong}{\rightarrow}\Theta(X)\op \Theta(X').
\]
This isomorphism is obtained by rearranging the parenthesis in the summations. With this in mind, we define $\alpha\boxplus \beta$ to be the following composite
\[
\Theta(X\boxplus_{\c}X')\stackrel{\mu_{X,X'}}{\rightarrow}\Theta(X)\op \Theta(X')\stackrel{\alpha\op \beta}{\rightarrow} \Theta(Y)\op \Theta(Y')\stackrel{\mu_{Y,Y'}^{-1}}{\rightarrow}\Theta(Y\boxplus_{\c}Y').
\]
This way defined we see that $\boxplus$ is strictly associative and that 
\[
\xymatrix{
X\boxplus_{\c} X'\ar[r]^-{\alpha\boxplus \beta}\ar[d]_-{\gamma^{\boxplus}_{\c}}& Y\boxplus_{\d} Y' \ar[d]^-{\gamma^{\boxplus}_{\d}}\\
X'\boxplus_{\c} X\ar[r]^-{\beta\boxplus \alpha}& Y'\boxplus_{\d} Y 
}
\]
is a commutative diagram.
\nin
In addition, note that by definition $\Theta(0_{\c})=0_{\Phi(\c)}$ and we have a coherent morphism $0_{f}:0_{\Phi(\c)}\to 0_{\Phi(\d)}$ in $\D$, therefore we define a morphism $0_{f}:0_{\c}\to 0_{\d}$ in $\D^{s}$ to be the morphism $0_{f}:0_{\Phi(\c)}\to 0_{\Phi(\d)}$ in $\D$. The coherences satisfied by $0_{f}$ in $\D$ imply that 
\[
0_{f}\boxplus \alpha =\alpha=\alpha\boxplus 0_{f}.
\]
Finally, since $\Lambda:\D\to \C$ is a fibered category, it follows easily that $\Lambda^{s}:\D^{s}\to \C^{s}$ is also a fibered category. We conclude then that $\Lambda^{s}:\D^{s}\to \C^{s}$ is a fibered bipermutative category. This proves the theorem.
\qed

\section{Fibered categories}
 
In this section we review briefly some aspects of fibered categories. In particular, we recall the fact that there is a one to one correspondence between fibered categories and contravariant functors to the category of small categories. We use this construction as motivation for our study of fibered bipermutative categories. Indeed, as mentioned before, using this correspondence, given a fibered bipermutative category $\Lambda:\D\to \C$ we can assign a functor $\C^{op}\to \P$. However, this functor does not capture all the information present in a fibered bipermutative category and because of this we need to enlarge the category $\C^{op}$. Thus we show that for a given fibered bipermutative category $\Lambda:\D\to \C$ we can correspond a functor $\Psi:\A\to \P$, where $\A$ is a wreath product for a functor $(\C^{op})^{*}:\textbf{Inj}\to \P$. The category $\C^{op}$ canonically includes into $\A$ and under this inclusion, the functor $\Psi$ recovers the construction for fibered categories. Moreover, the functor $\Psi$ recovers the multiplicative structure on $\D$. 
\nin 
The goal of this section is to construct the functor $\Psi$. We begin by recalling the definition of a fibered category. 
 
\begin{definition} If $F:\D\to \C$ is a functor and $c$ an object of $\C$, we denote by $\D_{c}$ the fiber $F^{-1}(c)$; that is, the subcategory of $\D$ whose objects are the objects $d$ such that $F(d)=c$ and the morphisms are the morphisms mapping to the identity of $c$.
\end{definition} 

\begin{definition}{\label{fiberedcategory}}
A fibered category over $\C$ is a functor $F:\D\to \C$ that satisfies the following properties:
\begin{itemize}
	\item for any morphism $f:c\to c'$ in $\C$ and any object $d'$ of $\D_{c'}$, we can find an object $d$ of $\D_{c}$ and a morphism $g:d\to d'$ such that $F(g)=f$,
	\item given any pair of morphisms $g:d\to d'$ and $g':d''\to d'$ in $\D$, let $f=F(g):c\to c'$ and $f':F(g'):c''\to c'$ their image in $\C$. Then for any $\tilde{f}:c''\to c$ such that $f\tilde{f}=f'$, there is a unique morphism $\tilde{g}:d\to d''$ such that $g\tilde{g}=g'$ and $F(\tilde{g})=h$.
\end{itemize}
\end{definition}
In general for a functor $F:\D\to \C$ there is the possibility that there is an isomorphism  $f:c\to c'$ with $\D_{c}$ and $\D_{c'}$ not equivalent categories. This does not happen in the case of fibered categories as any such isomorphisms $f$ induces an equivalence of categories $f^{*}:\D_{c'}\to \D_{c}.$

\begin{definition} Let $F:\D\to \C$ be a functor and $g:d\to d'$ an arrow in $\D$. Take $g':d''\to d'$ in $\D$, and let $f=F(g):c\to c'$, $f'=F(g'):c''\to c'$ their image in $\C$. We say that $g$ is a cartesian arrow if for any $\tilde{f}:c''\to c$ such that $f\tilde{f}=f'$, there is a unique morphism $\tilde{g}:d''\to d$ such that $g\tilde{g}=g'$ and $F(\tilde{g})=h$. Such a cartesian arrow $g$ over $f$ is called a pullback of $f$.
\end{definition}
This definition is summarized by the following diagram
 \[
   \xymatrix@R=6pt{
   {d''}\ar@{|->}[dd] \ar@{-->}[rd]_{\tilde{g}}
   \ar@/^/[rrd]^{g'} \\
   & {}d\ar@{|->}[dd]\ar[r]_{g}
   & {}d'\ar@{|->}[dd] \\
   {} c''\ar[rd]_{\tilde{f}}\ar@/^/[rrd]^{\ \ f'}\\
   & {} c \ar[r]_-{f}
   & {} c'. 
   }
   \]
Thus a fibered category over $\C$ is precisely a functor $F:\D\to \C$ such that given any morphism $f:c\to c'$ in $\C$ and $d'$ any object in $\D_{c'}$, there exists an object $d$ in $\D_{c}$ and a cartesian arrow $g:d\to d'$ such that $F(g)=f$.

\noindent{\bf{Example: }} If $u$ is an object of $\C$, then the comma category $\C/u$ is fibered category over $\C$ by considering $\theta_{u}:\C/u\to \C$ defined on objects by $\theta_{u}(f:v\to u)=v$ and if $\phi:(g:w\to u)\to (f:v\to u)$ is a morphism in $\C/u$; that is, $\phi:w\to v$ is such that $f\phi=g$, then $\theta_{u}(\phi)=\phi$.

\begin{definition}
If $F:\D\to \C$ and $G:\E\to \C$ are two fibered categories over $\C$, then a morphism of fibered categories, $\phi:\D\to \E$, is a functor such that $F=G\phi$ and $\phi$ sends cartesian arrows to cartesian arrows. 
\end{definition}

Given two fibered categories $F:\D\to \C$ and $G:\E\to \C$, we denote by $\text{Hom}_{\C}(\D,\E)$ the category whose objects are the morphisms of fibered categories $\phi:\D\to \E$. If $\phi,\psi:\D\to \E$ are objects in $\text{Hom}_{\C}(\D,\E)$, then a morphism $\alpha:\phi\to \psi$ is a base point preserving natural transformation; that is, a natural transformation $\alpha$ such that for every object $d$ of $\D$ over $c$, the morphism $\alpha_{d}$ is a morphism in the category $\E_{c}$.

Note that the definition of a fibered category over $\C$ states that given any morphism $f:c\to c'$ and any object $d'$ of $\D_{c'}$, then we can choose an object $f^{*}d'$ over $c$ and a morphism $\eta:f^{*}d'\to d'$ over $f$. We call such a morphism a pullback of $f$. The morphism $\eta$ is uniquely determined up to composition with an isomorphism in the fiber $\D_{c}$. So for each $f$ and each object $d'$ over $c'$ we can fix such a pullback $\eta$ over $f$. Also, if $\alpha:d'\to e'$ is a morphism in $\D_{c'}$, then by the uniqueness part in the definition of a cartesian arrow, we see that there is a unique morphism $f^{*}\alpha$ in the category $\D_{c}$ such that the following diagram commutes
\[
\xymatrix{
f^{*}d'\ar[r]^-{\eta_{d'}}\ar[d]_-{f^{*}\alpha}&d'\ar[d]^-{\alpha}\\
f^{*}e'\ar[r]^-{\eta_{e'}}&e'.
}
\]
This defines a functor $f^{*}:\D_{c'}\to \D_{c}$. Notice that if $f$ and $g$ are composable morphisms in $\C$, then $(fg)^{*}$ does not necessarily agree with $g^{*}f^{*}$. However, we can find a canonical isomorphism between $(fg)^{*}$ and $g^{*}f^{*}$. Thus for a fibered category $F:\D\to \C$, and a choice of pullbacks for each morphism $f$ in $\C$,  we can associate the following correspondence:
\begin{itemize}
	\item for an object $c$ of $\C$ we associate the category $\D_{c}$,
	\item for a morphism $f:c\to c'$ we associate the functor $f^{*}:\D_{c'}\to \D_{c}$.
\end{itemize}
At a first glance it would seem that this defines a functor $\C^{op}\to \cat$, but as we just pointed out, this is not the case. However, $\cat$ has the structure of a $2$-category and the previous assignment gives rise to a lax $2$-functor and we have the following theorem.

\begin{theorem}\label{2lax}
The above assignment defines a one to one correspondence between isomorphism classes of fibered categories with a choice of pullbacks and isomorphisms of contravariant lax $2$-functors $\C\to\cat$.
\end{theorem}
\Proof
See \cite[3.1.2 and 3.1.3]{Vistolo}.
\qed

In fact up to isomorphism, we can replace any given fibered category $F:\D\to \C$ with an isomorphic fibered category $F':\D'\to \C$ such that the corresponding assignment $\C^{op}\to\cat$ is indeed a functor. Conversely, given any functor $\C^{op}\to \cat$, we can associate a fibered category $F:\D\to \C$. This defines a one to one correspondence between isomorphism classes of fibered categories and isomorphism classes of functors $\C^{op}\to \cat$. We show this in the following theorem.

\begin{theorem}\label{functor}
There is a one to one correspondence between isomorphism classes of fibered categories and isomorphism classes of functors $\C^{op}\to \cat$.
\end{theorem}
\Proof
We only sketch the proof of this theorem. For the details we refer the reader to \cite[Chapter 3]{Vistolo}. Suppose first that $F:\D\to \C$ is a fibered category. We want to define a functor $f:\C^{op}\to \cat$. For an object $u$ define $f(u)$ to be the category $\text{Hom}_{\C}(\C/u,\D)$. Given a morphism $g:u\to v$ in $\C$, by composing with $g$ we obtain a morphism of fibered categories $g_{*}:\C/u\to \C/v$. This induces $f(g):\text{Hom}_{\C}(\C/v,\D)\to \text{Hom}_{\C}(\C/u,\D)$.  It follows easily that this defines a functor $f:\C^{op}\to \cat$. On the other hand, given a functor $f:\C^{op}\to \cat$ we can construct a fibered category $F:\D\to \C$ in the following way. The objects of $\D$ are the pairs $(x,c)$, where $c$ is an object in $\C$ and $x$ is an object in $F(c)$. If $(x,c)$ and $(y,d)$ are two such pairs, then a morphism in $\D$ from $(x,c)$ to $(y,d)$ is a pair $(\alpha,f)$, where $f:c\to d$ is a morphism in $\C$ and $\alpha:x\to F(f)y$ is a morphism in the category $F(c)$. If $(\alpha,f):(x,c)\to (y,d)$, $(\beta,g):(y,d)\to (z,e)$, then 
\[
(\beta,g)\circ (\alpha,f):=(F(f)(\beta)\circ\alpha,g\circ f):(x,c)\to (z,e).
\]  
It's easy to see that this defines a category. In addition we have a functor $F:\D\to \C$ that sends a pair $(x,c)$ to $c$ and a morphism $(\alpha,f)$ to $f$. This makes $\D$ into a fibered category over $\C$ whose fiber $\D_{c}$ over an object $c$ of  $\C$ is naturally isomorphic to the category $f(c)$. We leave to the reader to show that this is indeed a one to one correspondence.
\qed

We can apply this procedure to a bipermutative category on fibers  $\Lambda:\D\to \C$. This way we obtain a functor $f:\C^{op}\to \cat$ that has the further property that each $f(u)$ has the structure of a permutative category and the functor $f$ can be seen as a functor $f:\C^{op}\to \P$. However, this functor does not behaves well under the multiplicative structure on $\D$. We show however, that the fibered category $\Lambda:\D\to \C$ determines and is determined by a functor $\Psi:\A\to \P$ that behaves well under multiplication. Here $\A$ is a category that is obtained as a wreath product $\textbf{Inj}\int(\C^{op})^{*}$ for a functor $(\C^{op})^{*}:\textbf{Inj}\to \P$. The category $\C$ naturally embeds into $\A$ and under this embedding we recover the functor $f:\C^{op}\to \P$ arising from the fibered category $\Lambda:\D\to \C$. After the construction of the functor $\Psi$, we show in theorem \ref{verifyconditions} that it satisfies conditions $(\bf{c.1})-(\bf{c.14})$ of theorem \ref{multiEsigma}. We will explain how to construct the functor $\Psi:\A\to \P$ in what follows.

We begin by constructing the category $\A$. This is a straight forward generalization of \cite[Definition 5.1]{Elmendorf}. Let us denote by $\bf{Inj}$ the category whose object set is the set of integers $n\ge 0$, where we identify the integer $n$ with the set $\n=\{1,...,n\}$. The morphism set from $n$ to $m$ is the set of injective functions from $\n=\{1,...,n\}$ to $\m=\{1,...,m\}$ with composition the composition of functions. Define the functor 
\[
(\C^{op})^{*}:\bf{Inj}\to \P
\] 
such that $(\C^{op})^{*}(n)=(\C^{op})^{n}$. If $q:\n\to \m$ is a morphism in $\textbf{Inj}$; that is, $q:\n\to \m$ is an injective map, then $(\C^{op})^{*}(q):(\C^{op})^{n}\to (\C^{op})^{m}$ is the functor that for an $n$-tuple $\u=(u_{1},...,u_{n})$ corresponds $q_{*}\u=(u'_{1},...,u'_{m})$, where 
\begin{equation*}
u_{j}'=\left\{ 
\begin{array}{cc}
u_{i} & \text{if }q^{-1}(j)=\{i\}, \\ 
1 & \text{if }q^{-1}(j)=\emptyset,
\end{array}
\right. 
\end{equation*}
and if $\f=(f_{1},...,f_{n}):\u\to \v$ is a morphism in $(\C^{op})^{n}$, then $\C^{*}(\f)=\f':q_{*}\u\to q_{*}\v$, where $\f'=(f_{1}',...,f_{m}')$ and 
\begin{equation*}
f_{j}'=\left\{ 
\begin{array}{cc}
f_{i} & \text{if }q^{-1}(j)=\{i\}, \\ 
\text{id}_{1} & \text{if }q^{-1}(j)=\emptyset.
\end{array}
\right. 
\end{equation*}
Then associated to this functor, there is a wreath product category $\A=\bf{Inj}\int (\C^{op})^{*}$. More explicitly, the objects of $\A$ are the sequences of the form $(u_{1},...,u_{n})$ where $n\ge 0$ and a morphism from $\u=(u_{1},...,u_{n})$ to  $\v=(v_{1},...,v_{n})$ is a pair $(q,\underline{f})$, where $q:\n\to \m$ is an injection and $\underline{f}:q_{*}\u\to \v$ is a morphism in $(\C^{op})^{n}$. Note that $\A$ is a permutative category by concatenation. Indeed, for objects $\u=(u_{1},...,u_{n})$ and $\v=(v_{1},...,v_{m})$, we can define $\u\odot\v=(u_{1},...,u_{n},v_{1},...,v_{m})$ and similarly on morphisms. The unit of $\A$ is the empty tuple $()$.
\nin
Suppose that $\u=(u_{1},...,u_{n})$ is an object of $\A$. We will denote by $\C/\u$ the product category $\C/u_{1}\x\cdots \x\C/u_{n}$. We can see $\C/\u$ as a category over $\C$ by defining the functor $\theta_{\u}:\C/\u\to \C$ as follows. Given objects $f_{i}:w_{i}\to u_{i}$ in $\C/u_{i}$, then $\theta(f_{1},...,f_{n})=w_{1}\ot\cdots \ot w_{n}$. Also, if $\phi_{i}:w_{i}\to w_{i}'$ is a morphism in $\C/u_{i}$, then $\theta(\phi_{1},...,\phi_{n})=\phi_{1}\ot\cdots\ot\phi_{n}$. As explained above, this makes $\theta_{u_{1}}:\C/u_{1}\to \C$ into a fibered category in the case where $n=1$. This is not the case in general for $n>1$. When $\u=()$ is the empty tuple, then we will understand by $\C/()$ the trivial category with only object $\text{id}:1\to1$ and over $\C$. With this convention then the category $\text{Hom}_{\C}(\C/(),\D)$ is isomorphic to the category $\D_{1}$. Note that since we are assuming that every morphism in $\C$ is an isomorphism, then a morphism of fibered categories $\C/u\to \D$ is just a functor $F:\C/u\to \D$ such that $\Lambda\circ F=\theta_{u}$.

Let us construct now the functor $\Psi:\A\to \P$. Take $\u=(u_{1},...,u_{n})$ an object of $\A$. Define $\Psi(\u)$ to be the category whose objects are the functors $F:\C/\u\to \D$ such that $\Phi\circ F=\theta_{\u}$ and such that $F$ is of the form 
\[
F=\sum_{i=1}^{r}F_{i1}\ot\cdots\ot F_{in},
\]
where each $F_{ij}$ is an object of $\text{Hom}_{\C}(\C/u_{i},\D)$. This means that given an object $\f=(f_{1},...,f_{n})$ in $\C/\u$, then 
\[
F(\f)=F_{11}(f_{1})\ot\cdots\ot F_{1n}(f_{n})\op\cdots\op F_{r1}(f_{1})\ot\cdots\ot F_{rn}(f_{n})
\]
and given $\underline{\phi}=(\phi_{1},...,\phi_{n})$ a morphism in $\C/\u$, then 
\[
F(\underline{\phi})=F_{11}(\phi_{1})\ot\cdots\ot F_{1n}(\phi_{n})\op\cdots\op F_{r1}(\phi_{1})\ot\cdots\ot F_{rn}(\phi_{n}).
\]
This defines the objects of $\Psi(\u)$. Given 
\[
F=\sum_{i=1}^{r}F_{i1}\ot\cdots\ot F_{in} \text{\ and  \ } G=\sum_{j=1}^{s}G_{j1}\ot\cdots\ot G_{jn}
\] 
two such functors, then a morphism $\alpha:F\to G$ in $\Psi(\u)$ is a base point preserving natural transformation $F\to G$; that is, $\alpha$ is a natural transformation such that for an object $\f$ of $\C/\u$, with $f_{i}:w_{i}\to u_{i}$, then $\alpha_{\f}:F(\f)\to G(\f)$ is a morphism in $\D_{\ot w_{i}}$. As the composition of a base preserving natural transformation is also a base preserving natural transformation we obtain this way a well defined category $\Psi(\u)$. When $\u=()$ is the empty tuple, then $\Psi(())$ is then the category whose objects are functors $F:\C/()\to \D$ such that $F=\sum_{i}^{n}F_{i}$ with $F_{i}$ an object of $\text{Hom}_{\C}(\C/(),\D)$ and with morphisms the base preserving natural transformations. This category is canonically isomorphic to $\D_{1}$ under the isomorphism $F\mapsto F(\text{id}:1\to 1)$.

We want to show that if $\Lambda:\D\to \C$ is a fibered bipermutative category, then the category $\Psi(\u)$ has the structure of a permutative category. To do so, we need to define a functor $\op:\Psi(\u)\x\Psi(\u)\to \Psi(\u)$ and show that is satisfies the respective coherences. Suppose then that $F=\sum_{i=1}^{r}F_{i1}\ot\cdots\ot F_{in}$ and $G=\sum_{j=1}^{s}G_{j1}\ot\cdots\ot G_{jn}$ are two objects of $\Psi(\u)$. Define
\[
F\op G=\sum_{i=1}^{r}F_{i1}\ot\cdots\ot F_{in}\op \sum_{j=1}^{s}G_{j1}\ot\cdots\ot G_{jn}.
\]
This way defined we see that given $\f$ an object in $\C/\u$
\[
(F\op G)(\f)=F(\f)\op G(\f).
\]
Given $\underline{\phi}=(\phi_{1},...,\phi_{n}):\f\to \f'$ a morphism in $\C/\u$, then
\[
(F\op G)(\underline{\phi})=F(\underline{\phi})\op G(\underline{\phi}) :F(\f)\op G(\f)\to F(\f')\op G(\f') 
\]
is well defined as both $F(\underline{\phi})$ and  $G(\underline{\phi})$ are morphisms over $\ot\phi_{i}$ and we can add two such morphisms. It is easy to see that since the sum of morphisms over $\ot\phi_{i}$ is natural then $F\op G$ is also an object of $\Psi(\u)$.
This defines $\op$ on the objects of $\Psi(\u)$. Suppose now that $\alpha:F\to G$ and $\beta:H\to K$ are two morphism in $\Psi(\u)$. Given $\f$ an object in $\C/\u$, define 
\[
(\alpha\op \beta)_{\f}=\alpha_{\f}\op \beta_{\f}:F(\f)\op H(\f)\to G(\f)\op K(\f). 
\]
It follows at once that this way defined, $\alpha\op \beta$ is a base preserving natural transformation and that $\op:\Psi(\u)\x\Psi(\u)\to \Psi(\u)$ is a functor. We claim that $\op$ is strictly associative. Indeed, if $F,G$ and $H$ are objects of $\Psi(\u)$, then both $(F\op G)\op H$ and $F\op (G\op H)$ agree as functors as by assumption each $\op_{u}$ is strictly associative and the addition of morphisms is also strictly associative. 
\nin
Next we define an additive unit. We need to define an object $\Os_{\u}$ of $\C/\u$. To do so, given $u$ an object of $\C$ we define $\Os_{u}$ to be the object of $\text{Hom}_{\C}(\C/u,\D)$ defined by as follows. If $f:w\to u$ is an object in $\C/u$, then $\Os(f)=0_{w}$, where $0_{w}$ is the unit of the permutative category $\D_{w}$. If $\phi:w\to w'$ is a morphism in $\C/u$ then $\Os(\phi)=0_{\phi}:0_{w}\to 0_{w'}$ is the coherent morphism over $\phi$ as whose existence is guaranteed by hypothesis. 
Then we define $\Os_{\u}=\Os_{u_{1}}\ot\cdots \ot \Os_{u_{n}}$. This way defined $\Os_{\u}$ is an object of $\Psi(\u)$. We claim this is a strict unit;  that is, we want to see that for every $F$ object of $\Psi(\u)$
\[
F\op \Os_{\u}=F=\Os_{\u}\op F.
\]
Indeed, if $\f$ is an object of $\C/\u$, then by definition
\begin{align*}
(F\op \Os_{\u})(\f)&=F(\f)\op \Os_{\u}(\f)=F(\f)\op 0_{\ot w_{i}}=F(\f)\\
(\Os_{\u}\op F)(\f)&=\Os_{\u}(\f)\op F(\f)=0_{\ot w_{i}}\op F(\f)=F(\f)
\end{align*}
and similarly on morphisms. 
\nin
We also need to construct a symmetry isomorphism. Suppose then that $F$ and $G$ are two objects of $\Psi(\u)$. Define $\gamma:F\op G\to G\op F$ to be the transformation that for each object $\f$ of $\C/\u$ assigns
\[
\gamma_{\f}=\gamma^{\op}_{\ot w_{i}}:F(\f)\op G(\f)\to G(\f)\op F(\f).
\]
This is a natural transformation as addition commutes with the different $\gamma^{\op}_{u}$ and is base preserving by definition. In order to conclude that $(\Psi(\u),\op,\gamma)$ is a permutative category we need to verify that the following diagrams are commutative
\[
\xymatrix{
F\op G\ar[rr]^-{=}\ar[dr]_-{\gamma} &&F\op G \\
&G\op F\ar[ur]_-{\gamma},&   
}
\]
\[
\xymatrix{
F\op G\op H\ar[rr]^-{\gamma}\ar[dr]_{\textsl{id}\op\gamma} &&H\op F\op G\\
&F\op H\ot G\ar[ur]_{\gamma\op \textsl{id}}.&   
}
\]
But the commutativity of these diagrams follows as we have the commutativity of the corresponding diagrams when we evaluate on objects and morphisms, therefore $\Psi(\v)$ is an object of $\P$.

We now want to show that we can extend the definition of $\Psi$ as to get a functor $\Psi:\A\to \P$.  To start, for a given morphism $f:u\to v$ in $\C$ we define a functor 
\begin{align*}
f^{*}:\text{Hom}_{\C}(\C/v,\D)&\to \text{Hom}_{\C}(\C/u,\D)\\
F&\mapsto f^{*}F\\
\alpha&\mapsto f^{*}\alpha.
\end{align*}

The functor $f^{*}F$ is defined as follows. For an object $g:w\to u$ of $\C/u$ $f^{*}F(g)=F(f\circ g)$ and if $\phi:f\to f'$ is a morphism in $\C/u$ then $\phi$ can be seen as a morphism in $\C/v$ between $f\circ g$ and $f'\circ g$, thus we define $f^{*}F(\phi)=F(\phi)$. Similarly we define $f^{*}\alpha$ for a base preserving natural transformation $\alpha:F\to G$. Suppose now that $\u=(u_{1},...,u_{n})$ and $\v=(v_{1},...,v_{m})$ are two objects of $\A$ and that $(q,\f):\u\to \v$ is a morphism in $\A$, thus $\f:q_{*}\u=\u'\to \v$ is a morphism in $(\C^{op})^{m}$; that is,  $f_{i}:v_{i}\to u'_{i}$ is a morphism in $\C$. We want to define a lax$_{*}$ map $\Psi(q,\f):\Psi(\u)\to \Psi(\v)$. Suppose first that $m=n$ and thus $q=\sigma:\n\to \n$ is a permutation. Then each $f_{i}:v_{i}\to u_{\sigma^{-1}(i)}$ is a morphism in $\C$.
Take  
\[
F=\sum_{i=1}^{r}F_{i1}\ot\cdots\ot F_{in},
\]
an object of $\C/\u$, where each $F_{ij}$ is an object of $\text{Hom}_{\C}(\C/u_{i},\D)$ . Define
\[
\Psi(\sigma,\f)F=\sum_{i=1}^{r}f^{*}_{1}F_{i\sigma^{-1}(1)}\ot\cdots\ot f^{*}_{n}F_{i\sigma^{-1}(n)}.
\]
By definition it is clear that this is an object in $\Psi(\v)$. Suppose now that $\alpha:F\to G$ is a morphism in $\Psi(\u)$, where
\[
F=\sum_{i=1}^{r}F_{i1}\ot\cdots\ot F_{in} \ \ \ {\text{and}}\ \ G=\sum_{j=1}^{s}G_{j1}\ot\cdots\ot G_{jn}.
\]
Then we need to define 
\[
\Psi(\sigma,\f)\alpha:\Psi(\sigma,\f)F\to \Psi(\sigma,\f)G.
\]
For each $1\le i\le n$, the symmetry isomorphism $\gamma^{\ot}$ determines a coherent natural isomorphism
\[
\tau_{\sigma,f^{*}_{1}F_{i\sigma^{-1}(1)},...,f^{*}_{n}F_{i\sigma^{-1}(n)}}:f^{*}_{1}F_{i\sigma^{-1}(1)}\ot\cdots\ot f^{*}_{n}F_{i\sigma^{-1}(n)}\to f^{*}_{\sigma(1)}F_{i1}\ot\cdots\ot f^{*}_{\sigma(n)}F_{in}
\]
and similarly for $G$. Since $\Phi$ is assumed to be a strict map, then for each object $\g=(g_{1},...,g_{n})$, with $g_{i}:w_{i}\to v_{i}$ of $\C/\v$ the morphism 
\begin{align*}
&\tau_{\sigma,f^{*}_{1}F_{i\sigma^{-1}(1)},...,f^{*}_{n}F_{i\sigma^{-1}(n)}}(\g)\\
&:f^{*}_{1}F_{i\sigma^{-1}(1)}(g_{1})\ot\cdots\ot f^{*}_{n}F_{i\sigma^{-1}(n)}(g_{n})\to f^{*}_{1}F_{i\sigma(1)}(g_{\sigma(1)})\ot\cdots\ot f^{*}_{n}F_{i\sigma(n)}(g_{\sigma(n)})
\end{align*}
is a morphism over $\tau_{\sigma,w_{1},...,w_{n}}:w_{1}\ot\cdots\ot w_{n}\to w_{\sigma(1)}\ot\cdots\ot w_{\sigma(n)}$ which is the coherent isomorphism in $\C$ provided by $\gamma^{\ot}$.

Then we define $(\Psi(\sigma,\f)\alpha)_{\g}$, by the following commutative diagram
\[
\xymatrix@C=-8pc{
\sum_{i=1}^{r}f^{*}_{1}F_{i\sigma^{-1}(1)}(g_{1})\ot\cdots\ot f^{*}_{n}F_{i\sigma^{-1}(n)}(g_{n}) \ar[dr]^-{(\Psi(\sigma,\f)\alpha)_{\g}}\ar[dd]_-{a}& \\ &\sum_{j=1}^{s}f^{*}_{1}G_{j\sigma^{-1}(1)}(g_{1})\ot\cdots\ot f^{*}_{n}G_{j\sigma^{-1}(n)}(g_{\sigma(n)})
\\
\sum_{i=1}^{r}f^{*}_{\sigma(1)}F_{i1}(g_{\sigma(1)})\ot\cdots\ot f^{*}_{\sigma(n)}F_{in}(g_{\sigma(n)})\ar[dr]_-{\alpha_{\underline{\sigma^{*}(f\circ g)}}}&\\ &\sum_{j=1}^{s}f^{*}_{\sigma(1)}G_{j1}(g_{\sigma(1)})\ot\cdots\ot f^{*}_{\sigma(n)}G_{jn}(g_{\sigma(n)})\ar[uu]_-{b}. 
}
\]
Here 
\begin{align*}
a&=(\sum_{i=1}^{r}\tau_{\sigma,f^{*}_{1}F_{i\sigma^{-1}(1)},...,f^{*}_{n}F_{i\sigma^{-1}(n)}})(\g)\\
b&=(\sum_{i=1}^{r}\tau^{-1}_{\sigma,f^{*}_{1}G_{i\sigma^{-1}(1)},...,f^{*}_{n}G_{i\sigma^{-1}(n)}})(\g).
\end{align*}
The naturality of $\Psi(\sigma,\f)\alpha$ follows by the naturality of $\alpha$ and the maps 
\[
\tau_{\sigma,f^{*}_{1}F_{i1},...,f^{*}_{n}F_{in}} \text{\ \ and \ \ } \tau_{\sigma,f^{*}_{1}G_{j1},...,f^{*}_{n}G_{jn}}.
\] 
We need to check that $\Psi(\sigma,\f)\alpha$ is base point preserving. To do so, note that $\Lambda((\Psi(\sigma,\f)\alpha)_{\g})=\tau_{\sigma,w_{1},...,w_{n}}^{-1}\circ \tau_{\sigma,w_{1},...,w_{n}}=\text{id}$, and thus $(\Psi(\sigma,\f)\alpha)_{\g}$ is a morphism in the category $\D_{\ot w_{i}}$. This defines $\Psi(q,\f):\Psi(\u)\to \Psi(\v)$ in the case that $q:\n\to \m$ is an isomorphism. Let us define $\Psi(q,\f)$ for the case where $q:\n\to \underline{n+1}$ is the injective map that misses the value $n+1$, $\v=(u_{1},...,u_{n},1)$ and each $\f=\text{id}$. Thus in this case $u'_{i}=u_{i}$ for $1\le i\le n$ and $u'_{n+1}=1$. Let us fix $I:\C/1\to \D$ a morphism of fibered categories such that $I(\text{id})=1$.  Then for an object  
\[
F=\sum_{i=1}^{r}F_{i1}\ot\cdots\ot F_{in} 
\]
of $\Psi(\u)$, we define 
\[
\Psi(q,\f)F=\sum_{i=1}^{r}F_{i1}\ot\cdots\ot F_{in}\ot I. 
\]
This way defined we see that $\Psi(q,\f)F$ is an object of $\Psi(\v)$. On the other hand, if $\alpha:F\to G$ is a morphism in the category $\Psi(\u)$,
where
\[
F=\sum_{i=1}^{r}F_{i1}\ot\cdots\ot F_{in} \ \ \ {\text{and}}\ \ G=\sum_{j=1}^{s}G_{j1}\ot\cdots\ot G_{jn}.
\]
then we define $\Psi(q,\f)\alpha$ to be the morphism defined in the following diagram
\[
\xymatrix{
\sum_{i=1}^{r}F_{i1}\ot\cdots\ot F_{in}\ot I\ar[d]_-{\cong}\ar[r]^-{\Psi(q,\f)\alpha} & G=\sum_{j=1}^{s}G_{j1}\ot\cdots\ot G_{jn}\ot I\ar[d]^-{\cong}\\
(\sum_{i=1}^{r}F_{i1}\ot\cdots\ot F_{in})\ot I\ar[r]_{\alpha\ot \text{id}} &(\sum_{j=1}^{s}G_{j1}\ot\cdots\ot G_{jn})\ot I,
}
\] 
where the vertical arrows are obtained by application of the distributivity maps in $\D$. It is easy to see that this way defined $\Psi(q,\f)\alpha$ is a base point preserving natural transformations.
This defines $\Psi(q,\f)$ in this case. Note that any morphism in $\A$ can be factorize into a composition of the previous form, hence this way we define $\Psi(q,\f)$ for every morphism $(q,\f)$ in $\A$. By a direct inspection, we see that $\Psi(\sigma,\f)$ is a lax$_{*}$ map, in fact it is a strict map. We conclude then that $\Psi:\A\to \P$ is a functor. Thus we have proved the following theorem.

\begin{theorem}
To every fibered bipermutative category $\Lambda:\D\to \C$ we can associate a functor $\Psi:\A\to \P$. The functor determines and is determined by $\Lambda:\D\to \C$ up to canonical isomorphism.
\end{theorem}

In Theorem \ref{verifyconditions}, we will show that the functor $\Psi:\A\to \P$ has special properties arising from the multiplicative structure of the fibered bipermutative category $\Lambda:\D\to \C$. The functor $\Psi$ motivates the study in general of functors out of a permutative category to the multicategory of permutative categories or in general to a multicategory. We show in the next section that we can give such functors the structure of a multicategory.

\section{Multicategories and general construction}\label{secconstruction} 

In this section we review multicategories.  Roughly speaking a multicategory or colored operad, is a generalization of both operads and symmetric monoidal categories. Our main goal in this section is to show that given a permutative category $(\E,\ot,1)$ and a multicategory $\M$, we can give the structure of a multicategory to the category of functors $\E\to \M$.  

\begin{definition}
A multicategory $\M$ consists of the following:
\begin{enumerate}
	\item A collection of objects which we usually denote $\O_{\M}$.
	\item For $k\ge 0$ and any $(k+1)$-tuple of objects $a_{1},...,a_{k}$ and $b$, a set 
	\[
	\M(a_{1},...,a_{k};b),
	\]
	called the set of ``colored'' $k$-morphisms.
	\item A right action of $\Sigma_{k}$ on the collection of all $k$-morphisms, where for $\sigma\in \Sigma_{k}$
	\[
	\sigma^{*}:\M(a_{1},...,a_{k};b)\to \M(a_{\sigma(1)},...,a_{\sigma(k)},b).
	\]
	\item A 1-morphism $1_{a}\in \M(a;a)$ called the unit, for each object a of $\M$.
	\item A multiproduct function 
	\begin{align*}
	\Gamma_{\M}&:\M(b_{1},...,b_{n};c)\x \M(a_{11},...,a_{1j_{1}};b_{1})\x\cdots \x\M(a_{n1},...,a_{nj_{n}};b_{n})\\
	&\to \M(a_{11},...,a_{nj_{n}};c)
	\end{align*}
\end{enumerate}
that satisfies properties Multi(1)-Multi(4) in \cite[Definition 2.1]{Elmendorf} that generalize the properties of an operad as in \cite[Definition 1.1]{MayGeom}.
\end{definition}

As mentioned before a multicategory is a generalization of both operads and symmetric monoidal categories, to be more precise, an operad can be seen as a multicategory with only one object. Also, if $(\E,\op,0)$ is a symmetric permutative category, we can see $\E$ as a multicategory by considering as objects the objects of $\E$ and for $c_{1},...,c_{k},d$ objects of $\E$, the set of $k$-morphisms
\[
\E(c_{1},...,c_{k};d):=\E(c_{1}\op\cdots \op c_{k},d).
\]
Here $c_{1}\op\cdots \op c_{k}$ means the iterated application of $\op$ by inserting parenthesis in a consistent way. 

In \cite{Elmendorf}, Elmendorf and Mandell gave the category $\P$ the structure of a multicategory. To describe the $k$-morphisms of this multicategory we need the following definition. 

\begin{definition} Let $\C_{1},...,\C_{k}$ and $\D$ be small permutative categories. A functor 
\[
f:\C_{1}\x\cdots \x\C_{k}\to \D
\]
is said to be a $k$-linear map, if $f(c_{1},...,c_{k})=0$ whenever $c_{i}=0$ for some $1\le i\le k$ and  in addition, we have natural transformations $\delta^{i}$, for $1\le i\le k$, that are thought of as distributivity maps,
\[
\delta^{i}_{f}=\delta^{i}:f(c_{1},...,c_{i},...,c_{k})\op f(c_{1},...,c'_{i},...,c_{k})\to f(c_{1},...,c_{i}\op c'_{i},...,c_{k}).
\] 
These transformations are such that $\delta^{i}=\text{id}$ whenever either $c_{i}$ or $c'_{i}$ is 0, or if any of the other $c_{j}$'s are 0. In addition, these natural transformations are subject to commutativity of a suitable collection of diagrams as in definition \cite[3.2]{Elmendorf}.
\end{definition}

If $\C_{1},...,\C_{k},\D$ are small permutative categories, then $\P(\C_{1},...,\C_{k},\D)$, the set of $k$-morphisms in the multicategory $\P$, is precisely the set of all $k$-linear maps
\[
f:\C_{1}\x\cdots \x\C_{k}\to \D.
\]
The $k$-morphism set $\P(\C_{1},...,\C_{k},\D)$ forms a category. If 
\[
f,g:\C_{1}\x\cdots \x\C_{k}\to \D
\]
are two $k$-linear maps, then a morphism $\alpha:f\to g$ is a natural transformation such that $\alpha(c_{1},...,c_{k})$ is the identity map whenever $c_{i}=0$ for some $i$ and also the following diagram is commutative for each $1\le i\le k$
\[
\xymatrix{
f(c_{1},...,c_{i},...c_{k})\op f(c_{1},...,c'_{i},...c_{k})\ar[r]^-{\delta_{f}^{i}}\ar[d]_-{\alpha\op\alpha}   &f(c_{1},...,c_{i}\op c'_{i},...c_{k})\ar[d]^{\alpha}\\
g(c_{1},...,c_{i},...c_{k})\op g(c_{1},...,c'_{i},...c_{k})\ar[r]_-{\delta_{g}^{i}}         &g(c_{1},...,c_{i}\op c'_{i},...c_{k}).
}
\]

Throughout this paper we will be dealing with functors and multifunctors in a setup where the two concepts make sense (and are different). To avoid confusion we will start by setting up the notation that will be used through this paper. In this section $\M$ and $\mathcal{N}$ will denote general multicategories enriched over $\cat$ and $(\E,\op,0)$ a general permutative category.  Note that as explained before we can see $\E$ as a multicategory and also we can see $\M$ as a category by consider the objects and 1-morphisms and forgetting the rest of the data. 

\noindent{\bf{Notation:}}
\begin{itemize}
	\item We will denote by $\Fi(\E,\M)$ or $\M^{\E}$ the category whose objects are functors $\E\to \M$ and the morphisms are the natural transformation between them.
	\item Also, we will denote by $\Mi(\E,\M)$ the category whose objects are multifunctors $\E\to \M$ and the morphisms the natural transformations preserving the multiproduct.
	\item Finally, by $\Mi_{1}\Fi_{2}(\M,\E;\mathcal{N})$ we mean the category whose objects are assignments $F:\M\x \E\to \mathcal{N}$ that are multifunctors in the first component and are functors in the second component and the morphisms are the transformations preserving the multi-structure in the first component and are natural transformations in the second component.
\end{itemize}

\begin{remark}
As categories, $\Mi(\M,\mathcal{N}^{\E})$ and $\Mi_{1}\Fi_{2}(\M,\E;\mathcal{N})$ are naturally isomorphic. In what follows, we will identify those two categories without further comment.
\end{remark}

We will denote by $\O_{\E}$ and $\O_{\M}$ the object sets of $\E$ and $\M$ respectively. We are interested in studying the functors
\[
\E\to \M
\]
and want to construct a multicategory $\M^{\E}$, whose objects are precisely those functors. The goal of this section is then to prove the following theorem.

\begin{theorem}\label{generalization}
There is a multicategory $\Mo$ whose objects are the functors $\E\to \M$. This multicategory is enriched over $\cat$ if $\M$ is enriched over $\cat$.
\end{theorem}

Before proving this we need some remarks and definitions. First of all, note that given a functor $G:\E\to \M$ and any $k\ge 0$, we can see $G$ as a functor $G:\E^{k}\to \M$ by considering $G(c_{1},...,c_{k})=G(c_{1}\ot \cdots \ot c_{k})$. With this in mind we can give the following definition.

\begin{definition}\label{naturality}
Suppose that $F_{1},...,F_{k},G:\E\to \M$ are functors. A $k$-linear natural transformation from $F_{1},...,F_{k}$ to $G$, is an assignment that for each $k$-tuple $c_{1},..,c_{k}\in \O_{\E}$, corresponds a $k$-morphism 
\[
\phi_{c_{1},...,c_{k}}\in \M(F_{1}(c_{1}),...,F_{k}(c_{k}),G(\ot_{i=1}^{k}c_{i}))
\] 
satisfying the naturality condition described below.
\end{definition}

The naturality condition that we require is the following: suppose that $f_{i}:c_{i}\to c'_{i}$ is a morphisms in $\E$ for $1\le i \le k$, then we have $1$-morphisms $F_{i}(f_{i}):F(c_{i})\to F(c'_{i})$ for $1\le i\le k$ and $G(\ot_{i=1}^{k}f_{i}):G(\ot_{i=1}^{k}c_{i})\to G(\ot_{i=1}^{k}c'_{i})$. The multiproduct $\Gamma_{\M}$ of $\M$ gives maps
\begin{align*}
\Gamma_{\M}:&\M(G(\ot_{i=1}^{k}c_{i});G(\ot_{i=1}^{k}c'_{i}))\x \M(F_{1}(c_{1}),...,F_{k}(c_{k});G(\ot_{i=1}^{k}c_{i}))\\
&\to \M(F_{1}(c_{1}),...,F_{k}(c_{k});G(\ot_{i=1}^{k}c'_{i})),\\
\Gamma_{\M}:&\M(F_{1}(c'_{1}),...,F_{k}(c'_{k});G(\ot_{i=1}^{k}c'_{i}))\x \prod_{i=1}^{k}\M(F_{i}(c_{i});F_{i}(c'_{i}))\\
&\to \M(F_{1}(c_{1}),...,F_{k}(c_{k});G(\ot_{i=1}^{k}c'_{i})).
\end{align*}
We require that
\[
\Gamma_{\M}(G(\ot_{i=1}^{k}f_{i});\phi_{c_{1},...,c_{k}})=\Gamma_{\M}(\phi_{c'_{1},...,c'_{k}};F_{1}(f_{1}),...,F_{k}(f_{k})).
\]
\begin{remark} Note that if we take $k=1$ in the previous definition, then a $1$-linear natural transformation $\phi:F\to G$, is an assignment that for $c\in \O_{\E}$, a $1$-linear map $\phi_{c}:F(c)\to G(c)$ and for each morphism $f:c\to c'$ in $\E$ 
\[
\Gamma_{\M}(G(f);\phi_{c})=\Gamma_{\M}(\phi_{c'};F(f))
\]
This means that $\phi:F\to G$ is a natural transformation such that $\phi_{c}$ is $1$-linear for each object $c$ of $\E$.
\end{remark}

\noindent{\textit{Proof of Theorem \ref{generalization}.}} To prove Theorem \ref{generalization} we need to define a class of objects, a set of $k$-morphisms for $k\ge 0$, a right action of $\Sigma_{k}$ on the set of $k$-morphisms, a unit $1_{F}$ for every object $F$ of $\Mo$ and a multiproduct such that properties Multi(1)-Multi(4) of \cite[Definition 2.1]{Elmendorf} are satisfied. 
\nin
The objects of $\Mo$ are the functors $F:\E\to \M$. To define the $k$-morphisms in $\Mo$, suppose that 
$F_{1},...,F_{k},G:\E\to \M$ are objects in $\Mo$. Then  a $k$-morphism in $\Mo(F_{1},...,F_{k};G)$ is a $k$-linear natural transformation from $F_{1},...,F_{k}$ to $G$. When $k=0$, then a $0$-morphism in $\Mo(,G)$ is a $0$-morphism in $\M(,F(1))$, where $1$ is the coherent unit of $\E$.
\nin
The multicategory $\Mo$ will be an enriched category over $\cat$ whenever $\M$ is enriched over $\cat$. Thus for objects $F_{1},...,F_{k},G$ in $\Mo$, $\Mo(F_{1},...,F_{k},G)$ is the category whose objects are the $k$-linear natural transformations from $F_{1},...,F_{k}$ to $G$. If $\phi,\psi$ are two such transformations, then a morphism $f:\phi\to \psi$ in $\Mo(F_{1},...,F_{k},G)$, is an assignment, for objects $c_{i}$, $1\le i\le k$, a morphism $f_{c_{1},...,c_{k}}:\phi_{c_{1},...,c_{k}}\to \psi_{c_{1},...,c_{k}}$ in the category $\Mo(F_{1}(c_{1}),...,F_{k}(c_{k}),G(\ot_{i=1}^{k}c_{i}))$.
\nin
Let's define the multiproduct on $\Mo$. Suppose that $\phi^{i}\in \Mo(F_{i1},...,F_{ij_{i}};G_{i})$ is a $j_{i}$-morphism for $1\le i\le k$ and that $\psi\in \Mo(G_{1},...,G_{k};H)$ is an $k$-morphism. We want to define $\Gamma_{\Mo}(\psi;\phi^{1},...,\phi^{k})\in \Mo(F_{11},...,F_{kj_{k}};H)$, a $(j_{1}+\cdots+j_{k})$-morphism. To do so, take $c_{11},...,c_{kj_{k}}\in \O_{\E}$ and define
\[
\Gamma_{\Mo}(\psi;\phi^{1},...,\phi^{k})_{c_{11},...,c_{kj_{k}}}=\Gamma_{\M}(\psi_{\ot_{r=1}^{j_{1}}c_{1r},...,\ot_{r=1}^{j_{k}}c_{kr}};\phi^{1}_{c_{11},...,c_{1j_{1}}},...,\phi^{k}_{c_{k1},...,c_{kj_{k}}}).
\]
This definition  readily extends to morphisms in the case that $\M$ is enriched over $\cat$. We need to show that $\Gamma_{\Mo}(\psi;\phi^{1},...,\phi^{k})$ satisfies the naturality condition. Since $\M$ is a multicategory, then by definition we have that 
\[
\Gamma_{\M}(\psi_{\ot_{r=1}^{j_{1}}c_{1r},...,\ot_{r=1}^{j_{k}}c_{kr}};\phi^{1}_{c_{11},...,c_{1j_{1}}},...,\phi^{k}_{c_{k1},...,c_{kj_{k}}})
\] 
is a $(j_{1}+\cdots+j_{k})$-morphism. Suppose then that we have morphisms $f_{rs}:c_{rs}\to c_{rs}'$ in $\E$, for $1\le r\le k$ and $1\le s\le j_{r}$. The multiproduct of $\M$ gives us the maps 
\begin{align*}
\Gamma_{\M}:&\M(H(\ot_{r,s}c_{rs});H(\ot_{r,s}c'_{rs}))\x \M(F_{11}(c_{11}),...,F_{kj_{k}}(c_{kj_{k}});H(\ot_{r,s}c_{rs}))\\
&\to \M(F_{11}(c_{11}),...,F_{kj_{k}}(c_{kj_{k}});H(\ot_{r,s}c'_{rs}))\\
\Gamma_{\M}:&\M(F_{11}(c'_{11}),...,F_{kj_{k}}(c'_{kj_{k}}),H(\ot_{rs}c'_{rs}))\x \prod_{r,s}\M(F_{rs}(c_{rs}),F_{rs}(c'_{rs}))\\
&\to \M(F_{11}(c_{11}),...,F_{kj_{k}}(c_{kj_{k}}),H(\ot_{r,s}c'_{rs})).
\end{align*}
We need to check that 
\begin{align*}
&\Gamma_{\M}(H(\ot_{r,s}f_{rs});\Gamma_{\Mo}(\psi;\phi^{1},...,\phi^{k})_{c_{11},...,c_{kj_{k}}})\\
=&\Gamma_{\M}(\Gamma_{\Mo}(\psi;\phi^{1},...,\phi^{k})_{c'_{11},...,c'_{kj_{k}}};F_{11}(f_{11}),...,F_{kj_{k}}(f_{kj_{k}})).
\end{align*}
On the one hand by definition and property Multi(1) for $\M$ we have
\begin{align*}
&\Gamma_{\M}(H(\ot_{r,s}f_{rs});\Gamma_{\Mo}(\psi;\phi^{1},...,\phi^{k})_{c_{11},...,c_{kj_{k}}})\\
=&\Gamma_{\M}(H(\ot_{r,s}f_{rs});\Gamma_{\M}(\psi_{\ot_{r=1}^{j_{1}}c_{1r},...,\ot_{r=1}^{j_{k}}c_{kr}};\phi^{1}_{c_{11},...,c_{1j_{1}}},...,\phi^{k}_{c_{k1},...,c_{kj_{k}}}))\\
=&\Gamma_{\M}(\Gamma_{\M}(H(\ot_{r,s}f_{rs});\psi_{\ot_{r=1}^{j_{1}}c_{1r},...,\ot_{r=1}^{j_{k}}c_{kr}});\phi^{1}_{c_{11},...,c_{1j_{1}}},...,\phi^{k}_{c_{k1},...,c_{kj_{k}}})
\end{align*}
but as $\psi$ is an $k$-morphism in $\Mo$, we have that 
\begin{align*}
&\Gamma_{\M}(H(\ot_{r,s}f_{rs});\psi_{\ot_{r=1}^{j_{1}}c_{1r},...,\ot_{r=1}^{j_{k}}c_{kr}})=\\
&\Gamma_{\M}(\psi_{\ot_{r=1}^{j_{1}}c'_{1r},...,\ot_{r=1}^{j_{n}}c'_{kr}};G_{1}(\ot_{s}f_{1s}),...,G_{k}(\ot_{s}f_{ks}))
\end{align*}
hence
\begin{align}{\label{1}}
&\Gamma_{\M}(H(\ot_{r,s}f_{rs});\Gamma_{\Mo}(\psi;\phi^{1},...,\phi^{k})_{c_{11},...,c_{kj_{k}}})=\\
&\Gamma_{\M}(\Gamma_{\M}(\psi_{\ot_{r=1}^{j_{1}}c'_{1r},...,\ot_{r=1}^{j_{k}}c'_{kr}};G_{1}(\ot_{s}f_{1s}),...,G_{k}(\ot_{s}f_{ks}));\phi^{1}_{c_{11},...,c_{1j_{1}}},...,\phi^{k}_{c_{k1},...,c_{kj_{k}}}).
\end{align}
On the other hand, again using the definition and property Multi(1) for $\M$ we have
\begin{align*}
&\Gamma_{\M}(\Gamma_{\Mo}(\psi;\phi^{1},...,\phi^{k})_{c'_{11},...,c'_{kj_{k}}};F_{11}(f_{11}),...,F_{kj_{k}}(f_{kj_{k}}))\\
=&\Gamma_{\M}(\Gamma_{\M}(\psi_{\ot_{r=1}^{j_{1}}c'_{1r},...,\ot_{r=1}^{j_{k}}c'_{kr}};\phi^{1}_{c'_{11},...,c'_{1j_{1}}},...,\phi^{k}_{c'_{k1},...,c'_{kj_{k}}});F_{11}(f_{11}),...,F_{kj_{k}}(f_{kj_{k}}))\\
=&\Gamma_{\M}(\psi_{\ot_{r=1}^{j_{1}}c'_{1r},...,\ot_{r=1}^{j_{k}}c'_{kr}};\alpha_{1},...,\alpha_{k}),
\end{align*}
where 
\[
\alpha_{i}=\Gamma_{\M}(\phi^{i}_{c'_{i1},...,c'_{ij_{i}}};F_{i1}(f_{i1}),...,F_{ij_{i}}(f_{ij_{i}})).
\]
But as each $\phi^{i}$ is a $j_{i}$-morphism, by definition we get that 
\[
\alpha_{i}=\Gamma_{\M}(G(\ot_{s}f_{is});\phi^{i}_{c_{i1},...,c_{ij_{i}}})=\Gamma_{\M}(\phi^{i}_{c'_{i1},...,c'_{ij_{i}}};F_{i1}(f_{i1}),...,F_{ij_{i}}(f_{ij_{i}})),
\]
therefore
\begin{align}{\label{2}}
&\Gamma_{\M}(\Gamma_{\Mo}(\psi;\phi^{1},...,\phi^{k})_{c'_{11},...,c'_{kj_{k}}};F_{11}(f_{11}),...,F_{kj_{k}}(f_{kj_{k}}))=\\
&\Gamma_{\M}(\Gamma_{\M}(\psi_{\ot_{r=1}^{j_{1}}c'_{1r},...,\ot_{r=1}^{j_{k}}c'_{kr}};G_{1}(\ot_{s}f_{1s}),...,G_{k}(\ot_{s}f_{ks}));\phi^{1}_{c_{11},...,c_{1j_{1}}},...,\phi^{k}_{c_{k1},...,c_{kj_{k}}})
\end{align}
By (\ref{1}) and (\ref{2}) we see that $\Gamma_{\Mo}(\psi;\phi^{1},...,\phi^{k})$ also satisfies the naturality condition.

We define now a right action of $\Sigma_{k}$ on the set of all $k$-morphisms. Suppose then that $\phi \in \Mo(F_{1},...,F_{k};G)$ is a $k$-morphism in $\Mo$ and $\sigma\in \Sigma_{k}$. Define $\sigma^{*}\phi\in \Mo(F_{\sigma(1)},...,F_{\sigma(k)},G)$ as follows: for objects $c_{1},...,c_{k}$ of $\E$, define 
\[
\sigma^{*}\phi_{c_{\sigma(1)},...,c_{\sigma(k)}}=\sigma^{*}\Gamma(G(\tau_{\sigma,c_{1},...,c_{k}});\phi_{c_{1},...,c_{k}}).
\]
Here 
\[
\tau_{\sigma,c_{1},...,c_{k}}:c_{1}\ot \cdots\ot c_{k}\to c_{\sigma(1)}\ot \cdots \ot c_{\sigma(k)},
\]
is the natural isomorphism in $\E$ and 
\begin{align*}
\sigma^{*}:&\M(F_{1}(c_{1}),...,F_{k}(c_{k}),G(\ot_{i=1}^{k}c_{\sigma(i)}))\to\\
&\M(F_{\sigma(1)}(c_{\sigma(1)}),...,F_{\sigma(k)}(c_{\sigma(k)}),G(\ot_{i=1}^{k}c_{\sigma(i)}))
\end{align*}
is the action of $\Sigma_{k}$ on the $k$-morphisms of the multicategory $\M$. By a similar argument as in the definition of the multiproduct, it follows that this way defined, $\sigma^{*}\phi$ satisfies the naturality condition and thus $\sigma^{*}\phi$ is well defined. This definition extends easily to morphisms in the category of $k$-morphisms when $\M$ is enriched over $\cat$. This indeed defines an action as $\Sigma_{k}$ acts on the right on the set of $k$-morphisms of $\M$. If $F:\E\to \M$ is an object of $\Mo$, then the identity $1_{F}\in \Mo(F;F)$ is the $1$-linear natural transformation $1_{F}$, such that for each object $c$ of $\E$, $(1_{F})_{c}=1_{F(c)}\in \M(F(c);F(c))$. $1_{F}$ clearly satisfies the naturality condition and thus it defines a $1$-morphism in $\Mo(F;F)$ that is clearly a unit.
\nin
Let's check now that this data defines a multicategory; that is, we need to check the conditions Multi(1)-Multi(4) of \cite[Definition 2.1]{Elmendorf}. Conditions Multi(1) and Multi(2) are satisfied pointwise as $\M$ is a multicategory and thus these conditions hold for $\Mo$. Let's check that condition Multi(3) is satisfied. To do so, take $F_{ir},G_{i}$ and $H$ objects in $\Mo$ for $1\le i\le k$ and $1\le r \le j_{i}$. Also take $\sigma\in \Sigma_{k}$. We need to show the commutativity of the following diagram
\[
\xymatrix@C=-5pc{
\Mo(G_{1},...,G_{k};H)\x \prod_{i=1}^{k}\Mo(F_{i1},...,F_{ij_{i}};G_{i})\ar[rd]^{\Gamma}\ar[dd]_{\sigma^{*}\x \sigma^{-1}} &  \\
 &\Mo(F_{11},...,F_{kj_{k}};H)\ar[dd]^-{(\sigma_{<j_{\sigma(1)},...,j_{\sigma(k)}>})^{*}}\\
\Mo(G_{\sigma(1)},...,G_{\sigma(k)};H)\x \prod_{i=1}^{k}\Mo(F_{\sigma(i)1},...,F_{\sigma(i)j_{\sigma(i)}};G_{\sigma(i)})\ar[rd]_{\Gamma} &  \\
      &\Mo(F_{\sigma(1)1},...,F_{\sigma(k)j_{\sigma(k)}};H).
}
\]
To show the commutativity of this diagram, take objects $c_{11},...,c_{kj_{k}}$ of $\E$, $\phi^{i}\in \Mo(F_{i1},...,F_{ij_{i}};G_{i})$  for $1\le i\le k$ and $\psi\in \Mo(G_{1},...,G_{k};H)$ an $k$-morphism. Then by definition we have that 
\begin{align*}
&((\sigma_{<j_{\sigma(1)},...,j_{\sigma(k)}>})^{*}\Gamma(\psi,\phi^{1},...,\phi^{k}))_{c_{\sigma(1)1},...,c_{\sigma(k)j_{\sigma(k)}}}=\\
&(\sigma_{<j_{\sigma(1)},...,j_{\sigma(k)}>})^{*}\Gamma(H(\tau_{\sigma,c_{11},...,c_{kj_{k}}});\Gamma(\psi;\phi^{1},...,\phi^{k})_{c_{11},...,c_{kj_{k}}})=\\
&(\sigma_{<j_{\sigma(1)},...,j_{\sigma(k)}>})^{*}\Gamma(H(\tau_{\sigma,c_{11},...,c_{kj_{k}}});\Gamma(\psi_{\ot_{i=1}^{j_{1}}c_{1i},...,\ot_{i=1}^{j_{k}}c_{ki}};\phi^{1}_{c_{11},...,c_{1j_{1}}},...,\phi^{1}_{c_{k1},...,c_{kj_{k}}})).
\end{align*}
On the other hand
\begin{align*}
&\Gamma((\sigma^{*}\x \sigma^{-1})(\psi;\phi^{1},...,\phi^{k}))_{c_{\sigma(1)1},...,c_{\sigma(k)j_{\sigma(k)}}} =\\
&\Gamma(\sigma^{*}\psi;\phi^{\sigma(1)},...,\phi^{\sigma(k)})_{c_{\sigma(1)1},...,c_{\sigma(k)j_{\sigma(k)}}}=\\
&\Gamma(\sigma^{*}\psi_{\ot_{i=1}^{j_{\sigma(1)}}c_{\sigma(1)i},...,\ot_{i=1}^{j_{\sigma(k)}}c_{\sigma(k)i}}; \phi^{\sigma(1)}_{c_{\sigma(1)1},...,c_{\sigma(1)j_{\sigma(1)}}},...,\phi^{\sigma(k)}_{c_{\sigma(k)1},...,c_{\sigma(k)j_{\sigma(k)}}})=\\
&\Gamma(\sigma^{*}\Gamma(H(\tau_{\sigma,c_{11}},...,c_{kj_{k}});\psi_{\ot_{i=1}^{j_{1}}c_{1i},...,\ot_{i=1}^{j_{k}}c_{ki}});\phi^{\sigma(1)}_{c_{\sigma(1)1},...,c_{\sigma(1)j_{\sigma(1)}}},...,\phi^{\sigma(k)}_{c_{\sigma(k)1},...,c_{\sigma(k)j_{\sigma(k)}}}).
\end{align*}
Since $\M$ is a multicategory, by property Multi(3) we get
\begin{align*}
&\Gamma(\sigma^{*}\Gamma(H(\tau_{\sigma,c_{11}},...,c_{kj_{k}});\psi_{\ot c_{1i},...,\ot c_{ki}});\phi^{\sigma(1)}_{c_{\sigma(1)1},...,c_{\sigma(1)j_{\sigma(1)}}},...,\phi^{\sigma(k)}_{c_{\sigma(k)1},...,c_{\sigma(k)j_{\sigma(k)}}})\\
&=\Gamma((\sigma^{*}\x\sigma^{-1})(\Gamma(H(\tau_{\sigma,c_{11}},...,c_{kj_{k}};\psi_{\ot c_{1i},...,\ot c_{ki}});\phi^{1}_{c_{11},...,c_{1j_{1}}},...,\phi^{1}_{c_{k1},...,c_{kj_{k}}})))\\
&=(\sigma_{<j_{\sigma(1)},...,j_{\sigma(k)}>})^{*}\Gamma(\Gamma(H(\tau_{\sigma,c_{11}},...,c_{kj_{k}};\psi_{\ot c_{1i},...,\ot c_{ki}});\phi^{1}_{c_{11},...,c_{1j_{1}}},...,\phi^{1}_{c_{k1},...,c_{kj_{k}}})),
\end{align*}
but again using the fact that $\M$ is a multicategory, by property Multi(1) we get that 
\begin{align*}
&\Gamma(\Gamma(H(\tau_{\sigma,c_{11}},...,c_{kj_{k}};\psi_{\ot_{i=1}^{j_{1}}c_{1i},...,\ot_{i=1}^{j_{k}}c_{ki}});\phi^{1}_{c_{11},...,c_{1j_{1}}},...,\phi^{1}_{c_{k1},...,c_{kj_{k}}}))\\
&=\Gamma(H(\tau_{\sigma,c_{11},...,c_{kj_{k}}});\Gamma(\psi_{\ot_{i=1}^{j_{1}}c_{1i},...,\ot_{i=1}^{j_{k}}c_{ki}};\phi^{1}_{c_{11},...,c_{1j_{1}}},...,\phi^{1}_{c_{k1},...,c_{kj_{k}}})).
\end{align*}
Thus diagram Multi(3) commutes. In the case that $\M$ is enriched over $\cat$, a similar argument shows that the respective diagram commutes on morphisms. The commutativity of diagram Multi(4) is similar. This proves Theorem \ref{generalization}.
\qed

\noindent{\bf{Example: }}
Our main application of the previous theorem is the case where $\M=\P$. Explicitly, the multicategory $\Po$, has objects the functors
\[
\E\to \P.
\]
If $F_{1},...,F_{k}$ and $G$ are objects in $\Po$, then the $k$-morphism set, $\Po(F_{1},...,F_{k};G)$ forms a category, its objects are the $k$-linear natural maps; that is, the assignments $\phi$ that for each $c_{1},...,c_{k}\in \O_{\E}$ correspond a $k$-linear map
\[
\phi_{c_{1,...,c_{k}}}:F_{1}(c_{1})\x \cdots \x F_{k}(c_{k})\to G(c_{1}\ot \dots \ot c_{k})
\]
such that $\phi$ satisfies the naturality condition as described before. If $\phi$ and $\psi$ are two objects in $\Po(F_{1},...,F_{k};G)$, then a morphism $\alpha:\phi\to \psi$ is a morphism of natural $k$-linear transformations; that is, for $c_{1},...,c_{k}$
\[
\alpha_{c_{1},...,c_{k}}:\phi_{c_{1},...,c_{k}}\to \psi_{c_{1},...,c_{k}}
\]
is a natural map such that the following diagram is commutative
\[
\xymatrix{
\phi_{c_{1},...,c_{k}}(x_{1},...,x_{i},...,x_{k})\op \phi_{c_{1},...,c_{k}}(x_{1},...,x'_{i},...,x_{k}) \ar[r]^-{\delta^{i}_{\phi_{c_{1},...,c_{k}}}}\ar[d]_-{\alpha_{c_{1},...,c_{k}}\op \alpha_{c_{1},...,c_{k}}}   &\phi_{c_{1},...,c_{k}}(x_{1},...,x_{i}\op x'_{i},...,x_{k}) \ar[d]^-{\alpha_{c_{1},...,c_{k}}}\\
\psi_{c_{1},...,c_{k}}(x_{1},...,x_{i},...,x_{k})\op \psi_{c_{1},...,c_{k}}(x_{1},...,x'_{i},...,x_{k})\ar[r]_-{\delta^{i}_{\psi_{c_{1},...,c_{k}}}}         &\psi_{c_{1},...,c_{k}}(x_{1},...,x_{i}\op x'_{i},...,x_{k}).
}
\]
Here, $\delta^{i}_{\phi_{c_{1},...,c_{k}}}$ and $\delta^{i}_{\psi_{c_{1},...,c_{k}}}$ are the structural maps of the $k$-linear maps $\phi_{c_{1},...,c_{k}}$ and $\psi_{c_{1},...,c_{k}}$.
\nin
The action of $\Sigma_{k}$ on $k$-morphisms, is defined in the following way. If $\sigma\in \Sigma_{k}$ and $\phi\in \Po(F_{1},...,F_{k},G)$. For $c_{1},...,c_{k}\in \O_{C}$, define $\sigma^{*}\phi_{c_{\sigma(1)},...,c_{\sigma(k)}}$ is the $k$-linear map making the following diagram commutative
\[
\xymatrix{
\prod_{i}^{k}F_{\sigma(i)}(c_{\sigma(i)})\ \ \ \ \ \ar[r]^-{\sigma^{*}\phi_{c_{\sigma(1)},...,c_{\sigma(k)}}}\ar[d]_-{\sigma^{-1}}   &\ \ \ \ \ \ \ G(\ot_{i=1}^{k}c_{\sigma(i)})\\
\prod_{i}^{k}F_{i}(c_{i})\ \ \ \ar[r]_-{\phi_{c_{1},...,c_{k}}}         &\ \ \ \ G(\ot_{i=1}^{k}c_{i})\ar[u]_-{G(\tau_{\sigma,c_{1},...,c_{k}}).}
}
\]
Here $\tau_{\sigma,c_{1},...,c_{k}}$ is the coherent isomorphism $c_{1}\ot \cdots \ot c_{k}\stackrel{\approx}{\rightarrow} c_{\sigma(1)}\ot \cdots \ot c_{\sigma(k)}$ in $\E$ obtained by iterated applications of $\gamma$. On the other hand, the multiproduct in $\Po$ is defined as follows: suppose $\phi^{i}\in\Po(F_{i1},...,F_{ij_{i}};G_{i})$, $\psi\in \Po(G_{1},...,G_{n};H)$. For any objects $c_{is}$ of $\E$, for $1\le i\le n$, $1\le s\le j_{i}$, define
\[
\Gamma(\psi,\phi^{1},...,\phi^{n})_{c_{11},...,c_{nj_{n}}}=\psi_{c_{11}\ot\cdots \ot c_{1j_{1}},...,c_{n1}\ot\cdots \ot c_{nj_{n}}}(\phi^{1}_{c_{11},...,c_{1j_{1}}},...,\phi^{n}_{c_{n1},...,c_{nj_{n}}}).
\]
Suppose now that $\M_{1}$ and $\M_{2}$ are two multicategories and that $\F:\M_{1}\to \M_{2}$ is a multifunctor. Then $\F$ induces a functor $\F_{*}:\Mo_{1}\to \Mo_{2}$ defined by composition with $\F$. This functor preserves the multicategory structure; that is, $\F_{*}$ is a multifunctor. In the case where that both $\M_{1}$ and $\M_{2}$ are enriched over $\cat$ and that $\F$ is an enriched multifunctor, then $\F_{*}$ is also enriched. We prove this in the following theorem.

\begin{theorem}
The functor $\F_{*}$ is a multifunctor that is enriched if $\M_{1}$ and $\M_{2}$ are enriched over $\cat$ and $\F$ is an enriched multifunctor.
\end{theorem}
\Proof
Suppose that $F_{1},...,F_{k},G:\E\to \M_{1}$ are objects in $\Mo_{1}$, then $\F_{*}(F_{1})=\F\circ F_{1},...,\F_{*}(G)=\F\circ G:\E\to \M_{2}$. Take $\phi\in \Mo_{1}(F_{1},...,F_{k},G)$. Thus given objects $c_{1},...,c_{k}$ in $\E$, by definition $\phi_{c_{1},...,c_{k}}\in \M_{1}(F_{1}(c_{1}),...,F_{k}(c_{k}),G(\ot_{i=1}^{k}c_{i}))$ and thus $\F(\phi_{c_{1},...,c_{k}})\in \M_{2}(\F(F_{1})(c_{1}),...,\F(F_{k})(c_{k}),\F(G)(\ot_{i=1}^{k}c_{i}))$. We define 
\[
\F_{*}(\phi)\in \Mo_{2}(\F_{*}(F_{1}),...,\F_{*}(F_{k});\F_{*}(G))
\]
by
\[
\F_{*}(\phi)_{c_{1},...,c_{k}}=\F(\phi_{c_{1},...,c_{k}})\in \M_{2}(\F_{*}(F_{1})(c_{1}),...,\F_{*}(F_{k})(c_{k});\F_{*}(G)(\ot_{i=1}^{k}c_{i})).
\]  
We need to show that $\F_{*}(\phi)$ satisfies the naturality condition; that is, if we are given morphisms $f_{i}:c_{i}\to c'_{i}$ in $E$, we need to show that
\[
\Gamma(\F_{*}(G)(\ot f_{i});\F_{*}(\phi)_{c_{1},...,c_{k}})=\Gamma(\F_{*}(\phi)_{c'_{1},...,c'_{k}};\F_{*}(F_{1})(f_{1}),...,\F_{*}(F_{k})(f_{k})).
\]
But we know that this is true for $\phi$, thus we know that
\[
\Gamma(G(\ot f_{i});\phi_{c_{1},...,c_{k}})=\Gamma(\phi_{c'_{1},...,c'_{k}};F_{1}(f_{1}),...,F_{k}(f_{k})).
\]
Applying $\F$ to both sides and using that $\F$ is a multifunctor we see that
\begin{align*}
\Gamma(\F_{*}(G)(\ot f_{i});\F_{*}(\phi)_{c_{1},...,c_{k}})&=\F(\Gamma(G(\ot f_{i});\phi_{c_{1},...,c_{k}}))\\
\F(\Gamma(\phi_{c'_{1},...,c'_{k}};F_{1}(f_{1}),...,F_{k}(f_{k})))&=\Gamma(\F_{*}(\phi)_{c'_{1},...,c'_{k}};\F_{*}(F_{1})(f_{1}),...,\F_{*}(F_{k})(f_{k})),
\end{align*}
which proves that $\F_{*}(\phi)$ also satisfies the naturality condition.

Let's show now that $\F_{*}$ respects the multiproduct. To do so consider $F_{ir},G_{i}$ and $H$ objects in $\Mo_{1}$ for $1\le i\le k$ and $1\le r\le j_{i}$. Also  consider $\phi^{i}\in \Mo_{1}(F_{i1},...,F_{ij_{i}},G_{i})$ and $\psi \in \Mo_{1}(G_{1},...,G_{k},H)$ for $1\le i\le k$ and $c_{ir}$ objects in $\E$. Then by definition
\begin{align*}
&(\F_{*}\Gamma(\psi;\phi^{1},...,\phi^{k}))_{c_{11},...,c_{kj_{k}}}=\F_{*}\Gamma(\psi_{\ot_{r=1}^{j_{1}}c_{1r},...,\ot_{r=1}^{j_{k}}c_{kr}};\phi^{1}_{c_{11},...,c_{1j_{1}}},...,\phi^{k}_{c_{k1},...,c_{kj_{k}}})\\
=&\Gamma(\F_{*}(\psi)_{\ot_{r=1}^{j_{1}}c_{1r},...,\ot_{r=1}^{j_{k}}c_{kr}};\F_{*}(\phi^{1})_{c_{11},...,c_{1j_{1}}},...,\F_{*}(\phi^{k})_{c_{k1},...,c_{kj_{k}}})
\\
=&(\Gamma(\F_{*}(\psi);\F_{*}(\phi^{1}),...,\F_{*}(\phi^{k})))_{c_{11},...,c_{kj_{k}}},
\end{align*}
and thus $\F_{*}$ preserves the multiproduct. On the other hand, it's trivial to show that $\F_{*}$ preserves units. We still need to show that $\F_{*}$ is compatible with the action of $\Sigma_{k}$. Thus take $\sigma\in \Sigma_{k}$. We need to show that the following diagram is commutative
\[
\xymatrix{
\Mo_{1}(F_{1},...,F_{k};G)\ar[rr]^{\F_{*}\ \ \ \ }\ar[d]_{\sigma^{*}}  & &\Mo_{2}(\F_{*}(F_{1}),...,\F_{*}(F_{k});\F_{*}(G)) \ar[d]^{\sigma^{*}}\\
\Mo_{1}(F_{\sigma(1)},...,F_{\sigma(k)};G)\ar[rr]_{\F_{*}\ \ \ \ }       & &\Mo_{2}(\F_{*}(F_{\sigma(1)}),...,\F_{*}(F_{\sigma(k)});\F_{*}(G)).
}
\]
But by definition and the fact that $\F$ is a multifunctor we have that 
\begin{align*}
&\F_{*}(\sigma^{*}\phi)_{c_{\sigma_{1}},...,c_{\sigma_{k}}}=\F(\sigma^{*}\phi_{c_{\sigma_{1}},...,c_{\sigma_{k}}})=\F(\sigma^{*}\Gamma(G(\tau_{\sigma,c_{1},...,c_{k}});\phi_{c_{1},...,c_{k}}))\\
=&\sigma^{*}(\F_{*}(\Gamma(G(\tau_{\sigma,c_{1},...,c_{k}});\phi_{c_{1},...,c_{k}})))=\sigma^{*}(\Gamma(\F_{*}(G)(\tau_{\sigma,c_{1},...,c_{k}});\F_{*}(\phi_{c_{1},...,c_{k}})))\\
=&\sigma^{*}(\F_{*}(\phi))_{c_{\sigma(1)},...,c_{\sigma(k)}}.
\end{align*}
This proves that $\F_{*}$ is a multifunctor. In the case that $\M_{1}$, $\M_{2}$ are enriched over $\cat$ and $\F$ and enriched multifunctor, then by a similar argument we see that $\F_{*}$ is a multifunctor.
\qed

In \cite{Elmendorf} Elmendorff and Mandell constructed a multifunctor $K:\P\to \Sym$. Using this multifunctor and the previous theorem we get the following corollary.

\begin{corollary} The multifunctor $K:\P\to \Sym$ induces an enriched multifunctor $K_{*}:\Po\to \Sym^{\E}$
\end{corollary}

\section{Multifunctors to $\Po$}\label{multiP}

In this section we study the multicategory $\Po$; that is, the multicategory obtained by taking $\M=\P$ in the previous construction. In particular we characterize the multifunctors to $\Po$ out of certain parameters multicategories $\Sigma_{*}$ and $E\Sigma_{*}$ that are defined in this section. In \cite{Elmendorf} Elmendorf and Mandell proved that multifunctors $\Sigma_{*},E\Sigma_{*}\to \Per$ are determined by additional structure on the image of the only object of $\Sigma_{*}$ and $E\Sigma_{*}$ respectively. Following their idea, we show that multifunctors from $\Sigma_{*}$ and $E\Sigma_{*}$ to $\Po$ are determined by additional structure on the fibered category associated to the functor $\E\to \P$ that is the image of the only object of $\Sigma_{*}$ and $E\Sigma_{*}$.

We will denote by $\Sigma_{*}$ the associative operad; that is, $\Sigma_{*}$ is the operad whose value at $k\ge 0$ is the symmetric group $\Sigma_{k}$. This operad is relevant as its algebras on a symmetric monoidal category are the associative monoids. We regard $\Sigma_{*}$ as a multicategory with only one object. The goal of this section is to give an explicit description to the multifunctors 
\[
S:\Sigma_{*}\to \Po.
\]
Such a multifunctor maps the unique object of $\Sigma_{*}$ to an object $F$ of $\Po$. The multifunctor $\Sigma_{*}\to \Po$, then factors through the multicategory generated by $F$ in $\Po$. Suppose we have such a multifunctor $\Sigma_{*}\to \Po$ mapping the only object of $\Sigma_{*}$ to $F$. Consider $1_{2}\in \Sigma_{2}$ the identity element. This is a $2$-morphism in the multicategory $\Sigma_{*}$, hence it is mapped to a $2$-morphism $\ot\in \Po(F,F;F)$. This operation is strictly associative,  there is a unit $1$ that is an object of $F(1_{\E})$. The operation $\ot$ and the unit $1$ satisfy similar coherences as in \cite[Definition 3.3]{Elmendorf}. To be more precise, we have that $\ot\in \Po(F,F;F)$, this means that for every pair $c_{1},c_{2}\in \O_{\E}$ we have a $2$-linear map
\[
\oti_{c_{1},c_{2}}:F(c_{1})\x F(c_{2})\to F(c_{1}\ot c_{2}) 
\]
Each $\oti_{c_{1},c_{2}}$ comes equipped with distributivity maps
\begin{align*}
d_{c_{1},c_{2}}^{l}&:(x\oti_{c_{1},c_{2}}y)\op (x'\oti_{c_{1},c_{2}}y)\to (x\op x')\oti_{c_{1},c_{2}}y, \\
d_{c_{1},c_{2}}^{r}&:(x\oti_{c_{1},c_{2}}y)\op (x\oti_{c_{1},c_{2}}y')\to x\oti_{c_{1},c_{2}}(y\op y') 
\end{align*} 
that are the identity whenever $x$ or $x'$ equal $0_{c_{1}}$ or whenever $y$ or $y'$ equal $0_{c_{2}}$.
\nin  
By the definition of a $2$-linear map the following coherences are satisfied.
\[
\tag{\textbf{c.1}}
\xymatrix{
(x\oti_{c_{1},c_{2}}y)\op (x'\oti_{c_{1},c_{2}}y)\op (x''\oti_{c_{1},c_{2}}y)\ar[r]^{\ \ \ \  \ d^{l}_{c_{1},c_{2}}\op \text{id}}  \ar[d]_{\text{id}\op d^{l}_{c_{1},c_{2}}}           &\ \ \ \ ((x\op x')\oti_{c_{1},c_{2}}y)\op (x''\oti_{c_{1},c_{2}}y)\ar[d]^{d^{l}_{c_{1},c_{2}}}\\
(x\oti_{c_{1},c_{2}}y )\op ((x'\op x'')\oti_{c_{1},c_{2}}y)\ar[r]_{\ \ \ \ \ d^{l}_{c_{1},c_{2}}}       &(x\op x'\op x'')\oti_{c_{1},c_{2}}y
}
\]
is a commutative diagram and there is a similar commutative diagram for $d_{c_{1},c_{2}}^{r}$.
\nin
Also
\[
\tag{\textbf{c.2}}
\xymatrix{
(x\oti_{c_{1},c_{2}}y)\op (x'\oti_{c_{1},c_{2}}y)\ar[r]^{\ \ \ \ \ d^{l}_{c_{1},c_{2}}}\ar[d]_{\gamma^{\op}}   &(x\op x')\oti_{c_{1},c_{2}}y\ar[d]^{\gamma^{\op}\ot_{c_{1},c_{2}}\text{id}}\\
(x'\oti_{c_{1},c_{2}}y)\op (x\oti_{c_{1},c_{2}}y)\ar[r]_{\ \ \ \ \ d^{l}_{c_{1},c_{2}}}         &(x'\op x)\oti_{c_{1},c_{2}}y
}
\]
is a commutative diagram and there is a similar commutative diagram for $d_{c_{1},c_{2}}^{r}$.
\nin
We also need the following diagram to commute
\[
\tag{\textbf{c.3}}
\xymatrix@C=-7pc{&(x\oti_{c_{1},c_{2}}(y\op y'))\op (x'\oti_{c_{1},c_{2}}(y\op y'))
\ar[ddr]^-{d_{c_{1},c_{2}}^{l}}
\\
(x\oti_{c_{1},c_{2}}y)\op (x\oti_{c_{1},c_{2}}y')\op (x'\oti_{c_{1},c_{2}}y)\op (x'\oti_{c_{1},c_{2}}y')
\ar[ur]^-{d_{c_{1},c_{2}}^{r}\op d_{c_{1},c_{2}}^{r}\ \ \ }
\ar[dd]_-{\text{id}\op \gamma^{\op}\op \text{id}}
\\
&&(x\op x')\oti_{c_{1},c_{2}}(y\op y').
\\
(x\oti_{c_{1},c_{2}}y)\op (x'\oti_{c_{1},c_{2}}y)\op (x\oti_{c_{1},c_{2}}y')\op (x'\oti_{c_{1},c_{2}}y')
\ar[dr]_{d_{c_{1},c_{2}}^{l}\op d_{c_{1},c_{2}}^{l}}
\\
&((x\op x')\oti_{c_{1},c_{2}}y)\op ((x\op x')\oti_{c_{1},c_{2}}y')
\ar[uur]_-{d_{c_{1},c_{2}}^{r}}}
\]
In addition, we require that $\ot$ vanishes whenever one of the inputs is zero, more concretely, we require that

\begin{align*}\tag{\textbf{c.4}}
0_{c_{1}}\oti_{c_{1},c_{2}}x&=0_{c_{1}\otimes c_{2}},\\
x\oti_{c_{1},c_{2}}0_{c_{2}}&=0_{c_{1}\otimes c_{2}}.
\end{align*}

Diagrams $(\textbf{c.1}),(\textbf{c.2}),(\textbf{c.3})$ and $(\textbf{c.4})$ say that each $\ot_{c_{1},c_{2}}$ is a $2$-linear map. In order for $\ot$ to be a $2$-morphism in the multicategory $\Po$ we need the naturality condition to be satisfied; that is, given morphisms $f_{i}:c_{i}\to c'_{i}$ in $\E$ for $i=1$ and $i=2$, then the following diagram is commutative
\[
\tag{\textbf{c.5}}
\xymatrix{
 F(c_{1})\x F(c_{2})\ar[rr]^{F(f_{1})\x F(f_{2})}\ar[d]_{\oti_{c_{1},c_{2}}}        &                &F(c'_{1})\x F(c'_{2})\ar[d]^{\oti_{c'_{1},c'_{2}}}\\
 F(c_{1}\ot c_{2} )\ar[rr]_{F(f_{1}\ot f_{2})}				&							    &F(c'_{1}\ot c'_{2} ).
}
\]
We need this diagram to be a commutative diagram of $2$-linear maps; that is, given $x$ and $x'$ objects in $F(c_{1})$, $y$ object in $F(c_{2})$ we need the maps $(*)$ and $(**)$ to be the same
\begin{align*}\tag{*}
&F(f_{1})x\oti_{c'_{1},c'_{2}}F(f_{2})y\op F(f_{1})x'\oti_{c'_{1},c'_{2}}F(f_{2})y\stackrel{d^{l}_{c'_{1},c'_{2}}}{\rightarrow}F(f_{1})x\op F(f_{1})x'\oti_{c'_{1},c'_{2}}F(f_{2})y\\
&\stackrel{\lambda_{F(f_{1})}}{\rightarrow} (F(f_{1})(x\op x'))\oti_{c'_{1},c'_{2}}F(f_{2})y,
\end{align*}
\begin{align*}\tag{**}
&F(f_{1}\ot f_{2})(x\oti_{c_{1},c_{2}}y)\op F(f_{1}\ot f_{2})(x'\oti_{c_{1},c_{2}}y)\stackrel{\lambda_{F(f_{1}\ot f_{1})}}{\rightarrow}\\
&F(f_{1}\ot f_{2})(x\oti_{c_{1},c_{2}}y\op x'\oti_{c_{1},c_{2}}y)\stackrel{F(f_{1}\ot f_{2})(d^{l}_{c_{1},c_{2}})}{\rightarrow} F(f_{1}\ot f_{2})((x\op x')\oti_{c_{1},c_{2}}y);
\end{align*}
that is, we have that
\[
\tag{\textbf{c.6}}
\lambda_{F(f_{1})}\circ d^{l}_{c'_{1},c'_{2}}=F(f_{1}\ot f_{2})(d^{l}_{c_{1},c_{2}})\circ \lambda_{F(f_{1}\ot f_{2})},
\]
and a similar condition is satisfied for $d^{r}_{c_{1},c_{2}}$.
\nin
Also we have that $1$, which is an object of $F(1_{\E})$, is a strict unit, in this case this means that
\[
\tag{\textbf{c.7}}
1\oti_{1,c}x=x=x\oti_{c,1}1.
\]
We also need strict associativity, which in this case this means that for every objects $c_{1},c_{2}$ and $c_{3}$ of $\E$, and every objects $x,y$ and $z$ of $F(c_{1}),F(c_{2})$ and $F(c_{3})$ respectively, we have
\[
\tag{\textbf{c.8}}
(x\oti_{c_{1},c_{2}}y)\oti_{c_{1}\ot c_{2},c_{3}}z= x\oti_{c_{1},c_{2}\ot c_{3}}(y\oti_{c_{2},c_{3}}z).
\]
In addition we need the following diagram (\textbf{c.9}) to commute and a similar for $d^{r}$
\[
\tag{\textbf{c.9}}
\xymatrix@C=-8pc{&(x\oti_{c_{1},c_{2}\ot c_{3}}(y\oti_{c_{2},c_{3}}z))\op (x'\oti_{c_{1},c_{2}\ot c_{3}}(y\oti_{c_{2},c_{3}}z))
\ar[ddr]^-{d^{l}_{c_{1},c_{2}\ot c_{3}}}
\\
((x\oti_{c_{1},c_{2}}y )\oti_{c_{1}\ot c_{2},c_{3}}z)\op  ((x'\oti_{c_{1},c_{2}}y)\oti_{c_{1}\ot c_{2},c_{3}}z) 
\ar[ur]^-{=}
\ar[dd]_-{d_{c_{1}\ot c_{2},c_{3}}^{l}}
\\
&&(x\op x')\oti_{c_{1},c_{2}\ot c_{3}}(y\oti_{c_{2},c_{3}}z).
\\
((x\oti_{c_{1},c_{2}}y )\op (x'\oti_{c_{1},c_{2}}y))\oti_{c_{1}\ot c_{2}, c_{3}}z
\ar[dr]_{d_{c_{1},c_{2}}^{l}\ot_{c_{1},c_{2}\ot c_{3}} \text{id}\ \ \ \ \ \ }
\\
&((x'\op x)\oti_{c_{1},c_{2}}y)\oti_{c_{1}\ot c_{2},c_{3}}z
\ar[uur]_-{=}}
\]
Finally we need the diagram (\textbf{c.10}) to be commutative.
\[
\tag{\textbf{c.10}}
\xymatrix{
((x\oti_{c_{1},c_{2}}y)\oti_{c_{1}\ot c_{2},c_{3}}z)\op ((x\oti_{c_{1},c_{2}}y')\oti_{c_{1}\ot c_{2},c_{3}}z)\ar[r]^-{d^{l}_{c_{1}\ot c_{2},c_{3}}}\ar[d]_-{=}   &((x\oti_{c_{1},c_{2}}y)\op (x\oti_{c_{1},c_{2}}y))\oti_{c_{1}\ot c_{2},c_{3}}z\ar[d]^-{d^{r}_{c_{1},c_{2}}\ot_{c_{1}\ot c_{2},c_{3}}\text{id}}\\
(x\oti_{c_{1},c_{2}\ot c_{3}}(y\oti_{c_{2},c_{3}}z))\op (x\oti_{c_{1},c_{2}\ot c_{3}}(y'\ot_{c_{2},c_{3}}z))\ar[d]_-{d^{r}_{c_{1},c_{2}\ot c_{3}}}         &(x\oti_{c_{1},c_{2}}(y\op y'))\oti_{c_{1}\ot c_{2},c_{3}}z\ar[d]^-{=}\\
x\oti_{c_{1},c_{2}\ot c_{3}}((y\oti_{c_{2},c_{3}}z)\op (y'\oti_{c_{2},c_{3}}z))\ar[r]_-{\text{id} \ot_{c_{1},c_{2}\ot c_{3}} d^{l}_{c_{2},c_{3}}}					&x\oti_{c_{1},c_{2}\ot c_{3}}((y\op y')\oti_{c_{2},c_{3}}z).
}
\]
We have the following theorem.

\begin{theorem}\label{multisigma}
Having a multifunctor $\Sigma_{*}\to \Po$ mapping the only object of $\Sigma_{*}$ to a functor $F:\E\to \P$ is equivalent to having a functor $F:\E\to \P$, a functor 
\[
\oti_{c_{1},c_{2}}:F(c_{1})\x F(c_{2})\to F(c_{1}\ot c_{2}) 
\]
for each pair $c_{1},c_{2}\in \O_{\E}$, a unit $1$ which is an object in $F(1_{\E})$ and distributivity maps 
\begin{align*}
d_{c_{1},c_{2}}^{l}&:(x\oti_{c_{1},c_{2}}y)\op (x'\oti_{c_{1},c_{2}}y)\to (x\op x')\oti_{c_{1},c_{2}}y, \\
d_{c_{1},c_{2}}^{r}&:(x\oti_{c_{1},c_{2}}y)\op (x\oti_{c_{1},c_{2}}y')\to x\oti_{c_{1},c_{2}}(y\op y') 
\end{align*} 
satisfying conditions $(\textbf{c.1})-(\textbf{c.10})$.
\end{theorem}
\Proof
First suppose that we have a multifunctor $S:\Sigma_{*}\to \Po$. Then $S$ factors through the multicategory generated by $F:\E\to \P$, for some $F$. Thus, we may regard  $S$ as a multifunctor $S:\Sigma_{*}\to \{\Po(F,...,F;F)\}$, sending the only object of $\Sigma_{*}$ to $F$. Let $1_{n}\in \Sigma_{n}$ be the identity element. Then $1_{2}$ is a 2-morphism and we let $\ot=S_{2}(1_{2})\in \Po(F,F;F)$. Unraveling the definitions, we see that $\ot_{c_{1},c_{2}}$ is $2$-linear map for every $c_{1},c_{2}\in \O_{\E}$. Therefore conditions $(\textbf{c.1})-(\textbf{c.4})$ hold. In addition, properties $(\textbf{c.5})$ and $(\textbf{c.6})$ hold as $\ot$ satisfies the naturality condition, this follows from the fact that $\ot$ is a $2$-morphism in the multicategory $\Po$. Consider $1_{0}\in \Sigma_{0}$. This element is mapped under $S_{0}$ to a $0$-morphism in $\Po(;F)$. By definition, this is just an object in the category $F(1_{\E})$ which we define as the unit $1$. Condition $(\textbf{c.7})$ follows from the fact that, in $\Sigma_{*}$, we have the equality $\Gamma_{\Sigma_{*}}(1_{2};1_{1},1_{0})=1_{1}=\Gamma_{\Sigma_{*}}(1_{2};1_{0},1_{1})$.  Property $(\textbf{c.8})$, the associativity property, follows from the fact that in the multicategory $\Sigma_{*}$, we have $\Gamma_{\Sigma_{*}}(1_{2};1_{2},1_{1})=\Gamma_{\Sigma_{*}}(1_{2};1_{1},1_{2})$, thus we have that 
\[
\Gamma_{\Po}(\ot;\ot,1_{F})=\Gamma_{\Po}(\ot;1_{F},\ot).
\] 
The previous equation implies condition $(\textbf{c.8})$. Note that both $\Gamma_{\Po}(\ot;\ot,1_{F})$ and $\Gamma_{\Po}(\ot;1_{F},\ot)$ are $3$-morphisms in $\Po$. Therefore, for every triple $c_{1},c_{2},c_{3}\in \O_{\E}$, 
\[
\Gamma_{\Po}(\ot;\ot,1_{F})_{c_{1},c_{2},c_{3}}\ \ \text{ and }\ \ \Gamma_{\Po}(\ot;1_{F},\ot)_{c_{1},c_{2},c_{3}},
\]
are $3$-linear maps that agree, in particular their structural maps $\delta^{i}$ must agree. This shows that conditions $(\textbf{c.9})$ and $(\textbf{c.10})$ are satisfied.  

On the other hand, suppose we have a functor $F:\E\to \P$, a collection of $2$-linear maps $\ot_{c_{1},c_{2}}$ for each pair $c_{1}$ and $c_{2}\in \O_{\E}$ and a unit $1$ that satisfy properties $(\textbf{c.1})-(\textbf{c.10})$. We want to define a multifunctor $S:\Sigma_{*}\to \Po$. To begin, we define $S$ evaluated on the only object of $\Sigma_{*}$ as $F$. 

Note that conditions $(\textbf{c.1})-(\textbf{c.4})$ say precisely that $\ot_{c_{1},c_{2}}$ is a $2$-linear map 
\[
\oti_{c_{1},c_{2}}:F(c_{1})\x F(c_{2})\to F(c_{1}\ot c_{2}),
\]
with structural maps $d^{l}_{c_{1},c_{2}}$ and $d^{r}_{c_{1},c_{2}}$ for every pair $c_{1}$ and $c_{2}\in \O_{\E}$. In addition, by properties $(\textbf{c.5})$ and $(\textbf{c.6})$ we see that $\ot$ satisfies the naturality condition, therefore $\ot$ is a $2$-morphism in $\Po(F,F;F)$. We claim that properties $(\textbf{c.7})-(\textbf{c.8})$ imply that, as $3$-morphisms in $\Po(F,F,F;F)$, 
\begin{equation}\label{3morph}
\Gamma_{\Po}(\ot;\ot,1_{F})=\Gamma_{\Po}(\ot;1_{F},\ot).
\end{equation}
Indeed, for $c_{1},c_{2}$ and $c_{3}\in \O_{\E}$ we have
\begin{align*}
\Gamma_{\Po}(\ot;\ot,1_{F})_{c_{1},c_{2},c_{3}}&=\oti_{c_{1}\ot c_{2},c_{3}}\circ (\oti_{c_{1},c_{2}}\x \text{id}_{F(c_{3})}),\\
\Gamma_{\Po}(\ot;1_{F},\ot)_{c_{1},c_{2},c_{3}}&=\oti_{c_{1}, c_{2}\ot c_{3}}\circ (\text{id}_{F(c_{1})}\x \oti_{c_{2},c_{3}}).
\end{align*}
By equation $(\textbf{c.8})$, $\Gamma_{\Po}(\ot;\ot,1_{F})_{c_{1},c_{2},c_{3}}$ and $\Gamma_{\Po}(\ot;1_{F},\ot)_{c_{1},c_{2},c_{3}}$ agree as functors. On the other hand, the $3$-linear map $\Gamma_{\Po}(\ot;\ot,1_{F})_{c_{1},c_{2},c_{3}}$ has structural maps $\delta^{1},\delta^{2},\delta^{3}$, where $\delta^{1}$ is given by
\begin{equation}
\label{delta1}
\delta^{1}=(d_{c_{1},c_{2}}^{l}\ot_{c_{1}\ot c_{2}, c_{3}} \text{id})\circ d_{c_{1}\ot c_{2},c_{3}}^{l}.
\end{equation}
Similarly we can find equations $\delta^{2}$ and $\delta^{3}$. 
On the other hand, the $3$-linear map $\Gamma_{\Po}(\ot;1_{F},\ot)_{c_{1},c_{2},c_{3}}$ has structural maps $\delta^{1}_{*},\delta^{2}_{*},\delta^{3}_{*}$, where $\delta^{1}_{*}$ is given by
\begin{equation}
\label{delta2}
\delta^{1}_{*}=d^{l}_{c_{1},c_{2}\ot c_{3}}.
\end{equation}
Similarly we can find equations $\delta^{2}_{*}$ and $\delta^{3}_{*}$. By diagram $(\textbf{c.9})$, we see that (\ref{delta1}) and (\ref{delta2}) agree. The same statement is true for $\delta_{2}$, $\delta_{3}$ and $\delta^{2}_{*}$, $\delta^{3}_{*}$ and it follows by conditions $(\textbf{c.9})$ and $(\textbf{c.10})$.

To construct the multifunctor $S$, we need to define for $k\ge 0$ a functor
\[
S_{k}:\Sigma_{k}\to \Po(F,...,F;F)
\]
that respects the action of $\Sigma_{k}$ and the multiproduct. As $\Sigma_{k}$ is a discrete category, it suffices to define $S_{k}(\sigma)$ for all $\sigma\in \Sigma_{k}$. We will define first $S_{k}(1_{k})$ for $k\ge 0$. For $k=0$ we define $S_{0}(1_{0})=1$. To define $S_{k}(1_{k})$ for $k\ge 1$, we will use the following notation. Suppose that $c_{1},...,c_{k}\in \O_{\E}$ and that $x_{i}$ is an object in the permutative category $F(c_{i})$, for $1\le i\le k$. Define
\[
x_{1,...,k}=x_{1}\ot \cdots \ot x_{k}=(\cdots((x_{1}\oti_{c_{1},c_{2}} x_{2})\oti_{c_1\ot c_{2},c_{3}}x_{3}) \cdots\oti_{c_{1}\ot \cdots \ot c_{k-1},c_{k}} x_{k}).
\]
Note that by condition $(\textbf{c.8})$, this is well defined and equals any rearranging of the parenthesis using the respective operations $\ot$. With this in mind, define
\[
S_{k}(1_{k})_{c_{1},...,c_{k}}(x_{1},...,x_{k})=x_{1}\ot \cdots \ot x_{k}=x_{1,...,k}.
\]
We need to define the structural maps of the $S_{k}(1_{k})_{c_{1},...,c_{k}}$'s in order for it to be a $k$-linear map. To do this, note that by condition $(\textbf{c.8})$ we have that 
\[
S_{k}(1_{k})_{c_{1},...,c_{k}}(x_{1},...,x_{k})=x_{1,...,k}=x_{1,...,i-1}\ot x_{i}\ot x_{i+1,...,k}.
\] 
Then we define
\begin{align*}
\delta^{i}:&S_{k}(1_{k})_{c_{1},...,c_{k}}(x_{1},...,x_{i},...,x_{k})\op S_{k}(1_{k})_{c_{1},...,c_{k}}(x_{1},...,x'_{i},...,x_{k})\\
= &x_{1,...,i-1}\ot x_{i}\ot x_{i+1,...,k}\op x_{1,...,i-1}\ot x'_{i}\ot x_{i+1,...,k} \\
&\to x_{1,...,i-1}\ot (x_{i}\op x'_{i})\ot x_{i+1,...,k}\\
=&S_{k}(1_{k})_{c_{1},...,c_{k}}(x_{1},...,x_{i}\op x'_{i},...,x_{k}),
\end{align*}
to be the diagonal map in the following commutative diagram
\[
\xymatrix@C=-8pc{
 x_{1,...,i-1}\ot x_{i}\ot x_{i+1,...,k}\op x_{1,...,i-1}\ot x'_{i}\ot x_{i+1,...,k}\ar[rd]^-{d^{l}}\ar[dd]_-{d^{r}}   & \\
& (x_{1,...,i-1}\ot x_{i}\op x_{1,...,i-1}\ot x'_{i})\ot x_{i+1,...,k}\ar[dd]^-{d^{r}\ot\text{id}}\\
x_{1,...,i-1}\ot( x_{i}\ot x_{i+1,...,k}\op x'_{i}\ot x_{i+1,...,k})\ar[rd]_-{\text{id}\ot d^{l}}         & \\
& x_{1,...,i-1}\ot (x_{i}\op x'_{i})\ot x_{i+1,...,k}.
}
\]
The previous diagram is commutative by condition $(\textbf{c.10})$ and by a direct computation we see that the condition of commutation with $\op$ for $\delta^{i}$ is satisfied.
\nin
In general, we define for $\sigma\in \Sigma_{k}$
\[
S_{k}(\sigma)=\sigma^{*}S_{k}(1_{k}).
\]
This way defined, $S_{k}$ is clearly respects the $\Sigma_{k}$ action. Let's see that it respect the multiproduct. To do so, take $\phi_{i}\in \Sigma_{j_{i}}$ for $1\le i\le k$ and $\sigma\in \Sigma_{k}$. We need to verify the equality
\[
S_{j_{1}+\cdots +j_{k}}(\Gamma_{\Sigma_{*}}(\sigma;\phi_{1},...,\phi_{k}))=\Gamma_{\Po}(S_{k}(\sigma);S_{j_{1}}(\phi_{1}),...,S_{j_{k}}(\phi_{k})).
\]
In $\Sigma_{*}$, the multiproduct is defined by 
\[
\Gamma_{\Sigma_{*}}(\sigma;\phi_{1},...,\phi_{k})=\sigma_{\left\langle j_{1},...,j_{k}\right\rangle}\circ (\phi_{1}\op \cdots \op \phi_{k}).
\]
Here, $\phi_{1}\op \cdots \op \phi_{k}$ means the image of $\phi_{1},...\phi_{k}$ under the canonical embedding 
\[
\Sigma_{j_{1}}\x \cdots \x \Sigma_{j_{k}}\to \Sigma_{j_{1}+\cdots j_{k}}.
\] 
By definition we have
\begin{align*}
S_{j_{1}+\cdots +j_{k}}(\Gamma_{\Sigma_{*}}(\sigma;\phi_{1},...,\phi_{k}))=&(\Gamma_{\Sigma_{*}}(\sigma;\phi_{1},...,\phi_{k}))^{*}S_{j_{1}+\cdots +j_{k}}(1_{j_{1}+\cdots +j_{k}})\\
=& (\sigma_{\left\langle j_{1},...,j_{k}\right\rangle}\circ (\phi_{1}\op \cdots \op \phi_{k}))^{*}S_{j_{1}+\cdots +j_{k}}(1_{j_{1}+\cdots +j_{k}}).
\end{align*}
On the other hand, 
\[
\Gamma_{\Po}(S_{k}(\sigma);S_{j_{1}}(\phi_{1}),...,S_{j_{k}}(\phi_{k}))=\Gamma_{\Po}(\sigma^{*}S_{k}(1_{k});\phi_{1}^{*}S_{j_{1}}(1_{j_{1}}),...,\phi_{k}^{*}S_{j_{k}}(1_{j_{k}})).
\]
Using properties Multi(3) and Multi(4) in the multicategory $\Po$, we obtain
\begin{align*}
&\Gamma_{\Po}(S_{k}(\sigma);S_{j_{1}}(\phi_{1}),...,S_{j_{k}}(\phi_{k}))\\
&=(\sigma_{\left\langle j_{1},...,j_{k}\right\rangle}\circ (\phi_{1}\op \cdots \op \phi_{k}))^{*}\Gamma_{\Po}(S_{k}(1_{k});S_{j_{1}}(1_{j_{1}}),...,S_{j_{k}}(1_{j_{k}})).
\end{align*}
Thus it suffices to show that
\[
S_{j_{1}+\cdots +j_{k}}(1_{j_{1}+\cdots +j_{k}})=\Gamma_{\Po}(S_{k}(1_{k});S_{j_{1}}(1_{j_{1}}),...,S_{j_{k}}(1_{j_{k}})).
\]
But this follows by induction and equation (\ref{3morph}). Since $S$ respects identities by definition, we see that $S$ is a multifunctor mapping the only object of $\Sigma_{*}$ to $F$. 
Finally, the fact that these correspondences are inverses of each other follows by a direct computation.
\qed

We now discuss the case where the operations $\ot_{c_{1},c_{2}}$, satisfy some sort of commutativity up to coherent isomorphisms. Let's denote by $E\Sigma_{*}$ the category valued operad that for each $j$, corresponds $E\Sigma_{j}$ the translation category of $\Sigma_{j}$; that is, the object set of $E\Sigma_{j}$ is $\Sigma_{j}$ and there is exactly one morphism between any pair of objects of $E\Sigma_{j}$. We can see $E\Sigma_{*}$ as a multicategory with only one object and that is enriched over $\cat$. In the same spirit as in the previous theorem, we want to give an explicit description to the enriched multifunctors 
\[
S:E\Sigma_{*}\to \Po.
\]
Such a multifunctor maps the unique object of $E\Sigma_{*}$ to an object $F$ of $\Po$. The multifunctor $E\Sigma_{*}\to \Po$, then factors through the multicategory generated by $F$ in $\Po$. We can see each $\Sigma_{j}$ as the trivial category with only the identity morphisms. Thus we have an inclusion functor $i_{j}:\Sigma_{j}\to E\Sigma_{j}$ for every $j\ge 0$. The collection of all $i_{j}$ gives rise to a multifunctor 
\[
i:\Sigma_{*}\to E\Sigma_{*}.
\]
This multifunctor is trivially enriched over $\cat$. By composing $S$ with $i$, we see that for $F$ conditions $(\textbf{c.1})-(\textbf{c.10})$ are satisfied. In addition, we also have the following properties. Let $\xi\in \Sigma_{2}$ is the map that permutes the two elements of $\{1,2\}.$ Then we can find a map of $2$-morphisms  
\[
\mu:\ot\to \xi^{*}\ot;
\]
that is, for each pair of objects $c_{1}$ and $c_{2}$ of $\E$, we have a map of $2$-linear maps 
\[
\mu_{c_{1},c_{2}}:x\oti_{c_{1},c_{2}} y\to F(\tau_{\xi,c_{2},c_{1}})(y\oti_{c_{2},c_{1}} x)
\]
satisfying the following coherence diagrams
\[\tag{\textbf{c.11}}
\xymatrix{
x\oti_{c_{1},c_{2}}0_{c_{2}}\ar[r]^-{\mu_{c_{1},c_{2}}}\ar[d]_-{=}   &F(\tau_{\xi,c_{2},c_{1}})(0_{c_{2}}\oti_{c_{2},c_{1}}x)\ar[d]^-{=}\\
0_{c_{1}\ot c_{2}}\ar[r]_-{=}					&F(\tau_{\xi,c_{2},c_{1}})(0_{c_{2}\ot c_{1}}),
}
\]
\[\tag{\textbf{c.12}}
\xymatrix{
x\oti_{c_{1},c_{2}}y\ar[r]^-{\mu_{c_{1},c_{2}}}\ar[rd]_-{=}   &F(\tau_{\xi,c_{2},c_{1}})(y\oti_{c_{2},c_{1}}x)\ar[d]^-{F(\tau_{\xi,c_{2},c_{1}})(\mu_{c_{2},c_{1}})}\\
					&F(\tau_{\xi,c_{2},c_{1}})F(\tau_{\xi,c_{1},c_{2}})(x\oti_{c_{1},c_{2}}y),
}
\]
\[\tag{\textbf{c.13}}
\xymatrix@C=-6pc{&(x\op x')\oti_{c_{1},c_{2}}y
\ar[ddr]^-{\mu_{c_{1},c_{2}}}
\\
(x\oti_{c_{1},c_{2}}y)\op (x'\oti_{c_{1},c_{2}}y)
\ar[ur]^-{d^{l}_{c_{1},c_{2}}}
\ar[dd]_-{\mu_{c_{1},c_{2}}\op \mu_{c_{1},c_{2}}}
\\
&&F(\tau_{\xi,c_{2},c_{1}})(y\oti_{c_{2},c_{1}}(x\op x')),
\\
F(\tau_{\xi,c_{2},c_{1}})(y\oti_{c_{2},c_{1}}x)\op F(\tau_{\xi,c_{2},c_{1}})(y\oti_{c_{2},c_{1}}x')
\ar[dr]_-{\lambda_{F(\tau_{\xi,c_{2},c_{1}})}}
\\
&F(\tau_{\xi,c_{2},c_{1}})((y\oti_{c_{2},c_{1}}x)\op(y\oti_{c_{2},c_{1}}x'))
\ar[uur]_-{F(\tau_{\xi,c_{2},c_{1}})(d^{r}_{c_{2},c_{1}})}}
\]
\[\tag{\textbf{c.14}}
\xymatrix{
x\oti_{c_{1},c_{2}\ot c_{3}}(y\oti_{c_{2},c_{3}}z)\ar[r]^-{\text{id}\oti_{c_{1},c_{2}\ot c_{3}}\mu_{c_{2},c_{3}}}\ar[d]_-{=}   &x\oti_{c_{1},c_{2}\ot c_{3}}(F(\tau_{\xi,c_{3},c_{2}})z\oti_{c_{3},c_{2}}y)\ar[d]^-{=}\\
(x\oti_{c_{1},c_{2}}y)\oti_{c_{1}\ot c_{2},c_{3}}z\ar[d]_-{\mu_{c_{1}\ot c_{2},c_{3}}}         &F(1\ot \tau_{\xi,c_{3},c_{2}})(x\oti_{c_{1},c_{3}\ot c_{2}}(z\oti_{c_{3},c_{2}} y))\ar[d]^-{=}\\
F(\tau_{\xi,c_{3},c_{1}\ot c_{2}})(z\oti_{c_{3},c_{1}\ot c_{2}}(x\oti_{c_{1},c_{2}}y))\ar[d]_-{=}					&F(1\ot \tau_{\xi,c_{3},c_{2}})((x\oti_{c_{1},c_{3}}z)\oti_{c_{1}\ot c_{3},c_{2}} y)\ar[d]^-{F(\text{id}\ot\tau_{\xi,c_{3},c_{2}})(\mu_{c_{1},c_{3}})\ot \text{id}}\\
F(\tau_{\xi,c_{3},c_{1}\ot c_{2}})((z\oti_{c_{3},c_{1}}x)\oti_{c_{3}\ot c_{1},c_{2}}y)\ar[r]_-{=}					&F(1\ot \tau_{\xi,c_{3},c_{2}})(F(\tau_{\xi,c_{3},c_{1}})(z\oti_{c_{3},c_{1}} x)\oti_{c_{1}\ot c_{3},c_{2}}y).
}
\]
The equalities in the previous diagram are obtained by applying the naturality property of $\ot$ and the associativity property of $\ot$ which are satisfied as conditions $(\textbf{c.1}-\textbf{c.10})$ hold. We have the following theorem.

\begin{theorem}\label{multiEsigma}
Having a multifunctor $E\Sigma_{*}\to \Po$ mapping the only object of $E\Sigma_{*}$ to a functor $F:\E\to \P$ is equivalent to having a functor $F:\E\to \P$, for each pair $c_{1},c_{2}$ a functor 
\[
\oti_{c_{1},c_{2}}:F(c_{1})\x F(c_{2})\to F(c_{1}\ot c_{2}), 
\]
a unit $1$ which is an object in $F(1_{\E})$, distributivity maps 
\begin{align*}
d_{c_{1},c_{2}}^{l}&:(x\oti_{c_{1},c_{2}}y)\op (x'\oti_{c_{1},c_{2}}y)\to (x\op x')\oti_{c_{1},c_{2}}y \\
d_{c_{1},c_{2}}^{r}&:(x\oti_{c_{1},c_{2}}y)\op (x\oti_{c_{1},c_{2}}y')\to x\oti_{c_{1},c_{2}}(y\op y') 
\end{align*} 
and a natural transformation
\[
\mu_{c_{1},c_{2}}:x\oti_{c_{1},c_{2}} y\to F(\tau_{\xi,c_{2},c_{1}})(y\oti_{c_{2},c_{1}} x)
\]
satisfying conditions $(\textbf{c.1})-(\textbf{c.14})$.
\end{theorem}

\Proof
Suppose that $S:E\Sigma_{*}\to \Po$ is an enriched multifunctor sending the only object of $E\Sigma_{*}$ to a functor $F:\E\to \P$. Let  $\ot=S_{2}(1_{2})$, $\mu$ the image of the unique morphism in $E\Sigma_{2}$ between $1_{2}$ and $\xi$ under $S_{2}$ and $1$ the image of $1_{0}$ under $S_{0}$. Then we see by the way we derived conditions $(\textbf{c.1})-(\textbf{c.14})$ that these are satisfied by $F, \ot, \xi, 1$.  

On the other hand, suppose we have a functor $F:\E\to \P$, a collection of $2$-linear maps $\ot_{c_{1},c_{2}}$, natural transformations $\mu_{c_{1},c_{2}}$ for each pair $c_{1}$ and $c_{2}\in \O_{\E}$ and a unit $1$ that satisfy properties $(\textbf{c.1})-(\textbf{c.14})$. We want to define a multifunctor $S:E\Sigma_{*}\to \Po$. To begin, note that in particular conditions $(\textbf{c.1})-(\textbf{c.10})$ are satisfied and thus by Theorem \ref{multisigma} we have a multifunctor $S:\Sigma_{*}\to \Po$ sending the only object of $\Sigma_{*}$ to $F$. We want to see that this multifunctor extends to a multifunctor $S:E\Sigma_{*}\to \Po$. We claim that the natural transformations $\{\mu_{c_{1},c_{2}}\}$, form a map of $2$-morphisms from $\ot$ to $\xi^{*}\ot$. Indeed, if $c_{1}$ and $c_{2}\in \O_{\E}$, and $x_{i}$ is an object in the category $F(c_{i})$ for $i=1,2$, then $\mu_{c_{1},c_{2}}$ is a natural map
\[
\mu_{c_{1},c_{2}}(x_{1},x_{2}):x_{1}\oti_{c_{1},c_{2}} x_{2}\to (\xi^{*}\oti_{c_{1},c_{2}}) (x_{1}, x_{2}).
\]
Also, $\mu_{c_{1},c_{2}}(x_{1},x_{2})$ is the identity whenever $x_{1}$ or $x_{2}$ is zero. This is a consequence of condition $(\textbf{c.11})$. In order to conclude that $\mu_{c_{1},c_{2}}$ is a map of $2$-linear maps, we need to show that the following diagram is commutative
\[
\xymatrix{
x_{1}\oti_{c_{1},c_{2}} x_{2}\op x'_{1}\oti_{c_{1},c_{2}} x_{2}\ar[r]^-{d^{l}_{c_{1},c_{2}}}\ar[d]_-{\mu_{c_{1},c_{2}\op\mu_{c_{1},c_{2}}}}   &(x_{1}\op x'_{1})\oti_{c_{1},c_{2}} x_{2}\ar[d]^-{\mu_{c_{1},c_{2}}}\\
(\xi^{*}\oti_{c_{1},c_{2}})(x_{1}, x_{2})\op (\xi^{*}\oti_{c_{1},c_{2}}) (x'_{1}, x_{2})\ar[r]_-{\delta^{1}_{\xi^{*}\ot_{c_{1},c_{2}}}}         &(\xi^{*}\oti_{c_{1},c_{2}})(x_{1}\op x'_{1},x_{2})
}
\]
and a similar one for $\delta^{2}_{\xi^{*}\ot_{c_{1},c_{2}}}$. But by definition, 
\[
\delta^{1}_{\xi^{*}\ot_{c_{1},c_{2}}}=F(\tau_{\xi,c_{2},c_{1}})(d^{r}_{c_{2},c_{1}})\circ \lambda_{F(\tau_{\xi,c_{2},c_{1}})},
\]
and thus the commutativity of the previous diagram follows from condition $(\textbf{c.13})$. We need to define a functor $S_{k}:E\Sigma_{k}\to \Po(F,...,F;F)$. Suppose $\rho, \sigma\in \Sigma_{k}$ and $f:\rho\to \sigma$ is the only morphism in $E\Sigma_{k}$ from $\rho$ to $\sigma$. Using the action of $\Sigma_{k}$ we only need to define $S_{k}(f)$ in the case of $\rho=1_{k}$. Thus we need to define a morphism on the category $\Po(F,...,F;F)$ from $S_{k}(1_{k})$ to $S_{k}(\sigma)$. Let's consider first the case where $\sigma$ is a transposition. Thus assume $\sigma=(m,n)$, where $1\le m< n\le k$. Take $c_{1},...,c_{k}\in \O_{\E}$ and $x_{i}$ an object of $F(c_{i})$ for $i=1,...,k$. By definition we have
\[
S_{k}(\sigma)_{c_{1},...,c_{k}}(x_{1},...,x_{k})=F(\tau_{(m,n),c_{i},c_{j}})(x_{1,...,m-1}\ot x_{m}\ot x_{m+1,...,n-1}\ot x_{n}\ot x_{n+1,...,k})
\]
Define $S_{k}(f)_{c_{1},...,c_{k}}$ to be the following composite in the category $\Po(F,...,F;F)$
\begin{align*}
&x_{1,...,k}= x_{1,...,m-1}\ot x_{m}\ot x_{m+1,...,n-1}\ot x_{n}\ot x_{n+1,...,k}\\
&\to  F(\text{id}\ot \tau_{\xi}\ot \text{id})(x_{1,...,m-1}\ot x_{m+1,...,n-1}\ot x_{m} \ot x_{n}\ot x_{n+1,...,k})\\
&\to  F(\text{id}\ot \tau_{\xi}\ot \text{id})F(\text{id}\ot \tau_{\xi}\ot \text{id})\\
&(x_{1,...,m-1}\ot x_{m+1,...,n-1}\ot x_{n} \ot x_{m}\ot x_{n+1,...,k})\\
&\to  F(\text{id}\ot \tau_{\xi}\ot \text{id})F(\text{id}\ot \tau_{\xi}\ot \text{id})F(\text{id}\ot \tau_{\xi}\ot \text{id})\\
&(x_{1,...,m-1}\ot x_{n}\ot x_{m+1,...,n-1} \ot x_{m}\ot x_{n+1,...,k})\\
&= F(\tau_{\sigma})(x_{\sigma(1),...,\sigma(k)})
\end{align*}
By condition $(\textbf{c.14})$, this is well defined and agrees with any other composition of maps obtained by applying $\mu$ to a sequence that starts in $x_{1,...,k}$ and finish in $F(\tau_{\sigma,c_{1},...,c_{k}})(x_{\sigma(1),...,\sigma(k)})$. As any element of $\Sigma_{k}$ can be obtained by composition of transpositions, we see that we can  extend this definition by using the $\Sigma_{k}$ action. This way we define $S_{k}(f)$, where $f$ is a morphism in $E\Sigma_{k}$ from $1_{k}$ to $\sigma$. Using the $\Sigma_{k}$ action we can extend the definition to all the morphisms in $E\Sigma_{k}$. We need to verify that this way defined, we obtain functors 
\[
S_{k}:E\Sigma_{k}\to \Po(F,...,F;F).
\]
This is done by induction on $k$. For $k=0$ and $k=1$ it is trivial as $E\Sigma_{i}$ is a discrete category. For $k=2$ there are only two morphisms that are not the identity, namely, $f:1_{2}\to \xi$ and its inverse. By condition $(\textbf{c.12})$ we see that $S_{2}(f)\circ S_{2}(f^{-1})=S_{2}(f^{-1})\circ S_{2}(f)=\text{id}$. The case where $k\ge 2$ follows easily by induction using condition $(\textbf{c.14})$ and is left for the reader. We are left to prove that $S_{k}$ respects the multiproduct on the level of morphisms of the categories $E\Sigma_{k}$ and $\Po(F,...,F;F)$. To see this, note that by the $\Sigma_{k}$ action and using properties Multi(3)-Multi(4) as in Theorem \ref{multisigma} we only need to prove that
\[
\Gamma_{\Po}(S_{k}(g);S_{j_{1}}(f_{1}),...,S_{j_{k}}(f_{r}))=S_{j_{1}+\cdot+j_{k}}(\Gamma_{E\Sigma_{*}}(g;f_{1},...,f_{k}))
\]
in the case where $f_{i}:1_{j_{i}}\to \sigma_{i}$, and $g:1_{k}\to \rho$. But for each $c_{ir}\in \O_{\E}$, for $1\le i\le k$ and $1\le r\le j_{i}$ we have that both
\[
\Gamma_{\Po}(S_{k}(g);S_{j_{1}}(f_{1}),...,S_{j_{k}}(f_{r}))_{c_{11},...,c_{kj_{k}}}\ \  \text{ and } \ \ S_{j_{1}+\cdot+j_{k}}(\Gamma_{E\Sigma_{*}}(g;f_{1},...,f_{k}))_{c_{11},...,c_{kj_{k}}} 
\]
are obtained by successive applications of $\mu$ starting and ending in the same place. By $(\textbf{c.14})$ it follows that these morphism have to agree.
\nin
Finally, the fact that these correspondences are inverses of each other follows by a direct computation.
\qed

The application of Theorem \ref{multiEsigma} that we are looking for is the following theorem.

\begin{theorem}\label{verifyconditions}
If $\Lambda:\D\to \C$ is a fibered bipermutative category then for the associated functor $\Psi:\A\to \P$, we can find functors 
\[
\oti_{\u,\v}:\Psi(\u)\x\Psi(\v)\to \Psi(\u\odot\v)
\]
for each pair $\u$ and $\v$ of objects of $\A$, a unit $\textbf{1}$ which is an object in $\Psi(())$, distributivity maps 
\begin{align*}
d_{\u,\v}^{l}&:(F\oti_{\u,\v}G)\op (F'\oti_{\u,\v}G)\to (F\op F')\oti_{\u,\v}G, \\
d_{\u,\v}^{r}&:(F\oti_{\u,\v}G)\op (F\oti_{\u,\v}G')\to F\oti_{\u,\v}(G\op G') 
\end{align*} 
and a natural transformation
\[
\mu_{\u,\v}:F\oti_{\u,\v} G\to \Psi(\tau_{\xi,\u,\v})(G\oti_{\v,\u} F)
\]
satisfying conditions $(\textbf{c.1})-(\textbf{c.14})$ of Theorem \ref{multiEsigma}.
\end{theorem}
\Proof
Let's begin by defining $\ot_{\u,\v}$. Take $\u=(u_{1},...,u_{n})$ and $\v=(v_{1},...,v_{m})$ and 
\[
F=\sum_{i=1}^{r}F_{i1}\ot\cdots\ot F_{in},\ \ \text{and} \ \
G=\sum_{j=1}^{s}G_{j1}\ot\cdots\ot G_{jm},
\]
objects of $\Psi(\u)$ and $\Psi(\v)$ respectively. Define
\[
F\oti_{\u,\v} G=\sum_{j=1}^{s}\sum_{i=1}^{r}F_{i1}\ot\cdots F_{in}\ot G_{j1}\ot\cdots\ot G_{jm}.
\]
In the case where $\v=()$, we define
\[
F\oti_{\u,()} G=\sum_{j=1}^{s}\sum_{i=1}^{r}F_{i1}\ot\cdots \ot(F_{in}\ot G).
\]
Note that every object of $\C/\u\odot \v$ can be written in the form $\f\odot\g=(f_{1},...,f_{n},g_{1},...,g_{m})$, where $\f=(f_{1},...,f_{n})$ is an object of $\C/\u$ and $\g=(g_{1},...,g_{m})$ and object of $\C/\v$. Then using the distributivity maps $d^{l}$ and $d^{r}$ of $\D$ we can obtain a coherent isomorphism 
\[
\lambda_{\f,\g}:F(\f)\ot G(\g)\stackrel{\cong}{\rightarrow} (F\oti_{\u,\v}G)(\f\odot \g). 
\]
For morphisms $\alpha:F\to F'$ and $\beta:G\to G'$, define $\alpha\ot \beta$ to be the base point preserving natural transformation defined by the following diagram
\[
\xymatrix{
(F\oti_{\u,\v}G)(\f\odot \g) \ \ \ar[r]^-{(\alpha\ot\beta)_{\f\odot\g}}\ar[d]_-{\lambda_{\f,\g}^{-1}}&\ \ (F'\ot G')(\f\odot \g)\\
F(\f)\ot G(\g)\ \ \ar[r]_-{\alpha_{\f}\ot \beta_{\g}}&\ \ F'(\f)\ot G'(\g).\ar[u]_-{\lambda_{\f,\g}}
}
\]
Using the coherence of the distributivity map and the fact that these preserve the base point, we see that this way defined $\alpha\ot \beta$ is a base point preserving natural transformation. This defines the functor $\ot_{\u,\v}$.
\nin
We also need to define the distributivity maps 
\begin{align*}
d_{\u,\v}^{l}&:(F\oti_{\u,\v}G)\op (F'\oti_{\u,\v}G)\to (F\op F')\oti_{\u,\v}G, \\
d_{\u,\v}^{r}&:(F\oti_{\u,\v}G)\op (F\oti_{\u,\v}G')\to F\oti_{\u,\v}(G\op G'). 
\end{align*} 
Given an object $\f\odot \g$ of $\C/\u\odot\v$ we define $(d_{\u,\v}^{l})_{\f\odot\g}$ by the following diagram
\[
\xymatrix{
(F\oti_{\u,\v}G)(\f\odot\g)\op (F'\oti_{\u,\v}G)(\f\odot \g)\ \ \ar[r]^-{(d_{\u,\v}^{l})_{\f\odot\g}}\ar[d]_-{\lambda^{-1}_{\f,\g}\op\lambda^{-1}_{\f,\g}}&\ \ ((F\op F')\oti_{\u,\v}G)(\f\odot\g)\\
F(\f)\ot G(\g)\op F'(\f)\ot G(\g)\ \ \ar[r]_-{d^{l}}&\ \ (F(\f)\op F'(\f))\ot G(\g).\ar[u]_-{\lambda_{\f,\g}\op\lambda_{\f,\g}} 
}
\]
We define $d_{\u,\v}^{r}$ in a similar fashion.
\nin
This way defined $\ot_{\u,\v}$ satisfies conditions $(\bf{c.1})-(\bf{c.3})$ as similar conditions are satisfied in $\D$. We also have that for any objects 
\[
F=\sum_{i=1}^{r}F_{i1}\ot\cdots\ot F_{in}\ \ \text{and} \ \
G=\sum_{j=1}^{s}G_{j1}\ot\cdots\ot G_{jm}
\]
of $\Psi(\u)$ and $\Psi(\v)$ respectively, we have
\begin{align*} 
F\oti_{\u,\v}\Os_{\v}&=\sum_{i=1}^{r}F_{i1}\ot\cdots\ot F_{in}\ot \Os_{v_{1}}\ot\cdots \Os_{v_{m}}=\Os_{\u\odot\v}\\
\Os_{\u}\oti_{\u,\v} G&=\sum_{j=1}^{s}\Os_{u_{1}}\ot\cdots\ot \Os_{\u_{n}}\ot G_{j1}\ot\cdots\ot G_{jm}=\Os_{\u\odot\v}
\end{align*}
and thus condition $(\bf{c.4})$ is satisfied.
\nin
We want to show next that if $(q,\f):\u\to \u'$ and $(p,\g):\v\to \v'$ are morphisms in $\A$, then the following diagram is commutative
\[
\xymatrix{
\Psi(\u)\x\Psi(\v)\ \ \ \ar[r]^-{\Psi(\sigma,\f)\x\Psi(\rho,\g)}\ar[d]_-{\oti_{\u,\v}}&\ \ \ \Psi(\u')\x \Psi(\v')\ar[d]^-{\oti_{\u',\v'}}\\
\Psi(\u\odot\v)\ \ \ \ar[r]_-{\Psi(\sigma\odot\rho,\f\odot\g)}&\ \ \ \Psi(\u'\odot\v'),
}
\]
where $p\odot q:\underline{n+m}\to \underline{n'+m'}$ is the injective map that acts as $q$ in the first $n$-elements and acts as $p$ in the last $m$-elements. To show the commutativity of this diagrams suppose first that $q=\sigma\in \Sigma_{n}$ and $p=\rho\in \Sigma_{m}$ are permutations. Take 
\[
F=\sum_{i=1}^{r}F_{i1}\ot\cdots\ot F_{in} \ \ \text{and} \ \ G=\sum_{j=1}^{s}G_{j1}\ot\cdots\ot G_{jm}
\]
objects of  $\Psi(\u)$ and $\Psi(\v)$ respectively. Then by definition we have 
\begin{align*}
&\Psi(\sigma,\f)F\oti_{\u',\v'}\Psi(\rho,\g)G\\ 
=&\left(\sum_{i=1}^{r}f^{*}_{1}F_{i\sigma^{-1}(1)}\ot\cdots\ot f^{*}_{n}F_{i\sigma^{-1}(n)}\right)\oti_{\u',\v'}\left(\sum_{i=1}^{r}g^{*}_{1}G_{\rho^{-1}(1)}\ot\cdots\ot g^{*}_{n}G_{i\rho^{-1}(n)}\right)\\
=&\sum_{j=1}^{s}\sum_{i=1}^{r}\left(f^{*}_{1}F_{i\sigma^{-1}(1)}\ot\cdots\ot f^{*}_{n}F_{i\sigma^{-1}(n)}\ot g^{*}_{1}G_{\rho^{-1}(1)}\ot\cdots\ot g^{*}_{n}G_{i\rho^{-1}(n)}\right).
\end{align*}
On the other hand, we have
\begin{align*}
&\Psi(\sigma\odot \rho,\f\odot\g)\left(F\oti_{\u',\v'}G\right)\\
=&\Psi(\sigma\odot \rho,\f\odot\g)\left(\sum_{j=1}^{s}\sum_{i=1}^{r}F_{i1}\ot\cdots F_{in}\ot G_{j1}\ot\cdots\ot G_{jm}\right)\\
=&\sum_{j=1}^{s}\sum_{i=1}^{r}\left(f^{*}_{1}F_{i\sigma^{-1}(1)}\ot\cdots\ot f^{*}_{n}F_{i\sigma(n)}\ot g^{*}_{1}G_{\rho^{-1}(1)}\ot\cdots\ot g^{*}_{n}G_{i\rho^{-1}(n)}\right).
\end{align*}
This shows that condition $(\bf{c.5})$ is satisfied. In addition, as each $\psi(\sigma,\f)$ is strict map, and the distributivity maps are natural, this also shows that condition $(\bf{c.6})$ is satisfied. A similar argument shows that the same is true for a general morphism in $(q,\f)$ and $(p,\g)$ in $\A$.
\nin
We now want to construct a multiplicative unit $\bf{1}$. This is an object in $\Psi(())$. By definition, the objects of $\Psi(())$ are the functors of the form $F=\sum_{i=1}^{r}F_{i}$ with $F_{i}$ an object of $\text{Hom}_{\C}(\C/(),\D)$. In particular, the functor $\bf{1}:\C/()\to \D$ that sends the only object of $\C/()$ to the unit $1$ of $\D$ is an object of $\Psi(())$. We claim that this is a unit for $\ot_{\u,\v}$. Indeed, if $F=\sum_{i=1}^{r}F_{i1}\ot\cdots\ot F_{in}$ is an object in $\Psi(\u)$, then 
\[
F\oti_{\u,()}\textbf{1}=\sum_{i=1}^{r}F_{i1}\ot\cdots\ot (F_{in}\ot 1)=F
\]
and similarly $\textbf{1}\ot_{(),\v}G=G$ for any object $G$ of $\Psi(\v)$. This shows that condition $(\bf{c.7})$ is satisfied. On the other hand, from the coherences in the category $\D$ we see at once that conditions $(\bf{c.8})-(\bf{c.10})$ are also satisfied. 
\nin
We now want to study the symmetry of the functor $\ot_{\u,\v}$. Note that given objects $\u=(u_{1},...,u_{n})$ and $\v=(v_{1},...,v_{m})$ of $\A$, then the symmetry isomorphism $\tau_{\gamma,\v,\u}=\gamma^{\odot}:\v\odot\u\to \u\odot\v$ is given by the pair $(\xi,\underline{\text{id}})$, where $\xi\in \Sigma_{n+m}$ is the permutation that interchanges the first $m$-block with the last $n$-block and $\underline{\text{id}}$ is the $(n+m)$-tuple of copies of $\text{id}$. With this in mind note that if
\[
F=\sum_{i=1}^{r}F_{i1}\ot\cdots\ot F_{in} \ \ \text{and}\ \
G=\sum_{j=1}^{s}G_{j1}\ot\cdots\ot G_{jm},
\]
are objects of $\Psi(\u)$ and $\Psi(\v)$ respectively, then we obtain
\begin{align*}
\Psi(\tau_{\gamma,\v,\u})(G\oti_{\v,\u}F)&=\Psi(\tau_{\gamma,\v,\u})\left(\sum_{j=1}^{s}\sum_{i=1}^{r}G_{j1}\ot\cdots\ot G_{jm}\ot F_{i1}\ot\cdots\ot F_{in}\right)\\
&=\sum_{j=1}^{s}\sum_{i=1}^{r}F_{i1}\ot\cdots\ot F_{in}\ot G_{j1}\ot\cdots\ot G_{jm}
\end{align*}
and
\begin{align*}
F\ot_{\u,\v}G=\sum_{i=1}^{r}\sum_{j=1}^{s}F_{i1}\ot\cdots\ot F_{in}\ot G_{j1}\ot\cdots\ot G_{jm}.
\end{align*}
Define 
\[
\mu_{\u,\v}:F\oti_{\u,\v}G\to \Psi(\tau_{\gamma,\v,\u})(G\oti_{\v,\u}F)
\]
to be the natural isomorphism that rearranges the summation order; that is, $\mu_{\u,\v}$ is an iterated application of the natural isomorphism $\gamma^{\op}$ in $\D$. With this definition, it follows by the use of coherence theory that the conditions $(\bf{c.11})-(\bf{c.14})$ of theorem \ref{multiEsigma} are satisfied.
\qed

\begin{corollary}
A fibered bipermutative category $\Lambda:\D\to \C$ gives rise to a multifunctor $T_{1}:E\Sigma_{*}\to \P^{\A}$. By composing with the multifunctor $K_{*}$ we obtain a multifunctor $T:E\Sigma_{*}\to \Sym^{\A}$.
\end{corollary}

\section{Model Categories}

In this section we use the machinery of closed model categories to rigidify a given multifunctor $T:E\Sigma_{*}\to \Sym^{\E}$, as to obtain a multifunctor $T':*\to \Sym^{\E}$, here $\E$ is a general permutative category. The application that we have in mind is the multifunctor $T:E\Sigma_{*}\to \Sym^{\A}$ provided by the previous corollary. Here $*$ is the operad whose algebras are commutative monoids which can be seen as a multicategory with only one object. Later we show that having a multifunctor $*\to \Sym^{\E}$ is equivalent to having a lax map $\E\to \Sym$. Thus in the case of the $T:E\Sigma_{*}\to \Sym^{\A}$ provided by the previous corollary, we obtain a lax map $\A\to \Sym$. Finally, we show that this map induces a lax map $(\C^{op})^{-1}\C^{op}\to \Sym$ via a right adjoint in a Quillen equivalence coming from a lax map $\A\to (\C^{op})^{-1}\C^{op}$. 
\nin
The core of this section is to give a model structure to $\Mi_{1}\Fi_{2}(\M,\E;\Sym)$. This is done in \cite{Elmendorf} for the case that $\E$ is the trivial category. In this case, the category $\Mi_{1}\Fi_{2}(\M,\E;\Sym)$ is isomorphic to the category of multifunctors $\M\to \Sym$. The case where $\E$ is a general permutative category follows by a straight forward modification of the construction in \cite{Elmendorf}. We include the proofs here for completeness following very closely the treatment in \cite[Section 11]{Elmendorf}. Therefore no originality is claimed for this part of the work. The reader is suggested to read the beautiful work of Elmendorf and Mandell \cite{Elmendorf}.

To begin, we recall briefly cofibrantly generated model categories and the positive model structure on symmetric spectra. For a complete discussion about cofibrantly generated model categories we refer the reader to \cite{Shipley} and \cite{Hov} and for symmetric spectra  \cite{HSS} and \cite{MMSS}.
\nin
Let $\lambda$ an ordinal and $\C$ a cocomplete category. We can see $\lambda$ as a category with a unique arrow $f:x\to y$ whenever $x\le y$. By a $\lambda$-sequence we mean a colimit preserving functor $X:\lambda\to \C$. We refer to the map $X_{0}\to \text{colim}_{\beta<\lambda}X_{\beta}$  as the composition of the $\lambda$-sequence. 

\begin{definition}
Let $I$ be a set of maps in a category $\C$ that is cocomplete. A map $f:A\to B$ in $\C$ is said to be a relative $I$-cell complex if there exists an ordinal $\lambda$ and a $\lambda$-sequence $X:\lambda\to \C$ such that $f$ is the composition of $X$ and for each $\beta$ with $\beta+1<\lambda$ there is a pushout
\[
\xymatrix{
C_{\beta}\ar[r]\ar[d]_{g_{\beta}}   &X_{\beta}\ar[d]\\
D_{\beta}\ar[r]         &X_{\beta+1}
}
\]
such that $g_{\beta}\in I$. The collection of relative $I$-cell complexes will be denoted by $I$-cell.
\end{definition}

Thus a relative $I$-complex is a (possibly) transfinite composition of maps that are pushouts of maps in $I$.

\begin{definition}
Let $I$ be a class of maps in a category $\C$. A map $f$ is said to be $I$-injective if it has the right lifting property with respect to any map in $I$. We denote by $I$-inj. the class of $I$-injective maps in $\C$. Similarly, a  map $f$ is said to be $I$-projective if it has the left lifting property with respect to any map in $I$. We denote by $I$-proj. the class of $I$-projective maps in $\C$
\end{definition}

\begin{definition}
A model structure on a small category $\C$ is said to be cofibrantly generated if there are sets $I$ and $J$ of maps such that the following conditions are satisfied.
\begin{itemize}
	\item The domains of the maps of $I$ are small relative to $I$-cell.
	\item The domains of the maps of $J$ are small relative to $J$-cell.
	\item The class of fibrations is $J$-inj.
	\item The class of trivial fibrations is $I$-inj.
\end{itemize}
\end{definition}

Let us review now briefly the positive closed model structure on symmetric spectra. This is an example of  cofibrantly generated closed model structure with generating sets that we denote by $I^{+}$ and $J^{+}$. The set $I^{+}$ can explicitly described as follows. For $m\ge 0$, let $F_{m}$ be the the left adjoint functor to the evaluation functor $Ev_{m}$ from symmetric spectra to simplicial sets, then
\[
I^{+}=\{F_{m}\partial\Delta[n]_{+}\to F_{m}\Delta[n]_{+} |m>0, n\ge 0\}.
\]
Note that the maps in $I^{+}$ have small domain and codomain. The set $J^{+}$ can be described in a similar fashion, all we need to know is that all the maps in $J^{+}$ also have small domain and codomain. According to this, the positive closed model structure on symmetric spectra has as weak equivalences the stable maps of symmetric spectra; that is, the maps $f:X\to Y$ such that $f^{*}:[Y,E]\to [X,E]$ is a bijection for all injective $\Omega$-spectrum $E$. A map of symmetric spectra is a cofibration if and only if it is a retract of a relative $I^{+}$-complex and a map is an acyclic cofibration if and only if it is a retract of a relative $J^{+}$-complex. The fibrations are the maps that satisfy the right lifting property with respect to the maps in $J^{+}$.

As mentioned before we want to rectify a multifunctor $T:E\Sigma_{*}\to \Sym^{\E}$ as to obtain a multifunctor $T':*\to \Sym^{\E}$. As our first step, we give a closed model category structure to the category of multifunctors $E\Sigma_{*}\to \Sym^{\E}$. To avoid any confusion we work in a general context. Thus let us fix $\M$ a multicategory enriched over $\cat$ which we can see as enriched over simplicial sets by taking the nerve. We will give a model structure to the category $\Mi(\M,\Sym^{\E})$ of multifunctors $\M\to \Sym^{\E}$. As our first step, we identify the categories $\Mi(\M,\Sym^{\E})$  
and $\Mi_{1}\Fi_{2}(\M,\E;\Sym)$ under the canonical isomorphism that sends a multifunctor $F:\M\to \Sym^{\E}$ to the assignment $F':\M\x \E\to \Sym$ defined by $F'(m,c)=F(m)(c)$. $F'$ is a multifunctor in the first variable and a functor in the second variable and as explain before, this defines an isomorphism of categories. Next, we give a model structure to the category $\Mi_{1}\Fi_{2}(\M,\E;\Sym)$.

Let us denote by $\O_{\M}$ (resp. $\O_{\E}$) the object set of $\M$ (resp. $\E$). On the product category $\Sym^{\O_{\M}\x\O_{\E}}$ we have a closed model structure for every closed model structure on $\Sym$, in particular, we have a model category on $\Sym^{\O_{\M}\x\O_{\E}}$ coming from the positive model structure on $\Sym$. We will show that the category $\Mi_{1}\Fi_{2}(\M,\E;\Sym)$ can be seen as the category of algebras of a suitable monad on $\Sym^{\O_{\M}\x\O_{\E}}$. Using this monad we will see that we can lift the positive model category on $\Sym^{\O_{\M}\x\O_{\E}}$ to the category $\Mi_{1}\Fi_{2}(\M,\E;\Sym)$. This way we will obtain a cofibrantly generated model category on $\Mi_{1}\Fi_{2}(\M,\E;\Sym)$.
\nin
To begin, consider the functor $\iota_{(a,c)}:\Sym \to \Sym^{\O_{\M}\x\O_{\C}}$ for each object $(a,c)\in \O_{\M}\x\O_{\E}$, that is defined on objects as follows. Take $X$ a symmetric spectrum, then
\begin{equation*}
(\iota _{(a,c)}X)_{(b,d)}=\left\{ 
\begin{array}{cc}
X\text{ } & \text{if }(a,c)=(b,d) \\ 
\ast \text{ } & \text{else.}%
\end{array}%
\right. 
\end{equation*}
Similarly, $\iota_{(a,c)}$ is defined on morphisms. The functor $\iota_{(a,c)}$ is precisely the left adjoint to the projection functor $\pi_{(a,c)}:\Sym^{\O_{\M}\x\O_{\E}}\to \Sym$. The positive stable model structure on $\Sym^{\O_{\M}\x\O_{\E}}$ is cofibrantly generated with generating sets $\iota I^{+}$ and $\iota J^{+}$, where 
\begin{align*}
\iota I^{+}&=\{\iota_{(a,c)}f |f\in I^{+}, (a,c)\in \O_{\M}\x\O_{\E} \},\\
\iota J^{+}&=\{\iota_{(a,c)}f |f\in J^{+}, (a,c)\in \O_{\M}\x\O_{\E} \}.
\end{align*}
Thus a map in $\Sym^{\O_{\M}\x\O_{\E}}$ is a cofibration if and only if it is a retract of a relative $\iota I^{+}$-complex and a map in $\Sym^{\O_{\M}\x\O_{\E}}$ is an acyclic cofibration if and only if it is a retract of a relative $\iota J^{+}$-complex.
\nin
Given $(a_{i},d)\in \O_{\M}\x \O_{\E}$ for $1\le i\le n$, we are going to use the notation 
\[
F_{(a_{1},...,a_{n};d)}=F_{(a_{1},d)}\wedge \cdots \wedge F_{(a_{n},d)}.
\]
With this in mind we have the following definition. 

\begin{definition}
For $(b,c)\in \O_{\M}\x\O_{\E}$ and $F:\O_{\M}\x\O_{\E}\to \Sym$, let 
\[
(\Do F)_{(b,c)}=\left.\V_{d\in \O_{\E}, \E(d,c)}\V_{n\ge 0}\left(\V_{a_{1},...,a_{n}\in \O_{\M}}\M(a_{1},...,a_{n};b)_{+}\wedge F_{(a_{1},...,a_{n};d)}\right)\right/\Sigma_{n}. 
\]
Let $\eta:F\to \Do F$ be the map such that for $(b,c)\in \O_{\M}\x\O_{\E}$
\[
\eta_{(b,c)}:F_{(b,c)}\cong {1_{b}}_{+}\wedge F_{(b,c)}\to \M(b;b)_{+}\wedge F_{(b,c)}\to (\Do F)_{(b,c)}
\]
and $\mu:\Do\Do F\to \Do F$ is the map induced by the multiproduct in $\M$ and the composition in $\E$. 
\end{definition}

By a direct and standard computation, it follows that $(\Do,\mu, \eta)$ is a simplicial monad on the category $\Sym^{\O_{\M}\x\O_{\E}}$.  Also, we see that a $\Do$-algebra structure on an object of $\Sym^{\O_{\M}\x\O_{\E}}$ is equivalent to an assignment $\M\x\E\to \Sym$ that is a multifunctor in the first variable and a functor in the second variable and that the simplicial category of $\Do$-algebras is isomorphic to $\Mi_{1}\Fi_{2}(\M,\E;\Sym)$. Finally, if we see $\Do$ as a functor 
\[
\Do:\Sym^{\O_{\M}\x\O_{\E}}\to \Mi_{1}\Fi_{2}(\M,\E;\Sym),
\]
then $\Do$ is left adjoint to the forgetful functor $\Mi_{1}\Fi_{2}(\M,\E;\Sym)\to \Sym^{\O_{\M}\x\O_{\E}}$.
\nin
From now on, if $a_{1},...a_{n},x,b\in \O_{\M}$ and for $k\ge 0$, we will denote the $(n+k)$-morphism set $\M(a_{1},...,a_{n},x,...,x,b)$ with $k$-copies of $x$ by $\M(a_{1},...,a_{n},x^{k},b)$. With this notation, note that in particular if $F=\iota _{(x,y)}X$, for some $(x,y)\in \O_{\M}\x\O_{\E}$ and a symmetric spectrum $X$, then for $(b,c)\in \O_{\M}\x\O_{\E}$ 
\[
(\Do \iota _{(x,y)}X)_{(b,c)}=\V_{\E(y,c)}\V_{n\ge 0} \M(x^{n};b)_{+}\wedge_{\Sigma_{n}} X^{(n)}
\]
One advantage of seeing $\Mi_{1}\Fi_{2}(\M,\E;\Sym)$ as the category of $\Do$-algebras, is the following lemma.

\begin{lemma}{\label{bicomplete}}
The category $\Mi_{1}\Fi_{2}(\M,\E;\Sym)$ is bicomplete.
\end{lemma}

\Proof
Since $\Mi_{1}\Fi_{2}(\M,\E;\Sym)$ is the category of algebras over a monad on a complete category, $\Mi_{1}\Fi_{2}(\M,\E;\Sym)$ is complete. On the other hand, it follows as in \cite[Proposition 7.2]{EKMM} that $\Do$ preserves reflexive coequalizers, thus by \cite[Proposition 7.4]{EKMM} it follows that $\Mi_{1}\Fi_{2}(\M,\E;\Sym)$ is cocomplete with colimits formed in $\Sym^{\O_{\M}\x\O_{\E}}$ as in \cite[Proposition 7.4]{EKMM}; that is, given a diagram ${F_{i}}$ of $\Do$-algebras, with structural maps $\xi_{i}:\Do F_{i}\to F_{i}$, let $\text{colim}_{i}F_{i}$ be its colimit in $\Sym^{\O_{\M}\x\O_{\E}}$ and let $f_{i}:F_{i}\to \text{colim} F_{i}$ be the natural maps. Also consider 
\[
\alpha:\text{colim}\Do F_{i}\to \Do \text{colim} F_{i}
\]
the unique map in $\Sym^{\O_{\M}\x\O_{\E}}$ whose composite is the natural map $\Do F_{i}\to \text{colim}_{i} \Do F_{i}$. Then the underlying object of the colimit of the $F_{i}$'s in the category $\Mi_{1}\Fi_{2}(\M,\E;\Sym)$ is the (reflexive) coequalizer in $\Sym^{\O_{\M}\x\O_{\E}}$ of
\[
\Do(\text{colim} \Do F_{i})\rightrightarrows \Do (\text{colim} F_{i}),
\] 
where one of the maps is $\Do (\text{colim} \xi_{i})$ and the other map is $\mu\circ \Do\alpha$.
\qed

We are now ready to give the weak equivalences, fibrations and cofibrations for the positive stable model category in $\Mi_{1}\Fi_{2}(\M,\E;\Sym)$.

\begin{definition}
\ \ \
\begin{itemize}
 \item A morphism in $\Mi_{1}\Fi_{2}(\M,\E;\Sym)$ is said to be a  weak equivalence if the underlying morphism in $\Sym^{\O_{\M}\x\O_{\E}}$ is a stable equivalence.
	\item A morphism in $\Mi_{1}\Fi_{2}(\M,\E;\Sym)$ is said to be a positive stable fibration if the underlying morphism in $\Sym^{\O_{\M}\x\O_{\E}}$ is a positive stable fibration.
	\item A morphism in $\Mi_{1}\Fi_{2}(\M,\E;\Sym)$ is said to be a positive stable cofibration if it satisfies the left lifting property with respect to the acyclic stable fibrations.
\end{itemize}
\end{definition}

In order to show that under these definitions we obtain a closed model structure on $\Mi_{1}\Fi_{2}(\M,\E;\Sym)$, we will show that this description is precisely the description of a cofibrantly generated model structure on $\Mi_{1}\Fi_{2}(\M,\E;\Sym)$ with generating sets $\Do \iota I^{+}$ and $\Do \iota J^{+}$, where 
\begin{align*}
\Do \iota I^{+}=&\{\Do\iota_{(a,c)}f | f\in I^{+}, (a,c)\in \O_{\M}\x\O_{\E}\},\\
\Do \iota J^{+}=&\{\Do\iota_{(a,c)}f | f\in J^{+}, (a,c)\in \O_{\M}\x\O_{\E}\}.
\end{align*}
To do so, we will first characterize the positive stable fibrations and acyclic positive fibrations in the proposition below. 

\begin{proposition}\label{fib}
\ \ \
\begin{itemize}
 \item A map in $\Mi_{1}\Fi_{2}(\M,\E;\Sym)$ is an acyclic positive fibration if and only if it has the right lifting property with respect to retracts of relative $\Do\iota I^{+}$-complexes; that is, the class of acyclic positive fibrations equals $\Do\iota I^{+}$-inj.
 \item A map in  $\Mi_{1}\Fi_{2}(\M,\E;\Sym)$ is a positive fibration if and only if it has the right lifting property with respect to retracts of relative $\Do\iota J^{+}$-complexes; that is, the class of positive fibrations equals $\Do\iota J^{+}$-inj.
\end{itemize}
\end{proposition}

\Proof
The result follows by the adjunction of $\Do$ and the forgetful functor. We show the first result, the second part follows in a similar fashion. As explained before a map in $\Sym^{\O_{\M}\x\O_{\E}}$ is a cofibration if and only if it is a retract of a relative $\iota I^{+}$-complex. Thus the positive stable acyclic fibrations on $\Sym^{\O_{\M}\x\O_{\E}}$ are precisely the maps that satisfy the right lifting property with respect to retracts of relative $\iota I^{+}$-complexes. Thus by the adjunction of $\Do$ and the forgetful functor, a map in $\Mi_{1}\Fi_{2}(\M,\E;\Sym)$ satisfies the right lifting property with respect retracts of  relative $\Do\iota I^{+}$-complexes if and only if the underlying map in $\Sym^{\O_{\M}\x\O_{\E}}$ satisfies the right lifting property with respect to retracts of relative $\iota I^{+}$-complexes.
\qed

The hardest part of showing that under the given definitions of weak equivalences, fibrations and cofibrations we obtain a closed model structure on $\Mi_{1}\Fi_{2}(\M,\E;\Sym)$ is the proposition below. Once we prove it, all the properties of a closed model structure are easy to verify and follow by applications of Quillen's small object argument and by the adjunction of $\Do$ and the forgetful functor.

\begin{proposition}\label{J^{+}}
A relative $\Do \iota J^{+}$-complex is a stable equivalence.
\end{proposition}

To prove this proposition we follow the same ideas used to prove \cite[Lemma 11.7]{Elmendorf}; that is, we study pushouts in the category $\Mi_{1}\Fi_{2}(\M,\E;\Sym)$ of the form 
\[
F\coprod_{\Do\iota_{(x,y)}X} \Do\iota_{(x,y)}Y,
\] 
for some $(x,y)\in \O_{\M}\x\O_{\E}$, $f:X\to Y$ morphism in  $\Sym$ and $g:\iota_{(x,y)}X\to F$ morphism in $\Sym^{\O_{\M}\x\O_{\E}}$. 
\nin
We begin with the following lemma that helps us study certain coproducts in the category $\Mi_{1}\Fi_{2}(\M,\E;\Sym)$.

\begin{lemma}
Let $\C$ be a cocomplete category and $(\T,\mu,\eta)$ a monad in $\C$ that $\T$ preserves reflexive coequalizers. Given a $\T$-algebra $X$ with structural map $\xi_{X}:\T X\to X$ and an object $Y$ of $\C$, the coproduct $X\coprod \T Y$ in the category of $\T$-algebras $\C[\T]$ is computed as the reflexive coequalizer in $\C$
\[
\T(\T X\coprod Y)\underset{\beta }{\overset{\alpha }{\rightrightarrows }} \T(X\coprod Y)\to P,
\]
where 
\begin{align*}
\alpha&=\mu \circ \T(\T(i_{1}\coprod \T i_{2}\circ \eta_{Y})),\\
\beta&=\T(\xi_{X}\coprod i_{2}).
\end{align*}
Moreover, $P$ is the coequalizer of the previous diagram in $\C[\T]$.
\end{lemma}

\Proof
By \cite[Lemma 6.6]{EKMM}, we have that $P$ is a $\T$-algebra and that $P$ is the coequalizer in $\C[\T]$ of the diagram
\[
\T(\T X\coprod Y)\underset{\beta }{\overset{\alpha }{\rightrightarrows }} \T(X\coprod Y)\to P.
\]
As $X$ is a $\T$-algebra, we see that $X$ is the coequalizer of the diagram
\[
\T\T X\underset{\T\xi_{X} }{\overset{\mu_{X} }{\rightrightarrows }} \T X\to X.
\]
Moreover the maps 
\begin{align*}
h&=\T(i_{1 \T X}):\T\T X\to \T(\T X\coprod Y)\\
k&=\T(i_{1 X}):\T X\to \T(X\coprod Y)
\end{align*}
are so that the following diagrams are commutative
\[
\xymatrix{
\T\T X\ar[r]^-{\mu_{X}}\ar[d]_-{h}   &\T X\ar[d]^-{k}& & \T\T X\ar[r]^-{\T\xi_{X}}\ar[d]_-{h}   &\T X\ar[d]^-{k}\\
\T(\T X\coprod Y)\ar[r]_{\alpha}         &\T(X\coprod Y),& & \T(\T X\coprod Y)\ar[r]_-{\beta}         &\T(X\coprod Y).
}
\]
By the universal property of the coequalizer, the previous diagrams guarantee the existence of a map $j_{1}:X\to P$. Moreover, as both $X$ and $P$ are the coequalizers in the category of 
$\T$-algebras, we see that this map is a map of $\T$-algebras. On the other hand, we have the map $j_{2}=\T(i_{2}):\T Y\to \T (X\coprod Y)$ which is also a map of $\T$-algebras. By a direct computation one can see that $P$ together with the maps $j_{1}$ and $j_{2}$ satisfy the universal property of the coproduct in the category $\C[\T]$. 
\qed

Note that we can apply the previous lemma to the monad $\Do$, as $\Sym^{\O_{\M}\x\O_{\E}}$ is bicomplete and $\Do$ preserves reflexive coequalizers. Thus if $X$ is a symmetric spectrum, $(x,y)\in \O_{\M}\x \O_{\E}$ and $F$ an object in $\Mi_{1}\Fi_{2}(\M,\E;\Sym)$, then the coproduct $F\coprod \Do\iota_{(x,y)} X$ in the category $\Mi_{1}\Fi_{2}(\M,\E;\Sym)$ is computed as the reflexive coequalizer  in the category $\Sym^{\O_{\M}\x\O_{\E}}$

\begin{equation}\label{coequi}
\Do(\Do F\coprod \iota_{(x,y)}X)\underset{\beta }{\overset{\alpha }{\rightrightarrows }} \Do(F\coprod \iota_{(x,y)} X).
\end{equation}

We can rewrite this coequalizer in a better way using the fact that the functor $-\wedge X$ preserves coequalizers. Thus define $T_{0}F=F$ and  for $k\ge 1$, define $T_{k}F$ in $\Sym^{\O_{\M}\x\O_{\E}}$ so that for  $(b,c)\in \O_{\M}\x \O_{\E}$, $(T_{k}F)_{(b,c)}$ is the coequalizer in $\Sym$ of the following diagram
\begin{align*}
&\left.\V_{\E(y,c)}\V_{n\ge 0}\left(\V_{a_{1},...,a_{n}\in \O_{\M}}\M(a_{1},...,a_{n},x^{k};b)_{+}\wedge \Do F_{(a_{1},...,a_{n};y)}\right)\right/\Sigma_{n}\\
&\left.\rightrightarrows \V_{\E(y,c)}\V_{n\ge 0}\left(\V_{a_{1},...,a_{n}\in \O_{\M}}\M(a_{1},...,a_{n},x^{k};b)_{+}\wedge F_{(a_{1},...,a_{n};y)}\right)\right/\Sigma_{n},
\end{align*}
where one of the maps is induced by the algebra structural map $\xi_{F}:\Do F\to F$ and the other is induced by the multiproduct in $\M$ and composition in $\E$. Using the $T_{k}F$'s, we can rewrite the coequalizer (\ref{coequi}) as to obtain that for a symmetric spectrum $X$,  $(x,y)\in\O_{\M}\x\O_{\E}$ and $F$ an object in $\Mi_{1}\Fi_{2}(\M,\E;\Sym)$, the coproduct $F\coprod \Do \iota_{(x,y)}$ in $\Mi_{1}\Fi_{2}(\M,\E;\Sym)$ is given by
\[
(F\coprod \Do \iota_{(x,y)}X)_{(b,c)}=\V_{k\ge 0}(T_{k}F)_{(b,c)}\wedge_{\Sigma_{k}} X^{k}. 
\]
We also use the $T_{k}$'s as defined above to study pushouts of the form 
\[
F\coprod_{\Do\iota_{(x,y)}X} \Do\iota_{(x,y)}X.
\]
To do so, we will use the same construction $Q^{k}_{i}(g)$, for $k\ge 0$, $0\le i\le k$ and a map of symmetric spectra $g:X\to Y$ as in \cite[Section 12]{Elmendorf}. We include it for completeness. The $Q^{k}_{i}(g)$'s are inductively defined as follows: $Q^{k}_{0}(g)=X^{k}$, $Q^{k}_{k}(g)=Y^{k}$. For $0<i<k$, define $Q^{k}_{i}(g)$ as the pushout square
\[
\xymatrix{
\Sigma_{k+}\wedge_{\Sigma_{k-i}\x\Sigma_{i}}X^{k-i}\wedge Q^{i}_{i-1}(g)\ar[r]\ar[d]   &\Sigma_{k+}\wedge_{\Sigma_{k-i}\x\Sigma_{i}}X^{k-i}\wedge Y^{i}\ar[d]\\
Q^{k}_{i-1}(g)\ar[r]         &Q^{k}_{i}(g)
}
\]
The $Q_{i}^{k}$'s are defined so that essentially, $Q_{i}^{k}$ is the $\Sigma_{k}$-subspectrum of $Y^{k}$ with $i$ factors of $Y$ and $k-i$ factors of $X$. We are now ready to describe our final construction. For an object  $F$ of $\Mi_{1}\Fi_{2}(\M,\E;\Sym)$, let $F_{0}=F$ and for $k>0$ let $F_{k}$ to be the pushout in $\Sym^{\O_{\M}\x\O_{\E}}$ 
\[
\xymatrix{
T_{k}F\wedge_{\Sigma_{k}}Q^{k}_{k-1}(g)\ar[r]\ar[d]   &T_{k}F\wedge_{\Sigma_{k}}Y^{k}\ar[d]\\
F_{k-1}\ar[r]         &F_{k}.
}
\]
Where the map in the top is the map induced by $Q^{k}_{k-1}(g)\to Y^{k}$ and the map on the left is induced by $\iota_{(x,y)X}\to F$. Let $F_{\infty}=\text{colim}F_{k}$ where the colimit is computed on $\Sym^{\O_{\M}\x\O_{\E}}$. The relevance of the $F_{k}$'s is given in the following proposition which corresponds to \cite[Proposition 12.6]{Elmendorf}.

\begin{proposition}\label{filtration}
$F_{\infty}$ is isomorphic to the underlying object in $\Sym^{\O_{\M}\x\O_{\E}}$ of the pushout $\coprod_{\Do\iota_{(b,c)}X}\Do\iota_{(b,c)}Y$.
\end{proposition}

\Proof
The proof goes by showing that $F_{\infty}$ satisfies the universal property of the pushout 
\[
F\coprod_{\Do\iota_{(b,c)}X}\Do\iota_{(b,c)}Y
\] 
in the category $\Mi_{1}\Fi_{2}(\M,\E;\Sym)$.
\qed

With the given filtration $F=F_{0}\to F_{1}\to \cdots F_{\infty}$ of $F\coprod_{\Do\iota_{(b,c)}X}\Do\iota_{(b,c)}Y$ we can prove proposition \ref{J^{+}}. 

\noindent{\textit{Proof of proposition \ref{J^{+}}}.} 
We want to show that if $h:X\to Y$ is a relative $\Do \iota J^{+}$-complex then $h$ is a positive stable equivalence. We are going to consider first the special case where $h$ is obtained as a pushout $h:F\to F\coprod_{\Do\iota_{(x,y)}X}\Do\iota_{(x,y)}Y$, for a map $f:X\to Y$ of symmetric spectra in $J^{+}$ and $g:\iota_{(x,y)}X\to F$  a morphism in $\Sym^{\O_{\M}\x \O_{\E}}$. We have a filtration $F_{0}=F\to F_{1}\to \cdots\to F_{\infty}$ as in the previous proposition such that $h:F=F_{0}\to F_{\infty}$ is the transfinite composition of the maps $F_{k}\to F_{k+1}$. To prove that $h$ is a positive stable equivalence it's enough to prove then that each map $F_{k}\to F_{k+1}$ is a positive stable equivalence as the composition of (transfinite) acyclic positive cofibrations is again an acyclic positive cofibation in the positive stable model category. By \cite[Proposition 12.6]{Elmendorf} we have that $F_{k}\to F_{k+1}$ is an objectwise level cofibration of symmetric spectra.  Also note that the quotient $F_{k+1}/F_{k}$ is naturally isomorphic to $T_{k+1}F\wedge_{\Sigma_{k+1}}(Y/X)^{k+1}$. As $X\to Y$ is an acyclic positive cofibration it follows that $Y/X$ is positive cofibrant and stably equivalent to $*$. Hence $F_{k+1}/F_{k}$ is stably equivalent to $*$ in $\Sym^{\O_{\M}\x\O_{\E}}$. But as mentioned before $F_{k}\to F_{k+1}$ is an objectwise level cofibration and then we conclude that $F_{k}\to F_{k+1}$ is a stable equivalence. To finish the proof, note that we can write any relative $\Do \iota J^{+}$-complex as a retract of a possible transfinite composition of pushouts of the form we just considered. These are stable equivalences by the previous consideration. On the positive model structure a transfinite composition of stable equivalences is also a stable equivalence. This finishes the proof.
\qed

With the proof of the Proposition \ref{J^{+}} completed we can now move on to proving that we in fact have defined a model structure on $\Mi_{1}\Fi_{2}(\M,\E;\Sym)$. We start by constructing factorizations of morphisms in $\Mi_{1}\Fi_{2}(\M,\E;\Sym)$. To do so, note that as mentioned before, the domains of the maps in $I^{+}$ and $J^{+}$ are small, this implies that the domains of the maps in $\Do\iota I^{+}, \Do\iota J^{+}$ are also small. Therefore we can apply Quillen's small object argument as in \cite[Theorem 2.1.14]{Hov} to factorize a morphism $f$ in $\Mi_{1}\Fi_{2}(\M,\E;\Sym)$ in the form $f=pi$, where
\begin{itemize}
  \item  $i$ is in $\Do \iota I^{+}$-cell and $p$ is in $\Do \iota I^{+}$-inj. or
	\item  $i$ is in $\Do \iota J^{+}$-cell and $p$ is in $\Do \iota J^{+}$-inj.
\end{itemize}
By Proposition \ref{fib} a map in $\Do \iota I^{+}$-inj is precisely an acyclic positive fibration and a map in $\Do \iota J^{+}$-inj is a positive fibration. Thus we obtain the following proposition. 

\begin{proposition}{\label{Quillen}}
\ \ \ 
\begin{itemize}
  \item A map $f$ in $\Mi_{1}\Fi_{2}(\M,\E;\Sym)$  can be factored as $f=pi$, where $i$ is a relative $\Do \iota I^{+}$-complex and $p$ is an acyclic positive stable fibration. 
	\item A map $f$ in $\Mi_{1}\Fi_{2}(\M,\E;\Sym)$  can be factored as $f=pi$, where $i$ is a $\Do \iota J^{+}$-complex and $p$ is a positive stable fibration. 
\end{itemize}
\end{proposition}

In order for these to be the required factorizations we need to prove the following lemma.

\begin{lemma}\label{cofibrations}
\ \ \ 
\begin{itemize}
\item A map in $\Mi_{1}\Fi_{2}(\M,\E;\Sym)$ is a positive stable cofibration if and only if it is a retract of a relative $\Do \iota I^{+}$-complex.
\item A map in $\Mi_{1}\Fi_{2}(\M,\E;\Sym)$ is an  acyclic positive stable cofibration if and only if it is a retract of a relative $\Do \iota J^{+}$-complex.
\end{itemize}
\end{lemma}
\Proof
Let's prove the first assertion. Any retract of a relative $\Do\iota I^{+}$-complex satisfies the left lifting property with respect to the class $\Do\iota I^{+}$-inj. Thus by Proposition \ref{fib} and the definition any retract of a relative $\Do\iota I^{+}$-complex is a positive stable cofibration. On the other hand, take $f$ a map that is a positive stable cofibration. By the Proposition \ref{Quillen}, we can factor  $f=pi$, where $i$ is a relative $\Do \iota I^{+}$-complex followed and $p$ is an acyclic positive stable fibration. But since $f$ is a cofibration, by definition it follows that $f$ satisfies the left lifting property with respect to $p$. By \cite[Lemma 1.19]{Hov} we have that $f$ is a retract of $i$ and thus $f$ is a retract of a relative $\Do \iota I^{+}$-complex.

The second assertion is proved in the same way. If $f$ is a retract of a relative $\Do \iota J^{+}$-complex then by Proposition \ref{J^{+}} we have that $f$ is a stable equivalence. Also, $f$ satisfies the left lifting property with respect to the class $\Do\iota J^{+}$-inj. and thus by Proposition \ref{fib}, $f$ satisfies the left lifting property with respect to the class of positive stable fibrations, in particular $f$ is a positive stable cofibration. On the other hand, take $f$ a map that is an acyclic stable cofibration. By the Proposition \ref{Quillen}, we can factor  $f=pi$, where $i$ is a relative $\Do \iota J^{+}$-complex followed and $p$ is a positive stable fibration. By Proposition \ref{J^{+}}, $i$ is a stable equivalence and so is $f$. Since stable equivalences satisfy the two out of three property we see that $p$ is also a stable equivalence, hence $p$ is an acyclic stable fibration. In particular $f$ satisfies the left lifting property with respect to $p$ and by \cite[Lemma 1.19]{Hov} we have that $f$ is a retract of $i$ and thus $f$ is a retract of a relative $\Do \iota J^{+}$-complex.
\qed

Combining Propositions \ref{Quillen} and \ref{cofibrations} we obtain the desired factorizations as in the following lemma.

\begin{lemma}\label{factorizations}
A map $f$ in $\Mi_{1}\Fi_{2}(\M,\E;\Sym)$ can be factored as $f=pi$, where $i$ is a stable cofibration and $p$ is an acyclic positive stable fibration or we can take $i$ is an acyclic cofibration and $p$ is a positive stable fibration. 
\end{lemma}

After proving these factorizations we are now ready to show that this choice of weak equivalences, fibrations and cofibrations gives rise a closed model category structure on the category $\Mi_{1}\Fi_{2}(\M,\E;\Sym)$.

\begin{theorem}
There exists a model category structure on $\Mi_{1}\Fi_{2}(\M,\E;\Sym)$ whose weak equivalences are the morphisms that are stable equivalences, the fibrations are the positive stable fibrations and the cofibrations are the positive stable cofibrations.
\end{theorem}
\Proof
We need to verify axioms MC1.-MC5. of \cite[Definition 3.3]{Dwyer}. We need to note first that it is clear by the definition that the classes of weak equivalences, fibrations and cofibrations contain the identity and are closed under composition.
\begin{itemize}
	\item MC1. By Proposition \ref{bicomplete} we know that $\Mi_{1}\Fi_{2}(\M,\E;\Sym)$ has all small limits and colimits in particular those that are finite.
	\item MC2. Since in $\Sym^{\O_{\M}\x\O_{\E}}$ the two out of three property is satisfied, it follows at once that it is satisfied on $\Mi_{1}\Fi_{2}(\M,\E;\Sym)$ 
	\item MC3. By the same reason as in MC2. retracts of weak equivalences are weak equivalences. Also by the characterizations of Propositions \ref{fib}, \ref{cofibrations} we see that fibrations and cofibrations are closed under retracts.
	\item MC4. The cofibrations satisfy the left lifting property with respect to acyclic fibrations by definition. On the other hand by Propositions \ref{fib}, \ref{cofibrations} we have that fibrations satisfy the right lifting property with respect to acyclic fibrations.
	\item MC5. The factorization properties follow by Lemma \ref{factorizations}.
\end{itemize}
\qed

Our next step toward rectifying a multifunctor $T:E\Sigma_{*}\to \Sym^{\E}$ is to study the relationships between the previous model structure when we change the multicategory $\M$. To be more precise, we want to show that if $f:\M\to \M'$ is a multifunctor between multicategories enriched over simplicial sets, then 
\[
g^{*}:=(f\x id)^{*}:\Mi_{1}\Fi_{2}(\M',\E;\Sym)\to \Mi_{1}\Fi_{2}(\M,\E;\Sym)
\]
is a right adjoint in a Quillen adjunction. Also we want to show that in the case that $f$ is a weak equivalence between multicategories, then this Quillen adjunction is actually a Quillen equivalence. We will apply this to the particular case where $\M=E\Sigma_{*}$ and $\M'=*$ as to get the desired rectification of the multifunctor $T:E\Sigma_{*}\to \Sym^{\E}$. Once again we will follow very closely those ideas presented on \cite{Elmendorf} with only slight modifications as to adapt them to our situation. 

Suppose we have a multifunctor $f:\M\to \M'$ between two small multicategories, let $\O$ and $\O'$ be the object sets of $\M$ and $\M'$. Also let $\Do$ and $\Do'$ be the monads on $\Sym^{\O_{\M}\x\O_{\E}}$ and $\Sym^{\O_{\M'}\x\O_{\E}}$, induced by $\M$ and $\M'$ respectively.  By composing with $g:=f\x \text{id}$, we obtain a functor between the product categories
\[
\pi_{g}: \Sym^{\O_{\M'}\x\O_{\E}}\to \Sym^{\O_{\M}\x\O_{\E}}.
\]
This functor has a left adjoint which we denote by $\kappa_{g}$,
\[
\kappa_{g}:\Sym^{\O_{\M}\x\O_{\E}}\to \Sym^{\O_{\M'}\x\O_{\E}}.
\]
An explicit formula for $\kappa_{g}$ is given as follows: if $F:\O_{\M}\x \O_{\E}\to \Sym$, then
\[
(\kappa_{g}F)_{(b,c)}=\V_{a\in f^{-1}(b)}F_{(a,c).}
\]
Note that we can see the multifunctor $f$ as a natural transformation 
\[
\kappa_{g}\Do\to \Do'\kappa_{g}.
\]
With this definition we can prove the following theorem.

\begin{theorem}{\label{QuillenAd}}
Let $\M$ and $\M'$ be  multicategories enriched over simplicial sets and $f:\M\to \M'$. Then there exists a functor 
\[
g_{\#}:\Mi_{1}\Fi_{2}(\M,\E;\Sym)\to \Mi_{1}\Fi_{2}(\M',\E;\Sym)
\]
left adjoint to $g^{*}$ such that the pair $(g_{\#},g^{*})$ forms a Quillen adjunction.
\end{theorem}

\Proof
For an object $F$ of $\Mi_{1}\Fi_{2}(\M,\E;\Sym)$ define $g_{\#}F$ as the reflexive coequalizer of the following diagram
\[
\Do'\kappa_{g}\Do F \rightrightarrows \Do'\kappa_{g}F\to g_{\#}F,
\]
where one of the arrows comes from the multifunctor $f$ seen as a natural transformation $\kappa_{g}\Do\to \Do'\kappa_{g}$ and the other arrow comes from the algebra structure  map $\xi{F}:\Do F\to F$ of $F$. It follows at once by the universal property of this coequalizer and the adjunctions that this way defined $g_{\#}$ is left adjoint to 
\[
g^{*}:\Mi_{1}\Fi_{2}(\M',\E;\Sym)\to \Mi_{1}\Fi_{2}(\M,\E;\Sym).
\]
On the other hand, it follows by the definition that $g^{*}$ preserves fibrations and acyclic fibrations. Thus by \cite[Theorem 9.7]{Dwyer} we get that the pair $(g_{\#},g^{*})$ forms a Quillen adjunction.
\qed

We will see now that in the case that $f$ is a weak equivalences of multicategories then this Quillen adjunction actually forms a Quillen equivalence.
We begin by recalling the definition of a weak equivalence between multicategories from \cite[Section 12]{Elmendorf}.

\begin{definition}   
Let $\M$ and $\M'$ be two multicategories enriched over simplicial sets and $f:\M\to \M'$ an (enriched) multifunctor. We say that $f$ is a weak equivalence whenever the functor $\pi_{0}f$ is an equivalence between the category of components and for all objects $a_{1},...,a_{n},b$ of $\M$, the map of simplicial sets 
\[
\M(a_{1},...,a_{n},b)\to \M'(f(a_{1}),...,f(a_{n}),f(b))
\]
is a weak equivalence.
\end{definition}

We show now that in the case where $f:\M\to \M'$ is a weak equivalence then  the pair $(g_{\#},g^{*})$ forms a Quillen equivalence. 

\begin{theorem}{\label{Equivalence}}
If $f:\M\to \M'$ is a weak equivalence between multicategories enriched over simplicial sets, then the Quillen adjunction $(g_{\#},g^{*},)$ forms a Quillen equivalence.
\end{theorem}

\Proof
By \cite[Theorem 9.7]{Dwyer}, we need to show that for every cofibrant object $A$ of $\Mi_{1}\Fi_{2}(\M,\E;\Sym)$ and every fibrant object $B$ of $\Sym^{\M'\x \E}$, a map $h:g_{\#}A\to B$ is a weak equivalence if and only if its adjoint $j:A\to g^{*}B$ is a weak equivalence. By the Lemma \ref{weak} below, we know that $h$ is a weak equivalence if and only if $g^{*}h$ is a weak equivalence. By definition, $j$ equals the composite 
\[
A\stackrel{\psi}{\rightarrow} g^{*}g_{\#}A\stackrel{g^{*}h}{\rightarrow} g^{*}B
\]
where $\psi:A\to g^{*}g_{\#}A$ is the unit of the adjunction. We show below in Theorem \ref{unit} that whenever $A$ is cofibrant $\psi$ is positive stable equivalence. By the two out of three axiom we see that $h:A\to g^{*}B$ is a weak equivalence if and only if its adjoint $j:g_{\#}A\to B$ is a weak equivalence
\qed

\begin{lemma}{\label{weak}}
A map $h:F\to G$ in $\Mi_{1}\Fi_{2}(\M',\E;\Sym)$ is a weak equivalence if and only if $g^{*}h$ is a weak equivalence in $\Mi_{1}\Fi_{2}(\M,\E;\Sym)$.
\end{lemma}
\Proof
The proof is a straight forward generalization of the proof of \cite[Lemma12.4]{Elmendorf}.
\qed

Note that the initial object of $\Mi_{1}\Fi_{2}(\M,\E;\Sym)$ is the assignment 
\begin{align*}
\M\x\E&\to \Sym\\
(b,c)&\to \M(;b)_{+}\wedge S.
\end{align*}
\begin{definition}
Let $\lambda$ an ordinal and $F$ a cofibrant object of $\Mi_{1}\Fi_{2}(\M,\E;\Sym)$. We say that $F$ can be built in $\lambda$-stages, 
if we can find a $\lambda$-sequence 
\[
X:\lambda\to \Mi_{1}\Fi_{2}(\M,\E;\Sym)
\]
such that $X_{0}=\M(;-)_{+}\wedge S$, the initial object, and for each $\beta$ with $\beta+1<\lambda$ there is a pushout
\[
\xymatrix{
C_{\beta}\ar[r]\ar[d]_{g_{\beta}}   &X_{\beta}\ar[d]\\
D_{\beta}\ar[r]         &X_{\beta+1},
}
\]
such that $g_{\beta}$ is the coproduct of maps in $\Do\iota I^{+}$ and that $X_{0}\to F$ is isomorphic to the transfinite composition $X_{0}\to \textsl{colim}_{\beta<\lambda} X_{\beta}$. 
\end{definition}

Note that if $\beta\le \lambda$ and $A$ can be built in $\beta$-stages, then clearly $A$ can be built in $\lambda$-stages. We are now ready to prove the following theorem.

\begin{theorem}{\label{unit}}
If $A$ is a cofibrant object in $\Mi_{1}\Fi_{2}(\M,\E;\Sym)$, then the unit 
\[
\psi:A\to g^{*}g_{\#}A 
\] 
of the Quillen adjunction $(g_{\#},g^{*})$ is a stable equivalence.
\end{theorem}

\Proof
The proof of this theorem is taken from the proof of \cite[Theorem 12.5]{Elmendorf} with small modifications as to fit into our situation. Let $A$ be a cofibrant object in $\Mi_{1}\Fi_{2}(\M,\E;\Sym)$. Then $\M(;-)_{+}\wedge S\to A$ is a retract of a relative $\Do\iota I^{+}$-complex. It suffices then to prove the statement for the case where $A$ can be built in $\lambda$-stages for an ordinal $\lambda$. The proof goes by transfinite induction on the stages that $A$ can be built in. Thus we want to show that for every ordinal $\lambda$ and every cofibrant object that can be build in $\lambda$-stages, the map 
\[
\psi:A\to g^{*}g_{\#}A 
\] 
is a stable equivalence. For $\lambda=0$, we have that a cofibrant object that can be build in 0-stages is isomorphic to $\M(;-)_{+}\wedge S$. In this case
\begin{align*}
g^{*}g_{\#}(\M(;-)_{+}\wedge S)&=\M'(;-)_{+}\wedge S,\\
A=\M(;-)_{+}\wedge S&\to \M'(;-)_{+}\wedge S=g^{*}g_{\#}A
\end{align*}
and the statement follows as $f$ is weak equivalence. 

We also need the case $\lambda=1$. In this case, a cofibrant object that can be built in 1-stage is an object of the form $A=\Do F$, where $F$ is an object in $\Sym^{\O_{\M}\x \O_{\E}}$ that is objectwise cofibrant. In this situation the map $A\to g^{*}g_{\#}A$ equals
\[
\Do X\to g^{*}\Do\kappa_{g}X.
\]
By definition we see that for $(b,c)\in \O_{\M}\x\O_{\E}$
\begin{align*}
&(g^{*}\Do\kappa_{g}X)_{(b,c)}=\\
&\left.\V_{d\in \O_{\E},\E(d,c), n\ge 0}\left(\V_{a'_{1},...,a'_{n}\in \O_{\M'}}\V_{i=1}^{n}\V_{a_{i}\in f^{-1}(a'_{i})}(\M'_{a'_{1},...,a'_{n};f(b)}){+}\wedge F_{(a_{1},...,a_{n};d)}\right)\right/\Sigma_{n}. 
\end{align*}
Where here  $(\M'_{a'_{1},...,a'_{n};f(b)})$ means $\M'(a'_{1},...,a'_{n};f(b))$. Since $f$ is a weak equivalence, it follows that $\Do X\to g^{*}\Do\kappa_{g}X$ is a stable equivalence and thus the statement is true for $\lambda=1$. Suppose the statement is true for all ordinals $\beta<\lambda$ and that $\lambda$ is a limit ordinal. Take a $\lambda$-sequence $X:\lambda\to \Mi_{1}\Fi_{2}(\M,\E;\Sym)$ with $X_{0}=\M(;-)_{+}\wedge S$ and with $\M(;-)_{+}\wedge S\to \text{colim}_{\beta<\lambda}X_{\beta}$ isomorphic to $X_{0}\to A$. We have that 
\begin{equation}\label{colimX}
g^{*}g_{\#}X_{\beta}=\text{colim}_{\beta<\lambda}g^{*}g_{\#}X_{\beta}
\end{equation} 
and since each $X_{\beta}$ can be build in $\beta$-stages and $\beta< \lambda$, the statement is true for $X_{0}\to X_{\beta}$ and by (\ref{colimX}) we see that the statement is true for $X_{\lambda}$ and hence for $A$. Thus we are left to prove the statement for $\beta+1$ provided it's true for $\beta$ and $\beta\ge 1$. By a separate  transfinite induction, this case reduces to showing that if $A$ is a cofibrant object with $A\to g^{*}g_{\#}A$ a stable equivalence then the same is true for $B$, where $B$ is obtained as the pushout in $\Mi_{1}\Fi_{2}(\M,\E;\Sym)$ of the following diagram
\[
\xymatrix{
\Do\iota_{(x,y)} X\ar[r]\ar[d]   &A\ar[d]\\
\Do \iota_{(x,y)}Y\ar[r]         &B.
}
\]
By Proposition \ref{filtration}, we have filtration $A=A_{0}\to A_{1}\to \cdots \to A_{\infty}=B$, where each map $A_{k}\to A_{k+1}$ is an objectwise level cofibration. The associated graded to this filtration is 
\[
\V_{k\ge0}T_{k}\wedge_{\Sigma_{k}}(Y/X)^{k}.
\]
Note that this is isomorphic in $\Sym^{\O_{\M}\x \O_{\E}}$ to the underlying object of 
\[
A\coprod \Do\iota_{(x,y)}(Y/X).
\] 
Let $A'=g_{\#}A$ and $B'=g_{\#}B$. Note that 
\[
B'=A'\coprod_{\Do'\iota_{(f(x),y)}X} \Do'\iota_{(f(x),y)}X.
\]
For $B'$ we also have a filtration $A'=A'_{0}\to A'_{1}\to \cdots \to A'_{\infty}=B'$, where each map $A'_{k}\to A'_{k+1}$ is an objectwise level cofibration and whose associated graded is isomorphic in $\Sym^{\O_{\M'}\x \O_{\E}}$ to $A'\coprod \Do'\iota_{(f(x),y)}(Y/X)$. The map  $B\to g^{*}B'=\pi_{g}B'$ preserves the given filtrations. Also the map of associated graded is
\[
A\coprod \Do\iota_{(x,y)}(Y/X)\to \pi_{g}(A'\coprod \Do'\iota_{(f(x),y)}(Y/X))\approx g^{*}g_{\#}(A\coprod \Do\iota_{(x,y)}(Y/X)),
\]
and we that this is a stable equivalence. This follows from the fact that the cofibrant object $A\coprod \Do\iota_{(x,y)}(Y/X)$ can be built in $\beta$ stages. Indeed, 
suppose $W:\beta\to \Mi_{1}\Fi_{2}(\M,\E;\Sym)$ is a $\beta$-sequence such that $W_{0}=\M(;-)_{+}\wedge S$ and for each $\gamma$ with $\gamma+1<\beta$ there is a pushout
\[
\xymatrix{
C_{\gamma}\ar[r]\ar[d]_{g_{\gamma}}   &W_{\gamma}\ar[d]\\
D_{\gamma}\ar[r]         &W_{\gamma+1}
}
\]
such that $g_{\gamma}$ is the coproduct of maps in $\Do\iota I^{+}$ and that $W_{0}\to A$ is the transfinite composition $W_{0}\to \text{colim}_{\gamma<\beta} W_{\gamma}$. Let $W':\beta\to \Mi_{1}\Fi_{2}(\M,\E;\Sym)$ be the $\beta$-sequence defined as follows: Take $W'_{0}=\M(;-)_{+}\wedge S$, define $W'_{1}$ as the pushout of the diagram
\[
\xymatrix{
C_{0}\coprod \M(;-)_{+}\wedge S\ar[r]\ar[d]_{g_{0}\coprod i}   &W'_{0}\ar[d]\\
D_{0}\coprod \Do \iota_{(x,y)}(Y/X)\ar[r]         &W'_{1}.
}
\]
Thus $W'_{1}=W_{1}\coprod \Do \iota_{(x,y)}(Y/X)$. In general define for $\gamma\ge 2$ 
\[
W_{\gamma}'=W_{\gamma}\coprod \Do \iota_{(x,y)}(Y/X).
\] 
This way defined we see that for all $\gamma<\beta$, $W'_{\gamma+1}$ is the pushout of  
\[
\xymatrix{
C_{\gamma}\ar[r]\ar[d]_{g_{\gamma}}   &W'_{\gamma}\ar[d]\\
D_{\gamma}\ar[r]         &W'_{\gamma+1}.
}
\]
By the induction hypothesis, it follows that 
\[
A\coprod \Do\iota_{(x,y)}(Y/X)\to g^{*}g_{\#}(A\coprod \Do\iota_{(x,y)}(Y/X))
\] 
is a stable equivalence. Since each map $A_{k}\to A_{k+1}$, and $A'_{k}\to A_{k+1}$ are objectwise level cofibrations and the map of associated graded is a stable equivalence we have that each map $A{k}\to \pi_{g}A_{k}$ is a stable equivalence. From here we conclude that $B\to \pi_{g}B'\approx g^{*}g_{\#}B$ is a stable equivalence. 
\qed

\begin{corollary}\label{rectifying}
A multifunctor $T:E\Sigma_{*}\to \Sym^{\E}$ gives rise to a multifunctor $T':*\to \Sym^{\E}$
\end{corollary}
\Proof
Consider the multifunctor $f:E\Sigma_{*}\to *$ that sends every $k$-morphism in $E\Sigma_{k}$ to the only $k$-morphism in $*$. This defines trivially a multifunctor between multicategories enriched over simplicial sets. (We see $E\Sigma_{*}$ enriched over simplicial sets by taking the nerve of the category $E\Sigma_{k}$, for $k\ge 0$). Note that for every $k\ge 0$, the geometric realization of $E\Sigma_{k}$ is contractible and thus $f$ is a weak equivalence between the multicategories $\M=E\Sigma_{*}$ and $\M=*$. By Theorem \ref{Equivalence}, we have that the functor $g_{\#}:=(f\x\text{id})_{\#}$ is a left adjoint in a Quillen equivalence. Let $T'=g_{\#}(T)$, then $T':*\to \Sym^{\E}$ is the required multifunctor.
\qed

We show in the following proposition that such a multifunctor structure is equivalent as having a lax map $F:\E\to \Sym$.

\begin{proposition}\label{laxmul}
Having a multifunctor $L:*\to \Sym^{\E}$ is equivalent as having a lax map $F:\E\to \Sym$
\end{proposition}
\Proof
Suppose we have a multifunctor $L:*\to \Sym^{\E}$. The multicategory $*$ has only one object so we consider $F:\E\to \Sym$ the image under $L$ of this object. Thus  $F:\E\to \Sym$ is a functor. On $*$ we have a unique $2$-morphism, let $\phi\in \Sym^{\E}(F,F;F)$, be the image of this morphism. Unraveling the definitions we see that $\phi$ is a $2$-linear natural transformation, thus for every objects $c_{1},c_{2}$ of $\E$ we have a map of symmetric spectra
\[
\phi_{c_{1},c_{2}}:F(c_{1})\wedge F(c_{2})\to F(c_{1}\ot c_{2}).
\]
The collection of maps $\{\phi_{c_{1},c_{2}}\}$ satisfy the naturality condition as in Definition \ref{naturality}. Suppose that $f:c_{1}\to c_{1}'$, $g:c_{2}\to c_{2}'$ are morphisms in $\E$. Then the following naturality condition must be satisfied
\[
\xymatrix{
F(c_{1})\wedge F(c_{2})\ \ \ar[r]^-{F(f)\wedge F(g)}\ar[d]_-{\phi_{c_{1},c_{2}}}   &\ \ F(c_{1}')\wedge F(c_{2}')\ar[d]^-{\phi_{c'_{1},c'_{2}}}\\
F(c_{1}\ot c_{2})\ar[r]_-{F(f\ot g)}         &F(c_{1}'\ot c_{2}')
}
\]
This means that the map $\phi_{c_{1},c_{2}}:F(c_{1})\wedge F(c_{2})\to F(c_{1}\ot c_{2})$ is a natural map. On the other hand, $*$ has a unique $0$-morphism, let $\eta\in \Sym^{\E}$ be the image of this unique $0$-morphism. Then $\eta:S\to F(1)$, where here $S$ is the symmetric sphere spectrum. 
The maps $\eta,\phi$ satisfy some coherences coming from the fact that $L$ preserves the multifunctor structure, for example if $\sigma\in \Sigma_{2}$ is the nontrivial element, then $\phi=\sigma^{*}\phi$, this means that for every objects $c_{1},c_{2}$ in $\E$, the following diagram is commutative
\[
\xymatrix{
F(c_{1})\wedge F(c_{2})\ar[r]^-{\alpha}\ar[d]_-{\phi_{c_{1},c_{2}}}   &F(c_{2})\wedge F(c_{1})\ar[d]^-{\phi_{c_{2},c_{1}}}\\
F(c_{1}\ot c_{2})\ar[r]_-{F(\gamma)}         &F(c_{2}\ot c_{1})
}
\]
here $\alpha$ is the natural isomorphism in symmetric spectra that changes the factors and $\gamma$ is the natural isomorphism coming from the permutative structure on  $\E$. The other coherences for a lax map are obtained in a similar way.

Conversely, suppose that $F:\E\to \Sym$ is a lax map. Thus we have a natural maps $\eta:S\to F(1)$ and $\phi_{c_{1},c_{2}}:F(c_{1})\wedge F(c_{2})\to F(c_{1}\ot c_{2})$. Then using the coherences that these maps satisfy we see that the assignment
\[
L:*\to \Sym^{\E}
\] 
that sends the unique object of $*$ to $F$, the unique $0$-morphism to $\eta$ and the unique $2$-morphism to $\phi$ defines a multifunctor.
\qed

According to Corollary \ref{rectifying} and Proposition \ref{laxmul}, a multifunctor $T:E\Sigma_{*}\to \Sym^{\E}$ can be rectified as to obtain a lax map $\vartheta:\E\to \Sym$. In particular, the multifunctor $T:E\Sigma_{*}\to \Sym^{\A}$ induced by a bipermutative category of fibers $\Lambda:\D\to \C$, induces a lax map $\vartheta:\A\to \Sym$. We show now that this lax map induces a lax map $\phi:\C^{op}\to \Sym$. We show this in the following corollary.

\begin{corollary}
A fibered bipermutative category $\Lambda:\D\to \C$ induces a lax map $\phi:\C^{op}\to \Sym$.
\end{corollary}
\Proof
We have already seen that $\Lambda:\D\to \C$ a lax map $\vartheta:\A\to \Sym$. Recall that the category $\A$ has as objects the sequences $\u=(u_{1},...,u_{n})$ and as morphisms the tuples $(q,\f):\u\to \v$, where $q:\n\to \m$ is an injection and $\f:q_{*}\u\to \v$ is a morphism in $(\C^{op})^{m}$. We have a functor $\W:\A\to \C^{op}$ defined as follows. Given an object $\u$ of $\A$, then 
\[
\W(\u)=u_{1}\ot\cdots \ot u_{n}.
\]
If $(q,\f):\u\to \v$ is a morphism in $\A$, then $\f:q_{*}\u\to \v$ is a morphism in $(\C^{op})^{m}$, thus $f_{i}:u_{i}'\to v_{i}$ is a morphism in $\C^{op}$, where $u_{i}'$ is either $1$ or one of the $u_{j}$'s. In particular note that $\ot_{i=1}^{m}u_{i}'=\ot_{j=1}^{n}u_{\sigma(j)}$ for some permutation $\sigma\in \Sigma_{n}$. Define $\W((\sigma,\f)):\W(\u)\to \W(\v)$ to be the following composite in $\C^{op}$
\[
u_{1}\ot\cdots \ot u_{n}\stackrel{\tau_{\sigma,u_{1},...,u_{n}}}{\rightarrow}u_{\sigma(1)}\ot\cdots \ot u_{\sigma(n)}=u_{1}'\ot\cdots\ot u_{m}'\stackrel{f_{1}\ot\cdots\ot f_{m}}{\rightarrow} v_{1}\ot\cdots\ot v_{m}.
\]
Where $\tau_{\sigma,u_{1},...,u_{n}}:u_{1}\ot\cdots \ot u_{n}\rightarrow u_{\sigma(1)}\ot\cdots \ot u_{\sigma(n)}$ is the coherence isomorphism given by $\gamma$ in $\C^{op}$.

We claim that this functor is a strict map of permutative categories. Note that by definition it follows that $\W(\u\odot\v)=\W(\u)\ot\W(\v)$. On the other hand, given $\u$ and $\v$ object in $\A$, the symmetry isomorphism in $\A$ is given by $\gamma=(\xi_{n,m},\underline{\text{id}}):\u\odot \v\to \v\odot \u$, where $\xi_{n,m}\in \Sigma_{n,m}$ is the permutation of $(n+m)$-letters that interchanges the first $n$-block with the last $j$-block. By definition and the coherence of the  we see that $\W((\xi_{n,m},\underline{\text{id}})):u_{1}\ot\cdots\ot u_{n}\ot v_{1}\ot\cdots\ot v_{m}\to v_{1}\ot\cdots\ot v_{m}\ot u_{1}\ot\dots\ot u_{n}$ equals the isomorphism given by $\gamma$ in $\C^{op}$. Thus we see that $\W$ is a strict map. In particular, we can see $\A$ and $\C^{op}$ as multicategories and a lax map between them is precisely a multifunctor. Thus $\W$ is an object in $\Mi(\A,\C^{op})$. If we take $\M=\A$, $\M'=\C^{op}$, $f=\W$ and $\E$ the trivial category in Theorem \ref{QuillenAd} (which in this case agrees with \cite[Corollary 12.3]{Elmendorf}) then we see that the functor 
\[
\W^{*}:\Mi(\C^{op},\Sym)\to \Mi(\A,\Sym)
\] 
has a left adjoint 
\[
\W_{\#}:\Mi(\A,\Sym)\to\Mi(\C^{op},\Sym)
\] 
such that the pair $(\W_{\#},\W^{*})$ forms a Quillen adjunction. Thus the image of the lax map $\vartheta:\A\to \Sym$ which we call $\phi:\C^{op}\to \Sym$ is a lax map.
\qed

\section{Group completion}

In this section we show that a lax map $\phi:\C^{op}\to \Sym$ gives rise to lax map $\phi':(\C^{op})^{-1}\C^{op}\to \Sym$. (Here $(\C^{op})^{-1}\C^{op}$ denotes the Grayson-Quillen group completion of $\C^{op}$). In particular the image of the identity is a strictly commutative ring spectrum which is the ring spectrum we are looking for. We begin by recalling the definition of $(\C^{op})^{-1}\C^{op}$.

\begin{definition} 
Let $(\D,\op,0)$ be a small symmetric monoidal category. Then the group completion of $\D$ is the category $\D^{-1}\D$, whose objects are the pairs $(a,b)$ , where $a,b$ are object of $\D$. The pair $(a,b)$ is thought of as the formal difference $b-a$. If $(a,b), (c,d)$ are two objects in $\D^{-1}\D$, then a morphism in $\D^{-1}\D$, $f:(a,b)\to (c,d)$ is an equivalence class of data of the form $(s,\alpha, \beta)$, where $s$ is an object in $\D$ and $\alpha:a\op s\to c$, $\beta:b\op s\to d$ are morphisms in $D$. Two pairs $(s,\alpha, \beta)$, $(s',\alpha', \beta')$ are equivalent if there exists a morphism $\gamma:s\to s'$, such that the following diagrams are commutative
\[
\xymatrix{
a\op s\ar[r]^-{\alpha}\ar[d]_-{1\op \gamma}   &c,& &b\op s\ar[r]^-{\beta}\ar[d]_-{1\op \gamma}   &d\\
a\op s'\ar[ru]_-{\alpha'}         & & &b\op s'\ar[ru]_-{\beta'}         &
}
\]
\end{definition}

If in $\D$ every morphism is an isomorphism and for every object $a$ of $\D$, the functor $a\op-:\D\to \D$ is faithful then the functor $i:\D\to \D^{-1}\D$ that assigns to every object $a$ the object $(0,a)$, gives rise to a group completion on the level of classifying spaces. See \cite{Thomason} and \cite{Grayson} for a complete treatment. The following lemma follows by a straight forward computation.

\begin{lemma}{\label{Permutative}}
Let $(\D,\op,0)$ be a small permutative category for which every morphism is an isomorphism. Then  $\D^{-1}\D$ is also a permutative category and the map $i:\D\to \D^{-1}\D$ is a (strong) lax map.
\end{lemma}

As a particular case of the previous lemma, we have that the map $i:\C^{op}\to (\C^{op})^{-1}\C^{op}$ is a lax map for a permutative category $\C$. Since $\C^{op}$, $(\C^{op})^{-1}\C^{op}$ are permutative categories, we can see them as multicategories and a lax map between them is precisely a multifunctor. Thus $i$ is an object in the category $\Mi(\C^{op},(\C^{op})^{-1}\C^{op})$. If we take $\M=\C^{op}$, $\M'=(\C^{op})^{-1}\C^{op}$, $f=i$ and $\C^{op}$ in Theorem \ref{QuillenAd} (which in this case agrees with \cite[Corollary 12.3]{Elmendorf}) then we see that the functor 
\[
i^{*}:\Mi((\C^{op})^{-1}\C^{op},\Sym)\to \Mi(\C^{op},\Sym)
\] 
has a left adjoint 
\[
i_{\#}:\Mi(\C^{op},\Sym)\to\Mi((\C^{op})^{-1}\C^{op},\Sym)
\] 
such that the pair $(i_{\#},i^{*})$ forms a Quillen adjunction. We have proved the following theorem

\begin{theorem}\label{assignment}
A discrete fibered symmetric bimonoidal category $\Lambda:\D\to \C$ gives rise to a lax map  $\phi':(\C^{op})^{-1}\C^{op}\to \Sym$. In particular we can associate to $\Lambda:\D\to \C$ the image of 1 under $\phi'$ which a strictly commutative symmetric ring spectrum.
\end{theorem}

\section{Functoriality}

In this section we show that the assignment of Theorem \ref{assignment} is functorial.

\begin{definition}
Suppose that $\Lambda:\D\to \C$ and $\Lambda':\D'\to \C$ are two discrete fibered symmetric bimonoidal categories. A morphism from $\Lambda:\D\to \C$ to $\Lambda':\D'\to \C$ is a functor $\Theta:\D\to \D'$ such that 
\[
\xymatrix{
\D\ar[r]^{\Theta}\ar[d]_{\Lambda}   &\D'\ar[dl]^-{\Lambda'}\\
\C         &
}
\]
is a commutative diagram. In addition, we require that $\Theta$ is a lax map with respect to the operations $\ot$ and for each object $c$ of $\C$, the restriction $\Theta_{|\D_{c}}:\D_{c}\to \D'_{\Theta(c)}$ is a lax map with respect to the operations $\op$ and $\op_{c}$.
\end{definition}

\begin{definition}
We denote by $\FS$ the category of discrete fibered symmetric bimonoidal categories and morphisms between them.
\end{definition}

We want to show that the assignment of Theorem \ref{assignment} is functorial. To begin, suppose that $\Lambda:\D\to \C$ and $\Lambda:\D'\to \C$ are two objects in $\FS$ and $\Theta:\D\to \D'$ a morphism between them. Then after the streefication process, the functor $\Theta$ induces a strict map $\Theta^{s}:\D^{s}\to \D'^{s}$. The  next step of the construction is to correspond functors $\Psi$ and $\Psi':\A\to \P$ to $\Lambda^{s}:\D^{s}\to \C^{s}$ and $\Lambda^{s}:\D'^{s}\to \C^{s}$ respectively. Here $\A$ is a certain wreath product category as defined on section 3 and $\P$ is the multicategory of small permutative category. Using the map $\Theta^{s}$, we can obtain a multifunctor $\Xi:\P\to \P$ by defining it to be trivial for those permutative categories not in the image of $\Psi$. The multifunctor $\Xi$ is such that the following diagram is commutative
\[
\xymatrix{
\A\ar[r]^-{\Psi}\ar[rd]_-{\Psi'}& \P\ar[d]^{\Xi}\\
 & \P.
}
\]
 
After this we apply the multifunctor $K$ and obtain a commutative diagram of multifunctors
\[
\xymatrix{
E\Sigma_{*}\ar[r]^-{T_{\D}}\ar[rd]_-{T_{\D'}}& \Sym^{\A}\ar[d]\\
 & \Sym^{\A}.
}
\]
The rectification performed on section 6 can be done in a functorial way and thus after passing to the completion we obtain a commutative diagram of lax maps
\[
\xymatrix{
\C^{op}\C\ar[r]^-{\phi'_{\D}}\ar[rd]_-{\phi'_{\D'}}& \Sym\ar[d]\\
 & \Sym.
}
\]
In particular, by looking at the image of the unit, we obtain a map of symmetric spectra
\begin{equation}\label{morphismfunc}
\phi'_{\D}(1)\to \phi'_{\D'}(1).
\end{equation}
Thus we obtain the following theorem.

\begin{theorem}
There is a functor
\[
\Zr:\FS\to \CRS
\]
that for a discrete fibered symmetric bimonoidal category $\Lambda:\D\to \C$ corresponds $\Zr(\Lambda)$, the image of the unit under the lax map
\[
\phi':\C^{op}\C\to \Sym.
\]
If $\Theta:\D\to \D'$ is a morphisms of discrete fibered categories, then $\Zr(\Theta)$ is the morphism of symmetric spectra
\[
\Zr(\Lambda)\to \Zr(\Lambda')
\]
as in (\ref{morphismfunc}).
\end{theorem}

\section{The definition}

In this section we define the concept of  topological fibered symmetric bimonoidal categories and show that to such each object we can correspond an $E_{\infty}$-ring spectrum in a functorial way. This spectrum is our ultimate goal.

\begin{definition}
A topological fibered category is a fibered category $\Lambda:\D\to \C$, where $\D$ and $\C$ are topological categories and the functor $\Lambda$ is a continuous functor that satisfies the same properties than a discrete fibered category in such a way that all the functors insight are continuous.
\end{definition}

In a similar way as in the discrete case we can form the category $\mathcal{TFS}$ of topological fibered symmetric bimonoidal categories.

\begin{remark} Whenever we talk about a topological category we mean a small category $\C$ whose object set is discrete and its morphism set is topological space in such a way that the structural maps are continuous.
\end{remark}

Suppose now that $\Lambda:\D\to \C$ is a topological fibered symmetric bimonoidal category. By applying the singular functor from topological spaces to simplicial sets on the level of objects and morphism, we can see a topological fibered symmetric bimonoidal category as a functor
\[
\mathfrak{D}:\Delta^{op}\to \FS;
\]
that is, we can see a topological fibered symmetric bimonoidal category as a simplicial object in the category $\FS$. By composing $\mathfrak{D}$ with $\Zr$, we can correspond to a topological fibered symmetric bimonoidal category a functor
\[
\Delta^{op}\to \CRS,
\]
where $\CRS$ is the category of strictly commutative symmetric ring spectrum. By realizing this simplicial object we obtain an $E_{\infty}$-ring spectrum. Since this construction is functorial we finally arrive to the following theorem.

\begin{theorem} 
There is a functor
\[
\mathfrak{Z}:\mathcal{TFS}\to \mathcal{E}_{\infty}
\]
where $\mathcal{E}_{\infty}$ is the category of $E_{\infty}$-ring spectra, that for a topological fibered symmetric bimonoidal category corresponds the realization of the simplicial object in the category $\CRS$ as mentioned above.
\end{theorem}

\end{document}